\documentclass[english, 3p]{elsarticle}

\pdfoutput=1 %For building with PDFLaTeX on arXiv

\usepackage[latin1]{inputenc}
\usepackage[T1]{fontenc}
\usepackage{babel,geometry,lmodern}
\usepackage{amsmath}
\usepackage{amsfonts}
\usepackage{amssymb}
\usepackage{graphicx}
\usepackage[caption=false]{subfig}
\usepackage{color}
\usepackage[hang,flushmargin]{footmisc} 
\usepackage{tikz}
\usepackage{multicol}
\usepackage{todonotes}
\usepackage{multirow}
\usepackage{framed}
\usepackage{stmaryrd} % for jump operator
\usepackage{enumitem}
\usepackage{bm}
\usepackage[ruled]{algorithm2e}
\usepackage{tablefootnote}
\usepackage{datatool}

\usepackage{amsthm}

\newtheorem{remark}{Remark}

\makeatletter
\def\ps@pprintTitle{%
	\let\@oddhead\@empty
	\let\@evenhead\@empty
	\def\@oddfoot{}%
	\let\@evenfoot\@oddfoot}
\makeatother

\numberwithin{remark}{section}
\numberwithin{equation}{section}
\newcommand{\lp}{\left(}
\newcommand{\rp}{\right)}
\newcommand{\Fh}[1]{%
	\ifthenelse{\isempty{#1}}{\mathcal{F}_h}{\mathcal{F}_{h,#1}}%
}
\newcommand{\bfn}{\mathbf{n}}
\newcommand{\bft}{\mathbf{t}}
\newcommand{\bfu}{\mathbf{u}}
\newcommand{\bfv}{\mathbf{v}}
\newcommand{\bfx}{\mathbf{x}}

\newcommand{\bfK}{\bm{\kappa}}
\newcommand{\bfI}{\mathbf{I}}
\newcommand{\pO}{\partial\Omega}
\newcommand{\Kh}{\mathcal{K}_h}
\newcommand{\pK}{\partial K}
\newcommand{\jump}[1]{\llbracket #1 \rrbracket}
\newcommand{\average}[1]{\{\hspace{-3pt}\{ #1 \}\hspace{-3pt}\}}
\newcommand{\ip}[3]{\lp#1,#2\rp_{#3}}
\newcommand{\diff}{\textup{d}}
\newcommand{\dt}{\diff t}

\newcommand{\dy}{\diff y}

\newcommand{\dd}[2]{\frac{\partial #1}{\partial #2}}
\newcommand{\tin}{\textup{in}}
\newcommand{\tout}{\textup{out}}
\newcommand{\tD}{{\textsc{d}}}
\newcommand{\tN}{{\textsc{n}}}
\newcommand{\tB}{{\textsc{b}}}
\newcommand{\tI}{{\textsc{i}}}
\newcommand{\tF}{{\textsc{f}}}
\newcommand{\tM}{{\textsc{m}}}
\newcommand{\tG}{{\textsc{g}}}
\newcommand{\ndof}{N_{\textup{dof}}}
\newcommand{\qmin}{\check{q}}
\newcommand{\qmax}{\hat{q}}
\newcommand{\qoi}{\textup{QOI}}
\newcommand{\powten}[2]{#1\cdot10^{#2}}

\begin{document}

\title{A simple embedded discrete fracture-matrix model for a coupled flow and transport problem in porous media}

\author[ntnu]{Lars H.~Ods\ae ter\corref{cor}}
\ead{lars.odsater@gmail.com}

\author[ntnu]{Trond Kvamsdal}
\author[mgl]{Mats G.~Larson}

\address[ntnu]{Department of Mathematical Sciences, NTNU Norwegian University of Science and Technology,\\ Alfred Getz' vei 1, 7491 Trondheim, Norway}
\address[mgl]{Department of Mathematics and Mathematical Statistics, Umeå University, SE-901 87 Umeå, Sweden}

\cortext[cor]{Corresponding author}

\begin{abstract}
	Accurate simulation of fluid flow and transport in fractured porous media is a key challenge in subsurface reservoir engineering. Due to the high ratio between its length and width, fractures can be modeled as lower dimensional interfaces embedded in the porous rock. We apply a recently developed embedded finite element method (EFEM) for the Darcy problem. This method allows for general fracture geometry, and the fractures may cut the finite element mesh arbitrarily. We present here a velocity model for EFEM and couple the Darcy problem to a transport problem for a passive solute. The main novelties of this work is a locally conservative velocity approximation derived from the EFEM solution, and the development of a lowest order upwind finite volume method for the transport problem. This numerical model is compatible with EFEM in the sense that the same computational mesh may be applied, so that we retain the same flexibility with respect to fracture geometry and meshing. Hence, our coupled solution strategy represents a simple approach in terms of formulation, implementation and meshing. We demonstrate our model by some numerical examples on both synthetic and realistic problems, including a benchmark study for single-phase flow. Despite the simplicity of the method, the results are promising.
\end{abstract}


\maketitle

\section{Introduction}
Modeling fluid flow in fractured porous media is an important yet challenging problem in subsurface engineering. Fractures are characterized as thin layers with either very high or low conductivity. They can therefore act as preferential paths or barriers and be essential to the fluid flow. The fracture width is typically several orders of magnitude smaller than any other characteristic sizes in the reservoir, and the flow rate can be orders of magnitude larger or smaller than in the surrounding matrix. This pose great challenges to the fracture model and the numerical method.

In this work we consider a discrete fracture-matrix (DFM) model, where the fractures are modeled as lower dimensional interfaces embedded in the rock matrix. We assume Darcy flow both in the matrix and the fracture, and we only consider the case where the permeability in the fractures are orders of magnitude larger than in the matrix. As flow model we consider incompressible single-phase flow governed by conservation of mass and Darcy's law. The flow problem is then coupled to a transport problem for a passive solute.

A common strategy to represent fractures in a DFM model is by averaging the governing equations across the fractures. The fracture width is then modeled as a coefficient in the equations rather than a geometrical property, and suitable coupling conditions between the fracture and matrix equations are applied. Common terminologies for such approximations are mixed-dimensional, hybrid-dimensional or reduced models. A method for the high permeability case was presented in \cite{alboin2000ddf,alboin2002mfa}, where the flow equations on mixed form are averaged over the interface assuming a continuous pressure. This model was later generalized by \cite{faille2002anf} and \cite{martin2005mfa} to also handle the low permeability case, where a Robin type condition on the pressure at the interface is enforced. More recent DFM models similar to \cite{faille2002anf,martin2005mfa} can be found in, e.g., \cite{angot2009aan,frih2012mfa,formaggia2014arm,schwenck2015drf,faille2016mra,boon2016rdo}.

To account for more general fracture shapes, the model in \cite{alboin2000ddf,alboin2002mfa} was extended to curved fractures in \cite{morales2010tnf,capatina2016nef}. These models are derived by considering the asymptotic limit in the weak formulation as the fracture width tends to zero and the fracture permeability tends to infinity. In this limit, the condition of a continuous pressure follows. We mention that this model was extended to allow for pressure jumps across the interface in \cite{capatina2016nef}.

An embedded finite element method (EFEM) for the model in \cite{capatina2016nef} was derived in \cite{burman2017asf}. Fractures are allowed to cut through the elements arbitrarily, and the fracture solution is represented by the restriction of the basis functions for the higher dimensional matrix elements. Contributions from the fracture to the fluid flow is included by superposition. A great advantage of this method is that it handles very general fracture geometry, including curved interfaces, bifurcations and intersections, and it is also easy to implement. The assumption of a continuous pressure along the fracture interface is a key ingredient as it allows for continuous elements. However, the normal flux is discontinuous across the fracture interface leading to loss of regularity, which can be resolved by local refinement close to the fractures based on an a priori error estimate. In this work, we apply EFEM for the pressure problem.

For the family of lower dimensional DFM models where the pressure is not assumed to be continuous, a great variety of numerical methods are suggested. We mention for instance finite element methods \cite{fard2003nso,martin2005mfa,angot2003amo}, finite volume methods \cite{angot2009aan,brenner2016gdo,karimi2016agg}, discontinuous Galerkin methods \cite{antonietti2016dga}, mimetic finite difference methods \cite{antonietti2016mfd,scotti2017aoa} and virtual element methods \cite{fumagalli2017dve}, all of which require a conforming mesh across the fracture interface. Different schemes with mortar coupling that allows for non-conforming meshes are presented in, e.g., \cite{frih2012mfa,boon2016rdo}, but still the mesh needs to explicitly represent the fractures. Fully non-conforming discretizations have been realized through extended finite elements \cite{dangelo2012amf,huang2011otu,schwenck2015drf,capatina2016nef}.

The embedded DFM (EDFM) introduced by \cite{li2008efs} is another approach where the fractures are allowed to cut arbitrarily through the higher dimensional mesh, see also \cite{hajibeygi2011ahf,moinfar2014doa,panfili2013ued,panfili2014som,filho2015ioa} for newer contributions. EDFM is based on the two-point flux approximation (TPFA) where the fracture-fracture and fracture-matrix transmissibilities are approximated from geometrical quantities. Similar to EFEM, they are only valid for high permeability fractures, but the recently introduced projection-based EDFM (pEDFM) \cite{tene2017pbe,jiang2017aip} are also able to handle low permeability fractures. An important difference to EFEM is that the fractures are represented by their own lower dimensional elements along the fracture interface.

A comprehensive comparison of several DFM models for single-phase flow was conducted in a recent benchmark study \cite{flemisch2018bfs}. We follow up on this study and evaluate EFEM on a selection of these benchmark problems.

The transport problem is also modeled by a reduced (or mixed-dimensional) model, see e.g., \cite{alboin2002mfa,fumagalli2013arm,fumagalli2017dve}, where a coupling term models flow between the fracture and the matrix. We assume advection dominated flow since we are primarily interested in the coupling with the flow problem. In this work we present a novel discretization that is compatible with EFEM in the sense that the same computational mesh can be applied. This ensures that we have the same flexibility in terms of fracture geometry and meshing for the coupled solution strategy as is the case for EFEM.
More specifically, we apply a zeroth order upwind finite volume method (FVM), where the fracture solution is represented by elements cut by the fracture, and where the coupling term is approximated in a non-standard way by evaluating the normal velocity (flux) on the boundary of such elements. We mention here that an alternative approach is to use the CutFEM technology, see \cite{burman2018cfe} for a stationary convection problem and \cite{burman2015cdg} for a general introduction to CutFEM.
For compatibility of the numerical solvers, we must require the velocity approximation from the flow problem to be locally conservative. This is not directly obtained by EFEM, but is achieved through a postprocessing step \cite{odsaeter2016pon}.

We mention that combining FEM and FVM for heterogeneous and fractured porous media have been studied before, e.g., the finite element--finite volume method \cite{geiger2004cfe,paluszny2007hfe}, where FEM is used for the flow problem and a node-centered finite volume method is used for transport. This method was extended in \cite{nick2011ahf} to allow for a discontinuous solution across interfaces for the transport problem. However, this approach requires fractures to me aligned with element boundaries and uses a dual mesh for the transport solver.

This paper is organized as follows. In Section \ref{sec:model} we describe the discrete fracture-matrix models and the governing equations for the coupled flow and transport problem. Next, in Section \ref{sec:methods}, we define the numerical methods. This includes EFEM for the flow problem, FVM for the transport problem, and the velocity approximation which couples the two subproblems. In Section \ref{sec:results} some numerical results are presented, including a realistic problem with a complex fracture network. Finally, we make some concluding remarks in Section \ref{sec:conclusion}.

\section{Model formulation}
\label{sec:model}
As model problem we consider incompressible single-phase flow with advective transport of a concentration in a fractured porous media. Let $\Omega\in\mathbb{R}^d$, with $d=2,3$, be a convex polygonal domain with an embedded interface $\Gamma$ representing the fractures. The bulk domain $\Omega\setminus\Gamma$ will be referred to as the matrix.  We denote by $\bfK$ and $\phi$ the symmetric positive definite permeability tensor and porosity of the matrix, respectively. The fracture permeability is assumed to be isotropic and is denoted $\kappa_\Gamma$. We denote by $w$ and $\phi_\Gamma$ the fracture aperture and the porosity of the fractures, respectively. Moreover, $k_\Gamma=w\kappa_\Gamma$ is the effective (scaled by fracture aperture) fracture permeability. Next, $q$ and $q_\Gamma$ denotes source or sink terms in the matrix and fractures, respectively. The primary unknowns are the fluid pressure $p$, from which we can derive the fluid velocity $\bfu$, and the concentration $c$.

For the flow problem, the boundary $\pO$ is partitioned into a Dirichlet and Neumann part, denoted $\pO_\tD$ and $\pO_\tN$, respectively. For the transport problem, we let 
$\pO_\tin := \{\bfx\in\pO \vert \bfu\cdot\bfn < 0\}$ denote the inflow boundary, and $\pO_\tout := \{\bfx\in\pO \vert \bfu\cdot\bfn \ge 0\}$ denote the outflow boundary, where $\bfn$ is the outward pointing unit normal.

We use the standard notation $H^s(\omega)$ for the Sobolev space of order $s$ on $\omega$ with the special cases $L^2(\omega) = H^0(\omega)$ and $H_0^1(\Omega) = \{ v\in H^1(\Omega) : v\vert_{\pO_D}=0 \}$. The space of continuous functions on $\omega$ is denoted $C(\omega)$. For a normed vector space $V$, we let $\Vert\cdot\Vert_V$ denote the norm on $V$. For $V=L^2(\omega)$, we use the notation $\Vert\cdot\Vert_{L^2(\omega)}=\Vert\cdot\Vert_\omega$, and denote by $(\cdot,\cdot)_\omega$ the $L^2$ scalar product.

In the following, we restrict this presentation to the two dimensional case, i.e., $d=2$, but most of the theory and methods considered herein can by extended to three dimensions in a straight forward manner.

\subsection{Fracture representation}

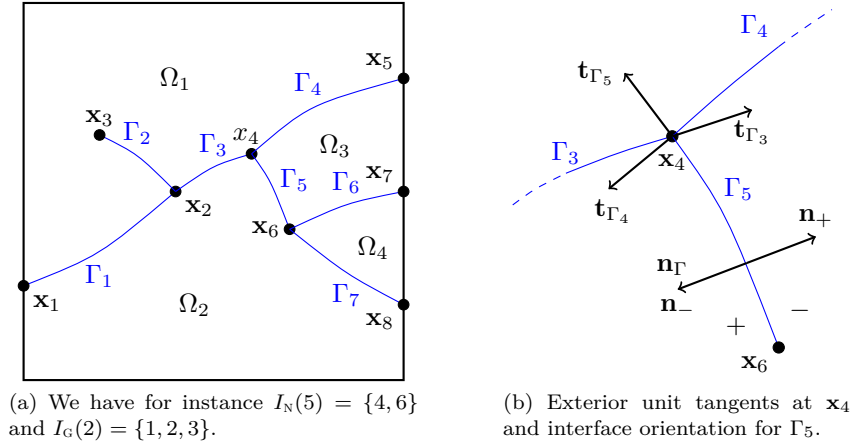
\begin{figure}[tbp]
	\centering
	\subfloat[We have for instance $I_\tN(5)=\{4,6\}$ and $I_\tG(2)=\{1,2,3\}$.]{
		\begin{tikzpicture}[scale=5]
		\draw[thick] (0,0) -- (1,0) -- (1,1) -- (0,1) -- cycle;
		
		\draw[blue] (0,0.25) .. controls (0.2,0.33) .. (0.4,0.5);
		\filldraw[black] (0,0.25) circle[radius=0.4pt] node[below right]{$\bfx_1$};
		\filldraw[black] (0.4,0.5) circle[radius=0.4pt] node[below right]{$\bfx_2$};
		
		\draw[blue] (0.2,0.65) .. controls (0.3,0.6) .. (0.4,0.5);
		\filldraw[black] (0.2,0.65) circle[radius=0.4pt] node[above]{$\bfx_3$};
		
		\draw[blue] (0.4,0.5) .. controls (0.5,0.57) .. (0.6,0.6);
		\filldraw[black] (0.6,0.6) circle[radius=0.4pt] node[above]{\bf$x_4\ $};
		
		\draw[blue] (0.6,0.6) .. controls (0.75,0.73) .. (1.0,0.8);
		\filldraw[black] (1.0,0.8) circle[radius=0.4pt] node[above left]{$\bfx_5$};
		
		\draw[blue] (0.6,0.6) .. controls (0.65,0.53) .. (0.7,0.4);
		\filldraw[black] (0.7,0.4) circle[radius=0.4pt] node[left]{$\bfx_6$};
		
		\draw[blue] (0.7,0.4) .. controls (0.85,0.47) .. (1.0,0.5);
		\filldraw[black] (1.0,0.5) circle[radius=0.4pt] node[above left]{$\bfx_7$};
		
		\draw[blue] (0.7,0.4) .. controls (0.85,0.28) .. (1.0,0.2);
		\filldraw[black] (1.0,0.2) circle[radius=0.4pt] node[below left]{$\bfx_8$};
		
		\node at (0.4,0.8) {$\Omega_1$};
		\node at (0.45,0.2) {$\Omega_2$};
		\node at (0.82,0.62) {$\Omega_3$};
		\node at (0.92,0.35) {$\Omega_4$};
		
		\node[blue] at (0.2,0.33) [below]{$\Gamma_1$};
		\node[blue] at (0.3,0.6) [above]{$\Gamma_2$};
		\node[blue] at (0.5,0.57) [above]{$\Gamma_3$};
		\node[blue] at (0.75,0.73) [above]{$\Gamma_4$};
		\node[blue] at (0.65,0.53) [right]{$\Gamma_5$};
		\node[blue] at (0.85,0.47) [above]{$\Gamma_6$};
		\node[blue] at (0.85,0.28) [below]{$\Gamma_7$};
		
		\end{tikzpicture}
		\label{fig:frac_prob_ill}
	}
	\hspace{10mm}
	\subfloat[Exterior unit tangents at $\bfx_4$ and interface orientation for $\Gamma_5$.]{
		\begin{tikzpicture}[scale=14]
		\draw[blue] (0.5,0.565) .. controls (0.55,0.585) .. (0.6,0.6);
		\draw[blue,dashed] (0.45,0.535) .. controls (0.47,0.55) .. (0.5,0.565) node[above]{$\Gamma_3$};
		\filldraw[black] (0.6,0.6) circle[radius=0.15pt] node[yshift=-3.5mm]{$\bfx_4$};
		\draw[blue] (0.6,0.6) .. controls (0.65,0.645) .. (0.7,0.685) node[above left]{$\Gamma_4$};
		\draw[blue,dashed] (0.7,0.685) .. controls (0.72,0.7) .. (0.75,0.72);
		\draw[blue] (0.6,0.6) .. controls (0.65,0.53) .. (0.7,0.4);
		\node[blue] at (0.64,0.55) [right]{$\Gamma_5$};
		\filldraw[black] (0.7,0.4) circle[radius=0.15pt] node[below left]{$\bfx_6$};
		
		\draw[thick,->] (0.6,0.6) -- (0.675,0.625) node[below]{$\bft_{\Gamma_3}$};
		\draw[thick,->] (0.6,0.6) -- (0.54,0.55) node[below]{$\ \bft_{\Gamma_4}$};
		\draw[thick,->] (0.6,0.6) -- (0.555,0.66) node[left]{$\ \bft_{\Gamma_5}$};
		
		\draw[thick,->] (0.67,0.48) -- (0.735,0.505) node[above]{$\bfn_+$};
		\draw[thick,->] (0.67,0.48) -- (0.605,0.455) node[below]{$\bfn_-$};
		\draw (0.60,0.46) node[above]{$\bfn_\Gamma$};
		
		\draw (0.66,0.42) node{$+$};
		\draw (0.72,0.435) node{$-$};
		
		\end{tikzpicture}
		\label{fig:gamma_directions}
	}
	\caption{Illustration of domain and fracture representation.}
\end{figure}

We allow for bifurcating fractures and represent $\Gamma$ as a graph with nodes $\mathcal{N}= \{\bfx_i\}_{i\in I_\tN}$ and edges $\mathcal{G}=\{\Gamma_j\}_{j\in I_\tG}$, where $I_\tN$ and $I_\tG$ are finite index sets, and each $\Gamma_j$ is a curve between two nodes with indices $I_\tN(j)$. For each $i\in I_\tN$, we let $I_\tG(i)$ be the set of indices corresponding to curves for which $\bfx_i$ is an end point. Furthermore, let $\{\Omega_i\}_{i=1}^{n_d}$ be a partition of $\Omega$ into $n_d$ subdomains defined by $\Gamma$. See Fig.~\ref{fig:frac_prob_ill}.

We fix an orientation of each $\Gamma_j$ such that the positive direction is from the node with lowest index towards the node with highest index. Then we define $\bfn_\Gamma$ to be the unit normal on $\Gamma$ pointing from the left side towards the right when facing the positive direction of $\Gamma$, see Fig~\ref{fig:gamma_directions}. For a scalar function $v$, possibly discontinuous at $\Gamma$, we define the jump as
\begin{align}
\jump{v} = v_+ - v_-, \quad \text{on } \Gamma,
\label{eq:jump_scalar}
\end{align}
where
\begin{align}
v_+(\bfx) &= \lim_{\epsilon\rightarrow 0^+} v(\bfx + \epsilon\bfn_\Gamma), \quad \bfx\in\Gamma, \\
v_-(\bfx) &= \lim_{\epsilon\rightarrow 0^+} v(\bfx - \epsilon\bfn_\Gamma), \quad \bfx\in\Gamma.
\end{align}
For a vector valued function $\bfv$, we define the jump in the normal component across $\Gamma$ as
\begin{align}
\jump{\bfv\cdot\bfn} = \bfv_+\cdot\bfn_+ + \bfv_-\cdot\bfn_-,
\label{eq:jump_vector}
\end{align}
where $\bfn_+ = -\bfn_\Gamma$ and $\bfn_- = \bfn_\Gamma$.

\subsection{Pressure problem}

The fractures are modeled as embedded surfaces with high permeability. Our model is essentially the same as the one presented in \cite{alboin2002mfa,capatina2016nef}. It was further studied in \cite{burman2017asf}, where it was extended to bifurcating cracks. The embedded model is derived as the asymptotic limit in the weak formulation as the fracture aperture goes to zero and the fracture permeability goes to infinity. We refer to \cite{alboin2002mfa,capatina2016nef} for further details.

\paragraph{Boundary value problem}

The governing equations for the pressure problem are
\begin{subequations}
\begin{align}
-\nabla\cdot \lp \bfK\nabla p \rp &= q, &\text{in } \Omega, \label{eq:darcy_omega} \\
-\nabla_{\Gamma}\cdot(k_{\Gamma}\nabla_{\Gamma}p) &= q_{\Gamma} + \jump{(\bfK\nabla p)\cdot\bfn}, 
&\text{on } \Gamma, \label{eq:darcy_gamma} \\
\jump{p} &=0, &\text{on } \Gamma, \label{eq:pressure_jump}
\end{align}
where $\nabla_{\Gamma}=\mathbf{P}\nabla$ is the tangential gradient with $\mathbf{P}=\mathbf{I}-\bfn_\Gamma\otimes\bfn_\Gamma$. 
The first equation \eqref{eq:darcy_omega} is the standard Darcy equation for single-phase flow describing conservation of mass. Eq.~\eqref{eq:darcy_gamma} governs conservation of mass in the fractures, where the last term on the right hand side represents a coupling term for the normal velocity across $\Gamma$. The last equation \eqref{eq:pressure_jump} is another coupling term, and represents a reasonable assumption for high conductive fractures. The model can also be extended to a non-zero pressure jump \cite{capatina2016nef}.

We equip the governing equations with Dirichlet and Neumann boundary conditions, i.e.,
\begin{align}
p &= p_\tD, &\text{on } \pO_\tD, \\
(\bfK\nabla p)\cdot\bfn &= u_\tN, &\text{on } \pO_\tN.
\end{align}
Furthermore, we enforce continuity of pressure at the interface nodes,
\begin{align}
p_{\Gamma_k}(x_i) &= p_{\Gamma_l}(x_i), &\forall k,l\in I_\tG(i), \forall i\in I_\tN,
\end{align}
and apply the Kirchhoff condition
\begin{align}
\sum_{j\in I_\tG(i)} ((k_{\Gamma_j}\nabla_{\Gamma_j}p_{\Gamma_j})\cdot\bft_{\Gamma_j})\vert_{x_i} &=0, &\forall i\in I_\tN\setminus I_{\tN,\tN}, \label{eq:kirchhoff} \\
((k_{\Gamma_j}\nabla_{\Gamma_j}p_{\Gamma_j})\cdot\bft_{\Gamma_j})\vert_{x_i} &= w u_\tN, & \forall i\in I_{\tN,\tN}, \label{eq:fracture_neumann}
\end{align}
\end{subequations}
where $\bft_{\Gamma_j}$ is the exterior unit tangent to $\Gamma_j$ (see Fig.~\ref{fig:gamma_directions}) and $I_{\tN,\tN}\subset I_\tN$ is the set of indices whose corresponding nodes belong to $\pO_\tN$. We remark that Eq.~\eqref{eq:kirchhoff} ensures mass balance at the interface nodes, while Eq.~\eqref{eq:fracture_neumann} is a Neumann condition for the part of $\Gamma$ that intersects with $\pO_N$, and where the fracture width $w$ is taken into account. Observe that Eq.~\eqref{eq:kirchhoff} implies a homogeneous Neumann condition at the fracture tip if the tip is in the interior of $\Omega$. This is a natural approximation that is commonly used, see e.g., \cite{angot2009aan,frih2012mfa,formaggia2014arm,boon2016rdo}.

The fluid velocity is defined by Darcy's law as $\bfu = -\bfK\nabla p$. We note that $\bfu_{\Gamma} = -k_{\Gamma}\nabla_{\Gamma}p$ gives the flow rate through the cross section of the fracture rather than the velocity.

\paragraph{Weak formulation}

We define the following function spaces,
\begin{align}
V_{\Gamma}(\Gamma) &= \{ v\in C(\Gamma) : v\in H^1(\Gamma_j),\,\forall j\in I_\tG \}, \\
V_0(\Omega) &= \{ v\in H_0^1(\Omega) : v\vert_\Gamma\in V_\Gamma(\Gamma) \}, \\
V_\tD(\Omega;v_\tD) &= \{ v\in H^1(\Omega) : v\vert_{\pO_\tD} = v_\tD,\, v\vert_\Gamma\in V_\Gamma(\Gamma) \}.
\end{align}
Then multiply \eqref{eq:darcy_omega} by a test function $v\in V_0(\Omega)$, integrate over $\Omega$ and apply Green's formula on each subdomain $\Omega_i$, to obtain
\begin{align}
\ip{q}{v}{\Omega} 
&= \sum_{i=1}^{n_d} \ip{-\nabla\cdot \lp \bfK\nabla p \rp}{v}{\Omega_i} \nonumber \\
&= \sum_{i=1}^{n_d} \lp \ip{\bfK\nabla p}{\nabla v}{\Omega_i} - \ip{(\bfK\nabla p)\cdot\bfn_i}{v}{\partial\Omega_i} \rp \nonumber \\
&= \ip{\bfK\nabla p}{\nabla v}{\Omega} - \ip{\jump{(\bfK\nabla p)\cdot\bfn}}{v}{\Gamma} - \ip{u_N}{v}{\partial\Omega_\tN} \nonumber \\
&= \ip{\bfK\nabla p}{\nabla v}{\Omega} - \ip{q_{\Gamma}}{v}{\Gamma} - \ip{\nabla_{\Gamma}\cdot(k_{\Gamma}\nabla_{\Gamma}p)}{v}{\Gamma} - \ip{u_\tN}{v}{\partial\Omega_\tN}.
\label{eq:pressure_weak_derivation}
\end{align}
Notice that $\bfn_i$ denotes the unit normal pointing out of $\Omega_i$. Then apply Green's formula on $\Gamma$ and use the Kirchoff conditions, Eq.~\eqref{eq:kirchhoff}--\eqref{eq:fracture_neumann}, to obtain
\begin{align}
-\ip{\nabla_{\Gamma}\cdot(k_{\Gamma}\nabla_{\Gamma}p)}{v}{\Gamma} 
&= \sum_{j\in I_\tG} -\ip{\nabla_{\Gamma_j}\cdot(k_{\Gamma_j}\nabla_{\Gamma_j}p)}{v}{\Gamma_j} \nonumber \\
&= \sum_{j\in I_\tG} \ip{k_{\Gamma_j}\nabla_{\Gamma_j}p}{\nabla_{\Gamma_j} v}{\Gamma_j}
- \sum_{j\in I_\tG} \sum_{i\in I_\tN(j)} \ip{(k_{\Gamma_j}\nabla_{\Gamma_j}p)\cdot\bft_{\Gamma_j}}{v}{x_i} \nonumber \\
&= \ip{k_{\Gamma}\nabla_{\Gamma}p}{\nabla v}{\Gamma} - \sum_{i\in I_\tN} \sum_{j\in I_\tG(i)} \lp (k_{\Gamma_j}\nabla_{\Gamma_j}p)\cdot\bft_{\Gamma_j} \rp\vert_{x_i} v(x_i) \nonumber \\
&= \ip{k_{\Gamma}\nabla_{\Gamma}p}{\nabla v}{\Gamma} - \sum_{i\in I_{\tN,\tN}} w u_\tN(x_i) v(x_i).
\end{align}
Combing these equations we get the following weak formulation. Find $p\in V_\tD(\Omega;p_\tD)$ such that
\begin{align}
a(p,v) = l(v), \qquad \forall v\in V_0(\Omega),
\label{eq:pressure_weak_form}
\end{align}
where
\begin{align}
a(u,v) &= \ip{\bfK\nabla u}{\nabla v}{\Omega} + \ip{k_{\Gamma}\nabla_{\Gamma}u}{\nabla v}{\Gamma}, \\
l(v)   &= \ip{q}{v}{\Omega} + \ip{q_{\Gamma}}{v}{\Gamma} 
- \ip{u_N}{v}{\partial\Omega_\tN} - \sum_{i\in I_{\tN,\tN}} w u_\tN(x_i) v(x_i).
\end{align}
Observe that the contribution from the fractures are included by superposition, i.e., by evaluating lower dimensional integrals along $\Gamma$.

\begin{remark}
When applying Green's formula in Eq.~\eqref{eq:pressure_weak_derivation}, we assume that there are no interfaces $\Gamma_j$ that terminates in the interior of $\Omega$ (as is the case for $\Gamma_2$ in Fig.~\ref{fig:frac_prob_ill}). In such cases we can divide in two the domains $\Omega_i$ with a terminating node such that the two new domains, denoted $\Omega_{i\textup{a}}$ and $\Omega_{i\textup{b}}$, are separated by the interface with the terminating node and an artifical line connecting the terminating node with either $\pO$ or another interface. We can then replace the contribution from $\Omega_i$ by the sum of the contributions from $\Omega_{i\textup{a}}$ and $\Omega_{i\textup{b}}$ in Eq.~\eqref{eq:pressure_weak_derivation}. For the situation in Fig.~\ref{fig:frac_prob_ill}, we may for instance divide $\Omega_1$ into two subdomains, $\Omega_{1\textup{a}}$ and $\Omega_{1\textup{b}}$, separated by $\Gamma_2$ and a line connecting $\bfx_3$ with the upper left corner of $\Omega$. We mention that a similar approach was used in \cite{angot2009aan}.
\end{remark}

\subsection{Transport problem}

Advective transport in the fractured domain is modeled by a reduced (or mixed-dimensional) model, see e.g., \cite{alboin2002mfa,fumagalli2013arm,fumagalli2017dve},
\begin{subequations}
\begin{align}
\phi \frac{\partial c}{\partial t} + \nabla\cdot(\bfu c) &= f(c), &\text{in } \Omega\times(0,T], 
\label{eq:transport_matrix} \\
w\phi_\Gamma \frac{\partial c_\Gamma}{\partial t} + \nabla_{\Gamma}\cdot(\bfu_\Gamma c_\Gamma) - \jump{\bfu\cdot\bfn c^*} &= f_\Gamma(c_\Gamma), &\text{in } \Gamma\times(0,T],
\label{eq:transport_fracture}
\end{align}
Initial and boundary conditions are given as
\begin{align}
c &= c_0, &\text{on } \Omega\times\{0\}, \\
c_\Gamma &= c_{\Gamma,0}, &\text{on } \Gamma\times\{0\}, \\
c &= c_\tB, &\text{on } \pO_\tin\times(0,T], \\
c_\Gamma &= c_{\Gamma,\tB}, &\text{on } \Gamma_\tin\times(0,T].
\end{align}
\label{eq:transport_problem}%
\end{subequations}
Here $c_0$ and $c_{\Gamma,0}$ are the initial concentrations in the matrix and fractures, respectively, while $c_\tB$ and $c_{\Gamma,\tB}$ are the inflow concentrations for the matrix and fractures, respectively.
The right hand sides denote source terms, defined as
\begin{subequations}
\begin{align}
f(c) &= \qmin c + \qmax c_\textup{w} = 
\begin{cases}
qc, & \text{if } q \le 0, 
\label{eq:rhs_def_matrix} \\
qc_\textup{w}, & \text{if } q > 0,
\end{cases} \\
f_\Gamma(c_\Gamma) &= \qmin_\Gamma c_\Gamma + \qmax_\Gamma c_\textup{w} = 
\begin{cases}
q_\Gamma c_\Gamma, & \text{if } q_\Gamma \le 0, \\
q_\Gamma c_\textup{w}, & \text{if } q_\Gamma > 0,
\end{cases}
\label{eq:rhs_def_fracture}
\end{align}
\label{eq:rhs_def}%
\end{subequations}
where $\qmin = \min(q,0)$ and $\qmax = \max(q,0)$, and $c_\textup{w}$ is the inflow concentration from the source term. The third term in the fracture equation \eqref{eq:transport_fracture}, $\jump{\bfu\cdot\bfn c^*}$, is a coupling term that models flow between the fracture and matrix. Here, $c^*$ is interpreted as
\begin{align}
c^*_\pm = 
\begin{cases}
c_\pm, &\text{if } (\bfu\cdot\bfn)_\pm \ge 0, \\
c_\Gamma, &\text{otherwise},
\end{cases}
\quad\text{on } \Gamma\times(0,T].
\label{eq:upwind_conc_gamma}
\end{align}

Observe that the transport problem is coupled to the pressure problem through the velocities $\bfu$ and $\bfu_\Gamma$.

\section{Numerical methods}
\label{sec:methods}
\subsection{Preliminaries}

\paragraph{Domain discretization}

Let $\Kh$ be a partition of $\Omega$, and denote by $K\in\Kh$ an element of the partition. We let $K\in\Kh$ be open such that $\bigcup_{K\in\Kh}\bar{K}=\bar{\Omega}$. The diameter of $K$ is denoted $h_K$, while $h$ is the maximum diameter of all elements. We assume $\Kh$ to be regular and quasi-uniform. By regular we mean that all elements are convex and that there exists $\rho>0$ such that each element $K\in\Kh$ contains a ball of radius $\rho h_K$ in its interior. Furthermore, $\Kh$ is quasi-uniform if there exists $\tau>0$ such that $h/h_K \le \tau$ for all $K\in\Kh$.

We denote by $\Fh{}$ the set of all element faces. This set is then divided into the set of interior faces, $\Fh{\tI}$, and boundary faces, $\Fh{\tB}$. We assume that each face in $\Fh{\tB}$ is either completely on the Dirichlet or Neumann part of the boundary, such that $\Fh{\tB}$ can be split into $\Fh{\tD}$ and $\Fh{\tN}$, i.e., the sets of faces on the Dirichlet and Neumann boundary, respectively. Similarly, let $\Fh{\tin}$ and $\Fh{\tout}$ be the set of boundary faces on the inflow and outflow boundary, respectively.
For each face $F\in\Fh{}$ we choose an orientation and let $\mathbf{n}_{F}$ denote the unit normal in this direction. The unit normal vector on $F\in\Fh{\tB}$ is chosen to coincide with the outward unit normal. Furthermore, $\bfn_K$ denotes the unit normal pointing out of $K$. We will also use the notation $\vert K\vert$ for the measure of $K$, and similarly $\vert F\vert$ for the measure of $F$.

For the transport solver, we further divide $\Kh$ into two subsets, $\Kh^\tM$ and $\Kh^\tF$, where $\Kh^\tF$ contains all fractured cells, $\Kh^\tF = \left\{ K\in\Kh : K\cap\Gamma\neq\emptyset \right\}$, and $\Kh^\tM = \Kh\setminus\Kh^\tF$, see Fig.~\ref{fig:fracture_elements}.
The set of interior faces, $\Fh{\tI}$ are then partitioned into three subsets, $\Fh{\tI}^\tF$, $\Fh{\tI}^\tM$ and $\Fh{\tI}^{\textsc{fm}}$, where $\Fh{\tI}^\tF$ are the set of faces between two fractured elements, $\Fh{\tI}^\tM$ is the set of faces between two matrix elements, and $\Fh{\tI}^{\tF\tM}$ are the sets of faces between a fracture and a matrix element.

\begin{remark}
	In the situation where $\Gamma$ coincides with an element face $F=\pK_-\cap \pK_+\in\Fh{\tI}$, we need to choose which of the neighboring elements that belong to $\Kh^\tF$. One possibility is to choose the element for which $\bfn_\Gamma$ is exterior, i.e., $K_-$. For the numerical examples presented in Section~\ref{sec:results} we have avoided this situation.
\end{remark}

\paragraph{Jump and average operators}

We define the jump operator $\jump{\cdot}$ over a face $F\in\Fh{}$ in the same way as we did for the jump over $\Gamma$, see Eq.~\eqref{eq:jump_scalar} and \eqref{eq:jump_vector}, where $\bfn_F$ now defines the orientation. 
Furthermore, we denote by $\average{\cdot}_\theta$ the weighted average operator on $F$, defined as
\begin{align}
\average{v}_\theta = \theta_Fv_- + (1-\theta_F)v_+,
\label{eq:average_def}
\end{align}
where $\theta_F=\theta\vert_F$ and $0\le\theta\le 1$. For $\theta=\tfrac{1}{2}$ we simply write $\average{\cdot}$ without any subscript. For $F\in\Fh{\tB}$ we define the jump and average operators as the one sided value, i.e.,
\begin{align}
\jump{v} = \average{v}_\theta = v_-.
\end{align}

\paragraph{Finite dimensional function spaces}

In our implementation we work with quadrilateral elements, but the all numerical methods can equally well be formulated on other elements, e.g., triangular. Denote by $\hat{K}=(0,1)^2$ the reference element with coordinates $(\xi,\eta)$, and by $M_K$ the mapping from $\hat{K}$ to $K$. With this, we denote by $\hat{Q}_r(\hat{K})$ the tensor product of polynomial spaces of degree less than or equal to $r$ in each spatial direction, i.e.,
\begin{align}
\textstyle
\hat{Q}_r(\hat{K}) = \left\{ v\in H^1(\hat{K}) : v = \lp \sum_{i=0}^r a_i\xi^i \rp \lp \sum_{i=0}^r b_i\eta^i \rp, a_i,b_i \in \mathbb{R} \right\}.
\end{align}
Next, $Q_r(K)$ denotes the reference element functions mapped to the actual element $K$,
\begin{align}
Q_r(K) = \{ \hat{v} \circ M_K^{-1} : \hat{v}\in \hat{Q}_r(\hat{K}) \}.
\end{align}
We may now define the following function spaces of piecewise polynomials of order $r>0$,
\begin{align}
Q_r(\Kh) &= \{ v\in C(\Omega) : v\vert_K \in Q_r(K), K\in \Kh \}, \\
Q_{r,\tD}(\Kh;v_\tD) &= \{ v\in Q_r(\Kh) : v\vert_{\pO_\tD} = v_\tD \}.
\end{align}
Moreover, we define the following spaces of piecewise constant functions ($r=0$),
\begin{align}
Q_0(\Kh) &= \{ v\in L^2(\Omega) : v\vert_K = a_K, a_K\in\mathbb{R}, K\in \Kh \}, \\
Q_0(\Fh{}) &= \{ v\in L^2(\Fh{}) : v\vert_F = a_F, a_F\in\mathbb{R}, F\in \Fh{} \}.
\end{align}
Finally, we denote by $\ip{\cdot}{\cdot}{\Kh}$ and $\ip{\cdot}{\cdot}{\Fh{}}$ the broken $L^2$ scalar products, i.e.,
\begin{align}
\ip{u}{v}{\Kh} = \sum_{K\in\Kh} (u,v)_K, \qquad \ip{u}{v}{\Fh{}} = \sum_{F\in\Fh{}} (u,v)_F.
\end{align}

\subsection{Pressure problem}

We follow \cite{burman2017asf} and approximate the pressure solution with piecewise bilinear functions by  restricting the weak formulation \eqref{eq:pressure_weak_form} to the finite dimensional subspace $Q_{1,\tD}(\Kh;p_\tD)\subset V_D(\Omega;p_\tD)$.
Find $p_h\in Q_{1,\tD}(\Kh;p_\tD)$ such that
\begin{align}
a(p_h,v) = l(v), \quad \forall v \in Q_{1,\tD}(\Kh;0).
\label{eq:efem}
\end{align}

The following a priori error estimate for the pressure approximation was proved in \cite{burman2017asf}. Let $\mathcal{N}_h(K)\subset\Kh$ be the set of all elements which are node neighbors of $K$, and let $h_\Gamma$ denote the mesh parameter in the vicinity of $\Gamma$ such that $h_K \le h_\Gamma$ for all $K\in \mathcal{N}_h(\Kh^\tF)$.
Then it holds that
\begin{align}
\Vert p-p_h\Vert_{\Omega} + \Vert p-p_h\Vert_{\Gamma} \lesssim (h_{\Gamma}+h^2) 
\lp \sum_{i=1}^{n_d} \Vert p\Vert_{H^2(\Omega_i)}\rp
+ h_{\Gamma}^2 \Vert p\Vert_{H^2(\Gamma)}.
\label{eq:apriori_l2_pressure}
\end{align}
As a consequence, one should refine locally around the fractures until $h_\Gamma\sim h^2$ to obtain the optimal order of convergence in terms of $h$.

\begin{remark}
The estimate \eqref{eq:apriori_l2_pressure} was proved in \cite{burman2017asf} for the simple geometry where $\Gamma$ is a smooth embedded interface in the interior of $\Omega$ without boundary. It can be extended to the case where $\Gamma$ is represented as a graph. However, depending on the geometry of $\Gamma$, we may loose some regularity of the solution (this is especially the case for terminating nodes in the interior of $\Omega$, e.g., $\bfx_2$ in Fig.~\ref{fig:frac_prob_ill}), so that we must replace the term $( h_{\Gamma}+h^2)$ by $(h_{\Gamma}+h^{s})$ for some $s\in[1,2]$. Yet, the condition $h_\Gamma\sim h^2$ is still sufficient to obtain the optimal convergence in terms of $h$ since refinement around $\Gamma$ also means refinement around the interface nodes.
\end{remark}

\subsection{Transport problem}
\label{sec:methods_transport}

We approximate the concentration solution by piecewise constants $c_h\in Q_0(\Kh)$, and let $c_h$ on $\Kh^\tF$ represent the concentration in the fractures, and $c_h$ on $\Kh^\tM$ represent the concentration in the matrix. We use an upwind approximation of the concentration on element faces. The numerical scheme can be formulated as a zeroth order finite volume method (FV), or equivalently as a zeroth order Discontinuous Galerkin method (DG). We only express the FV formulation here, and refer to Appendix \ref{app_DG} for the DG formulation.

\paragraph{FV formulation}

We integrate the matrix equation \eqref{eq:transport_matrix} over $K\in\Kh^\tM$ to obtain the integral formulation
\begin{align}
\int_K \phi\dd{c}{t} + \int_{\pK}\bfu\cdot\bfn_Kc = \int_K f(c), \quad K\in\Kh^\tM.
\label{eq:integral_eq_matrix}
\end{align}
Similarly, for $K\in\Kh^\tF$, we integrate the fracture equation \eqref{eq:transport_matrix} over $K\cap\Gamma$ to obtain
\begin{align}
\int_{K\cap\Gamma} w\phi_\Gamma\dd{c_\Gamma}{t} + \int_{\partial(K\cap\Gamma)} \bfu_\Gamma\cdot\bfn_{K\cap\Gamma} c_\Gamma - \int_{K\cap\Gamma} \jump{\bfu\cdot\bfn c^*} = \int_{K\cap\Gamma} f_\Gamma(c_\Gamma), \quad K\in\Kh^\tF.
\label{eq:integral_eq_fracture}
\end{align}

The lowest order finite volume method is then obtained by replacing $c$ by $c_h\in Q_0(\Kh)$. We use an upwind approximation on $\pK$, i.e.,
\begin{subequations}
\begin{align}
\bfu\cdot\bfn_K c_h\vert_{F=\pK\cap\partial\tilde{K}} = 
\begin{cases}
\bfu\cdot\bfn_K c_h\vert_K, & \text{if } \bfu\cdot\bfn_K \ge 0, \\
\bfu\cdot\bfn_K c_h\vert_{\tilde{K}}, & \text{if } \bfu\cdot\bfn_K < 0,
\end{cases}
\end{align}
where $\tilde{K}$ is a neighbor element of $K$. If $F\subset\pK$ is a boundary face, we have
\begin{align}
\bfu\cdot\bfn_K c_h\vert_{F=\pK\cap\pO} = 
\begin{cases}
\bfu\cdot\bfn_K c_h\vert_K, & \text{if } \bfu\cdot\bfn_K \ge 0, \\
\bfu\cdot\bfn_K c_\tB, & \text{if } \bfu\cdot\bfn_K < 0.
\end{cases}
\end{align}
\label{eq:upwind}%
\end{subequations}
An equivalent upwind approximation of $c_h$ is used on $\partial(K\cap\Gamma)$.

Recall the definition of the vector valued jump in Eq.~\eqref{eq:jump_vector}. The coupling term in Eq.~\eqref{eq:integral_eq_fracture}, with $c$ replaced by $c_h$, can be written as
\begin{align}
\jump{\bfu\cdot\bfn c_h^*} = (\bfu\cdot\bfn_\Gamma c_h^*)_- - (\bfu\cdot\bfn_\Gamma c_h^*)_+, \quad \textup{on } K\cap\Gamma.
\end{align}
Given a velocity approximation that is continuous in the interior of an element, we see that $\jump{\bfu\cdot\bfn c_h^*}$ vanish as long as $\Gamma$ is not aligned with the element faces. This would result in no coupling between the fracture and the matrix. To overcome this, we approximate the flow between matrix and fracture by evaluating $(\bfu\cdot\bfn_\Gamma c_h^*)_\pm$ on the part of the element boundary that borders to matrix elements, i.e.,
\begin{align}
\int_{K\cap\Gamma} \jump{\bfu\cdot\bfn c_h^*} \approx \sum_{\tilde{K}\in\Kh^\tM} \int_{\partial\tilde{K}\cap\pK} \bfu\cdot\bfn_K c_h,
\label{eq:coupling_approx}
\end{align}
where the upwind scheme \eqref{eq:upwind} applies. Since $c_h$ on $K\in\Kh^\tF$ represents the approximation in the fracture, \eqref{eq:coupling_approx} is compatible with the condition \eqref{eq:upwind_conc_gamma}. Given the low order method, the approximation \eqref{eq:coupling_approx} seems reasonable when combined with local refinement around $\Gamma$.

We apply the implicit Euler (IE) method as time integrator. For simplicity, we use constant time steps $\Delta t$, and let $c_h^n$ denote the approximation at $t=n\Delta t$, with $c_h^0 = Q_0c_0$, where $Q_0$ is a projection operator from $L^2(\Omega)$ to $Q_0(\Kh)$.

To sum up, the FV-IE scheme can be formulated as follows. Find $c_h^{n+1}\in Q_0(\Kh)$ such that
\begin{subequations}
\begin{align}
\int_K \phi\frac{c_h^{n+1}-c_h^n}{\Delta t} + \int_{\pK}\bfu\cdot\bfn_Kc_h^{n+1} 
&= \int_K f^{n+1}(c_h^{n+1}), 
&\forall K\in\Kh^\tM, 
\label{eq:FV-IE_matrix} \\
\int_{K\cap\Gamma} w\phi_\Gamma\frac{c_h^{n+1}-c_h^n}{\Delta t} 
+ \int_{\partial(K\cap\Gamma)} \bfu_\Gamma\cdot\bfn_{K\cap\Gamma} c_h^{n+1}
&- \sum_{\tilde{K}\in\Kh^\tM} \int_{\partial\tilde{K}\cap\pK} \bfu\cdot\bfn_K c_h^{n+1} \nonumber \\
&= \int_{K\cap\Gamma} f_\Gamma^{n+1}(c_h^{n+1}), 
&\forall K\in\Kh^\tF.
\label{eq:FV-IE_fracture}
\end{align}
\label{eq:FV-IE}%
\end{subequations}

\paragraph{Interpretation of solution}

\begin{figure}[tbp]
	\centering
	\subfloat[Partition of $\Kh$ into matrix ($\Kh^\tM$) and fracture ($\Kh^\tF$) elements.]{
	\includegraphics[height=0.25\textwidth]{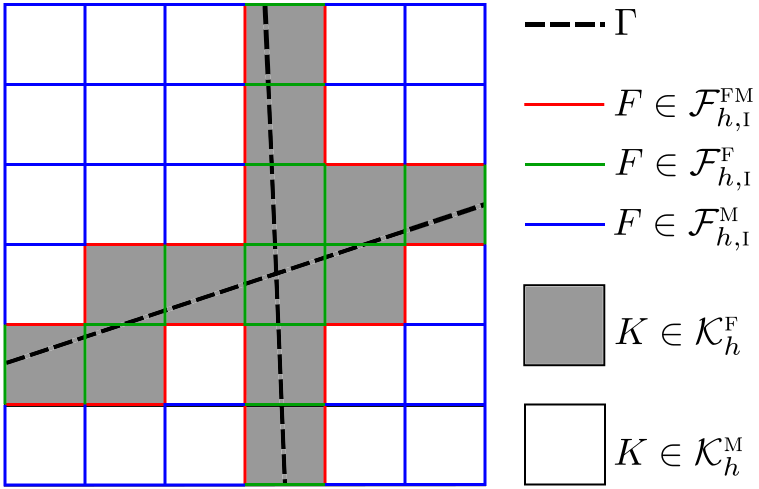}
	\label{fig:fracture_elements}}
	\hspace{3mm}
	\subfloat[Original numerical concentration solution, $c_h\in Q_0(\Kh)$.]{
	\includegraphics[width=0.25\textwidth]{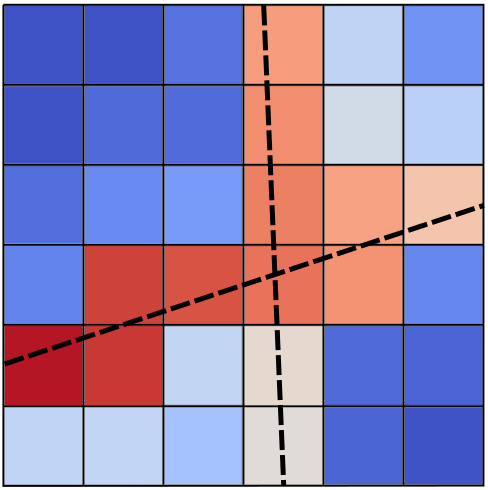}}
	\hspace{3mm}
	\subfloat[Interpreted concentration solution on fracture elements $K\in\Kh^\tF$.]{
	\includegraphics[width=0.25\textwidth]{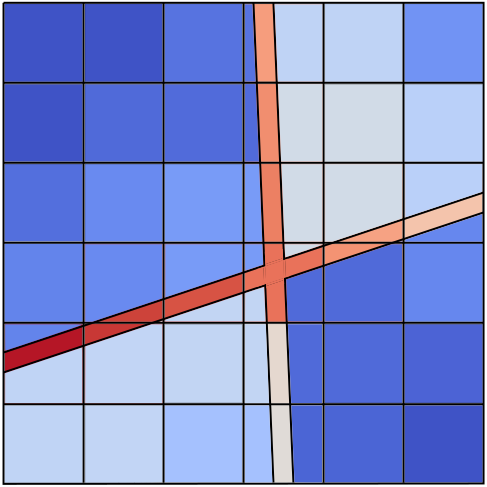}}
	\caption{Synthetic example of partition into matrix and fracture elements, and corresponding interpretation of concentration solution. The interpreted solution on interface elements $K\cap\Gamma$ is visualized with a fixed (exaggerated) thickness.}
	\label{fig:fracture_conc_interpretation}
\end{figure}

In the numerical method, $c\vert_{K\cap\Gamma}$, for $K\in\Kh^\tF$, is represented by the value $c_h\vert_K$. However, $K$ also contains subdomains belonging to the matrix, whose concentration solution we represent by the solution in the neighboring matrix elements, see Fig.~\ref{fig:fracture_conc_interpretation} for an illustrative example. For a single fracture the matrix concentration in $K$ to the left/right of $\Gamma$ is given by the solution in the left/right-neighboring matrix elements. For cells $K$ with intersecting or bifurcating fractures, this interpretation is slightly more complex, as $\Gamma$ divides $K$ into more than two subdomains. We refer to Appendix \ref{app_interpretation} for a well-defined interpretation.

\subsection{Velocity model}

We observe from the transport model, Eq.~\eqref{eq:FV-IE}, that we need the flux (normal velocity) over the element faces. We denote the flux by $U$ and define it as
\begin{align}
U =
\begin{cases}
-\bfK\nabla p\cdot\bfn_F, & \text{on } F\in\Fh{}\setminus\lp\Fh{\tI}^\tF\cup\Fh{\tN}\rp, \\
-\tfrac{1}{\vert F\vert} k_\Gamma\nabla_\Gamma p\cdot \bft_{\Gamma,F}, & \text{on } F\in\Fh{\tI}^\tF, \\
u_\tN, & \text{on } F\in\Fh{\tN},
\end{cases}
\label{eq:velocity}
\end{align}
where $\bft_{\Gamma,F}$ is the unit tangent to $\Gamma$ oriented in the same direction as $\bfn_F$.
Recall that $k_\Gamma=w\kappa_\Gamma$ so that $k_\Gamma\nabla_\Gamma p\cdot\bft_\Gamma$ gives the flow rate rather than the velocity.
The reason for multiplying by $1/\vert F\vert$ for faces $F\in\Fh{\tI}^\tF$ is that we want to work directly on  $\Fh{}$ so that when $U$ is integrated over $F$ we get the total flow rate through $F$. With this definition of $U$ we can simplify two terms in Eq.~\eqref{eq:FV-IE_fracture},
\begin{align}
\int_{\partial(K\cap\Gamma)} \bfu_\Gamma\cdot\bfn_{K\cap\Gamma} c_h^{n+1}
- \sum_{\tilde{K}\in\Kh^\tM} \int_{\partial\tilde{K}\cap\pK} \bfu\cdot\bfn_K c_h^{n+1}
= \int_{\pK} U c_h^{n+1} \bfn_K\cdot\bfn_F.
\label{eq:FV_simplification}
\end{align}
The term $\bfn_K\cdot\bfn_F$ is either plus or minus $1$ depending on the orientation of $F$.

For the coupled flow and transport problem, the velocity is derived from the pressure approximation $p_h$. The pressure gradient, $\nabla p_h$, is not continuous across element faces, so a naive flux approximation is to use the average value,
\begin{align}
U_h =
\begin{cases}
-\average{\bfK\nabla p_h\cdot\bfn_F}_\vartheta, & \text{on } F\in\Fh{}\setminus\lp\Fh{\tI}^\tF\cup\Fh{\tN}\rp, \\
-\tfrac{1}{\vert F\vert}\average{k_\Gamma\nabla_\Gamma p_h\cdot \bft_\Gamma}, & \text{on } F\in\Fh{\tI}^\tF, \\
u_\tN, & \text{on } F\in\Fh{\tN}.
\end{cases}
\label{eq:Uh}
\end{align}
Recall the definition of the average operator in Eq.~\eqref{eq:average_def}.
As weights we follow \cite{odsaeter2016pon} and use weights equal to the normal component of the permeability of the neighboring cell. Hence,
\begin{align}
\vartheta_F = \frac{\delta_{\textup{Kn}}^+}{\delta_{\textup{Kn}}^+ + \delta_{\textup{Kn}}^-}, \quad \delta_{\textup{Kn}}^\pm= \bfn_F\cdot(\bfK_\pm\bfn_F),
\end{align}
where $\bfK_\pm$ are the permeabilities of the two cells sharing $F$. For a fractured cell we use $\bfK=\kappa_\Gamma\bfI$, and for $F\in\Fh{\tB}$, $\omega\vert_F = 1/\delta_{Kn}^-$.
In the case of isotropic permeability, i.e., $\bfK=\kappa\bfI$, observe that $\vartheta_F=\kappa_+/(\kappa_+ + \kappa_-)$ such that
\begin{align}
U_h = -\average{\bfK\nabla p_h\cdot\bfn_F}_\vartheta = 
- \frac{\kappa_+\kappa_-}{\kappa_+ + \kappa_-} (\nabla p_h)_-\cdot \bfn_F
- \frac{\kappa_-\kappa_+}{\kappa_+ + \kappa_-} (\nabla p_h)_+\cdot \bfn_F
= - k_e\average{\nabla p_h \cdot \bfn_F},
\end{align}
where $k_e = \tfrac{2\kappa_+\kappa_-}{\kappa_+ + \kappa_-}$ is the effective face permeability (harmonic average).

We say that a flux approximation $U_h$ on $\Fh{}$ is locally conservative if
\begin{align}
\int_{\pK} U_h(\bfn_K\cdot\bfn_F) = \int_K q, \quad \forall K\in\Kh.
\end{align} 
As reported in \cite{odsaeter2016pon}, Eq.~\eqref{eq:Uh} does not define a locally conservative flux approximation. If coupled to the transport scheme one may get unphysical solutions due to artificial sinks and sources. To deal with this, we apply the postprocessing method presented in \cite{odsaeter2016pon}. This method was shown to preserve accuracy of the velocity solution and was demonstrated to be especially beneficial for highly heterogeneous media.

\paragraph{Variationally consistent postprocessing of fluxes}

The core idea of the postprocessing method is to add a piecewise constant correction to $U_h$ under the constraint that the correction is minimized in a weighted $L^2$ norm,
\begin{align}
\Vert v\Vert_{\omega,\Fh{}} = \sqrt{\ip{\omega v}{v}{\Fh{}}},
\end{align}
where $\omega$ are positive and bounded weights. We define the weight on $F\in\Fh{}$ as the inverse of the effective normal component of the permeability,
\begin{align}
\omega\vert_F = \frac{\delta_{\textup{Kn}}^+ + \delta_{\textup{Kn}}^-}{2\delta_{\textup{Kn}}^+\delta_{\textup{Kn}}^-}.
\end{align}
These weights were demonstrated to be a good choice for heterogeneous permeability as low permeable interfaces are better preserved compared to minimizing in the standard $L^2$ norm ($\omega=1$) \cite{odsaeter2016pon}.

Next, we define a residual operator, $\mathcal{R}:L^2(\Fh{})\rightarrow Q_0(\Kh)$, measuring the discrepancy from local conservation,
\begin{align}
\mathcal{R}(U)\vert_K =
\begin{cases}
\tfrac{1}{\vert K\vert}\lp \int_K q - \int_{\pK} U \bfn_F\cdot\bfn_K\rp, & K\in\Kh^\tM, \\
\tfrac{1}{\vert K\vert}\lp \int_{K\cap\Gamma} q_\Gamma - \int_{\pK} U \bfn_F\cdot\bfn_K\rp, & K\in\Kh^\tF.
\end{cases}
\end{align}
With this we define the postprocessed flux, $V_h$, as follows.
\begin{align}
V_h =
\begin{cases}
U_h + \omega^{-1}\jump{y}, & \text{on } F\in\Fh{}\setminus\Fh{\tN}, \\
u_\tN, & \text{on } F\in\Fh{\tN},
\end{cases}
\label{eq:velocity_pp}
\end{align}
where $y\in Q_0(\Kh)$ is the unique solution to
\begin{align}
\ip{\omega^{-1}\jump{y}}{\jump{w}}{\Fh{}\setminus\Fh{\tN}} = 
\ip{\mathcal{R}(U_h)}{w}{\Kh}, \quad \forall w\in Q_0(\Kh).
\end{align}
For further details on the postprocessing method, we refer to \cite{odsaeter2016pon}.

\paragraph{Coupled formulation}

Applying Eq.~\eqref{eq:FV_simplification} and using $V_h$ as an approximation to $U$, the FV-IE scheme, Eq.~\eqref{eq:FV-IE}, can be formulated as follows. Find $c_h^{n+1}\in Q_0(\Kh)$ such that
\begin{subequations}
\begin{align}
\int_K \phi\frac{c_h^{n+1}-c_h^n}{\Delta t} + \int_{\pK}V_h c_h^{n+1} \bfn_K\cdot\bfn_F
&= \int_K f^{n+1}(c_h^{n+1}), 
&\forall K\in\Kh^\tM, 
\label{eq:FV-IE_pp_matrix} \\
\int_{K\cap\Gamma} w\phi_\Gamma\frac{c_h^{n+1}-c_h^n}{\Delta t} 
+ \int_{\pK} V_h c_h^{n+1} \bfn_K\cdot\bfn_F
&= \int_{K\cap\Gamma} f_\Gamma^{n+1}(c_h^{n+1}), 
&\forall K\in\Kh^\tF.
\label{eq:FV-IE_pp_fracture}
\end{align}
\label{eq:FV-IE_pp}%
\end{subequations}
Note that the second term in Eq.~\eqref{eq:FV-IE_pp_fracture} contains both the flux along the fracture and the coupling term, cf.~Eq.~\eqref{eq:FV_simplification}.

\section{Numerical results}
\label{sec:results}
In this section we demonstrate the numerical methods presented in Section~\ref{sec:methods}. First, in Section~\ref{sec:results_explicit}, we consider a pure transport problem where the velocity is given explicitly and the exact solution is known. This is to verify our transport model, Eq.~\eqref{eq:FV-IE_pp}, and in particular our approximation of the coupling term, see Eq.~\eqref{eq:coupling_approx}.

Next, in Section~\ref{sec:results_pressure}, we consider two benchmark cases for the pure pressure problem presented in \cite{flemisch2018bfs}. The first case is a regular fracture network, while the second problem is a realistic case with a complex fracture network. The aim is to compare EFEM, Eq.~\eqref{eq:efem}, to other DFM models for single-phase flow.

Finally, in Section~\ref{sec:results_coupled}, we solve the coupled pressure and transport problem on the same cases as in Section~\ref{sec:results_pressure}. This will reveal the capabilities of our solution approach.

All implementation of the numerical methods are based on the open-source software deal.II \cite{bangerth2007dag}. We only consider 2D problems, but our fracture model and numerical methods can be applied to 3D problems as well. All meshes are built up of quadrilateral elements. The meshes may be locally refined by recursively dividing selected elements in four, but we allow for no more than one hanging node per element face. The number of degrees of freedom are denoted $\ndof$.

\subsection{Pure transport problem}
\label{sec:results_explicit}
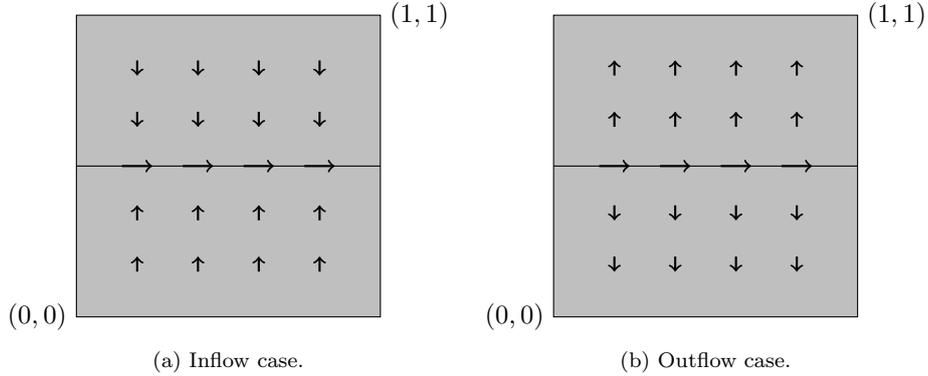
\begin{figure}[tbp]
	\centering
	\subfloat[Inflow case.]{
	\begin{tikzpicture}[scale=4]
	\draw[fill=lightgray] (0,0) node[left]{$(0,0)$} -- (1,0) -- (1,1) node[right]{$(1,1)$} -- (0,1) -- cycle;
	\draw (0,0.5) -- (1,0.5);
	\draw[->, thick] (0.15,0.5) -- (0.25,0.5);
	\draw[->, thick] (0.35,0.5) -- (0.45,0.5);
	\draw[->, thick] (0.55,0.5) -- (0.65,0.5);
	\draw[->, thick] (0.75,0.5) -- (0.85,0.5);
	\foreach \y in {0.85, 0.68} {
		\foreach \x in {0.2,0.4,0.6,0.8} {
			\draw[->, thick] (\x,\y) -- (\x,\y-0.05);
		}
	}
	\foreach \y in {0.15, 0.32} {
		\foreach \x in {0.2,0.4,0.6,0.8} {
			\draw[->, thick] (\x,\y) -- (\x,\y+0.05);
		}
	}
	\end{tikzpicture}
	}
	\subfloat[Outflow case.]{
		\begin{tikzpicture}[scale=4]
		\draw[fill=lightgray] (0,0) node[left]{$(0,0)$} -- (1,0) -- (1,1) node[right]{$(1,1)$} -- (0,1) -- cycle;
		\draw (0,0.5) -- (1,0.5);
		\draw[->, thick] (0.15,0.5) -- (0.25,0.5);
		\draw[->, thick] (0.35,0.5) -- (0.45,0.5);
		\draw[->, thick] (0.55,0.5) -- (0.65,0.5);
		\draw[->, thick] (0.75,0.5) -- (0.85,0.5);
		\foreach \y in {0.85, 0.68} {
			\foreach \x in {0.2,0.4,0.6,0.8} {
				\draw[<-, thick] (\x,\y) -- (\x,\y-0.05);
			}
		}
		\foreach \y in {0.15, 0.32} {
			\foreach \x in {0.2,0.4,0.6,0.8} {
				\draw[<-, thick] (\x,\y) -- (\x,\y+0.05);
			}
		}
		\end{tikzpicture}
	}
	\caption{Pure transport problem: Description of the two cases. Arrows describe the velocity field.}
	\label{fig:explicit_velocity}
\end{figure}

We consider first a pure transport problem with an explicitly given velocity field. Let $\Omega=(0,1)^2$ and $\Gamma=(0,1)\times\{0.5\}$. The fracture velocity, $u_\Gamma=\bfu_\Gamma\cdot(1,0)=10$, and we consider two cases for the matrix velocity, $\bfu$. Either $\bfu=(0,1)$ for $y<0.5$ and $\bfu=(0,-1)$ otherwise, or $\bfu=(0,-1)$ for $y<0.5$ and $\bfu=(0,1)$ otherwise. The two cases are depicted in Fig.~\ref{fig:explicit_velocity} and are denoted \emph{inflow} and \emph{outflow}, respectively. In both cases we set $w=1$ and use initial and boundary conditions $c_0=0$ and $c_\tB=1$.

The 1D advective transport equation describing the fracture concentration, $c_\Gamma$, is given as
\begin{subequations}
\begin{align}
w\frac{\partial c_\Gamma}{\partial t} + u_\Gamma \frac{\partial c_\Gamma}{\partial x} - \jump{\bfu\cdot\bfn c^*} &= 0, \qquad\text{on } \Gamma, \\
c_\Gamma &= 1, \qquad\text{at } x=0.
\end{align}
\end{subequations}
For the inflow case we have $c^*=c$, and for the outflow case $c^*=c_\Gamma$. Both cases are solved on uniform $N\times N$ meshes with time steps $\dt=0.001$, and the simulations are run until a steady-state solution is reached.

\paragraph{Inflow case}

\begin{figure}[tbp]
	\centering
	\subfloat[$t=0.05$.]{
		\includegraphics[width=0.48\textwidth]{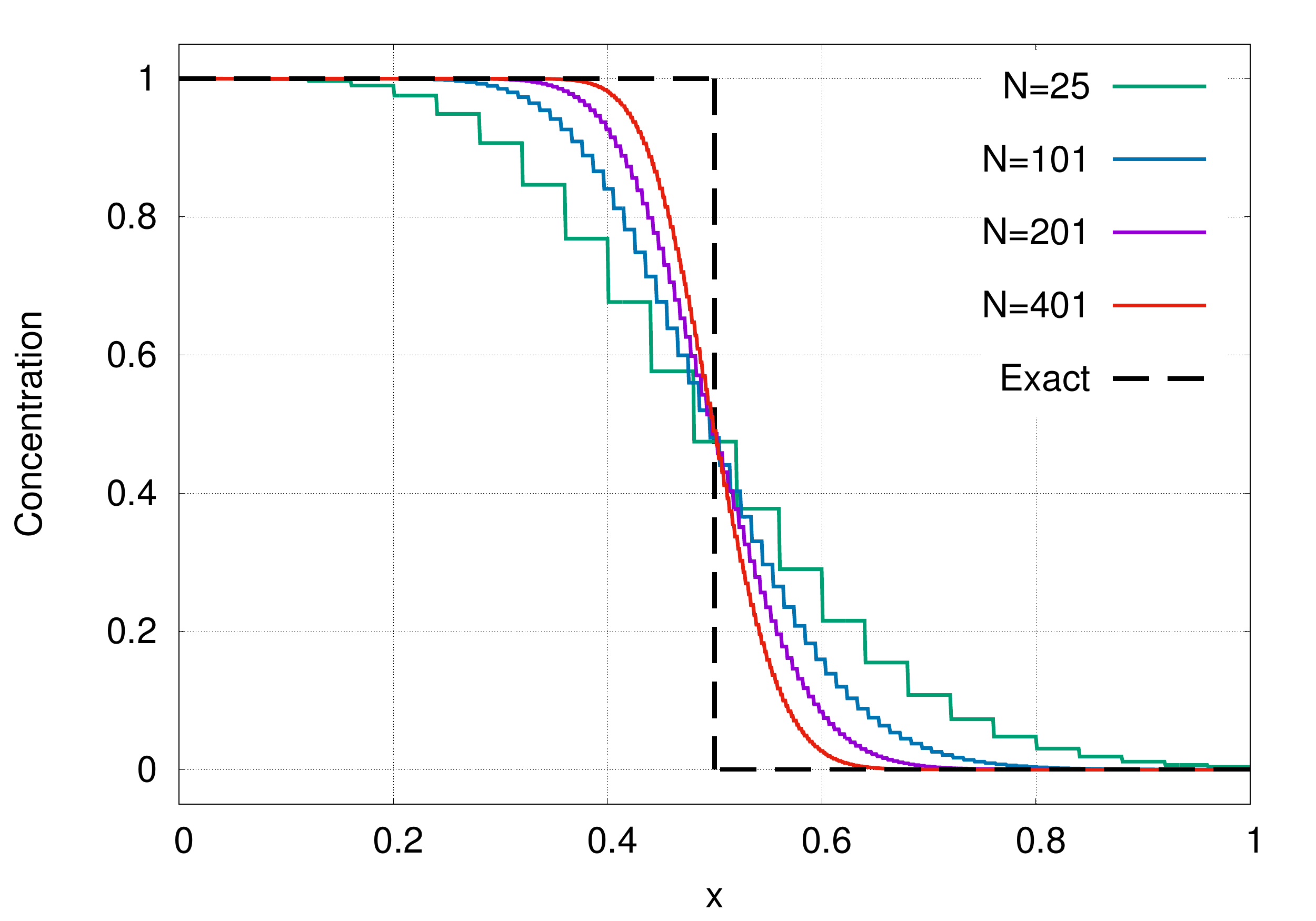}
	}
	\subfloat[Steady-state solution, $t=1.5$.]{
		\includegraphics[width=0.48\textwidth]{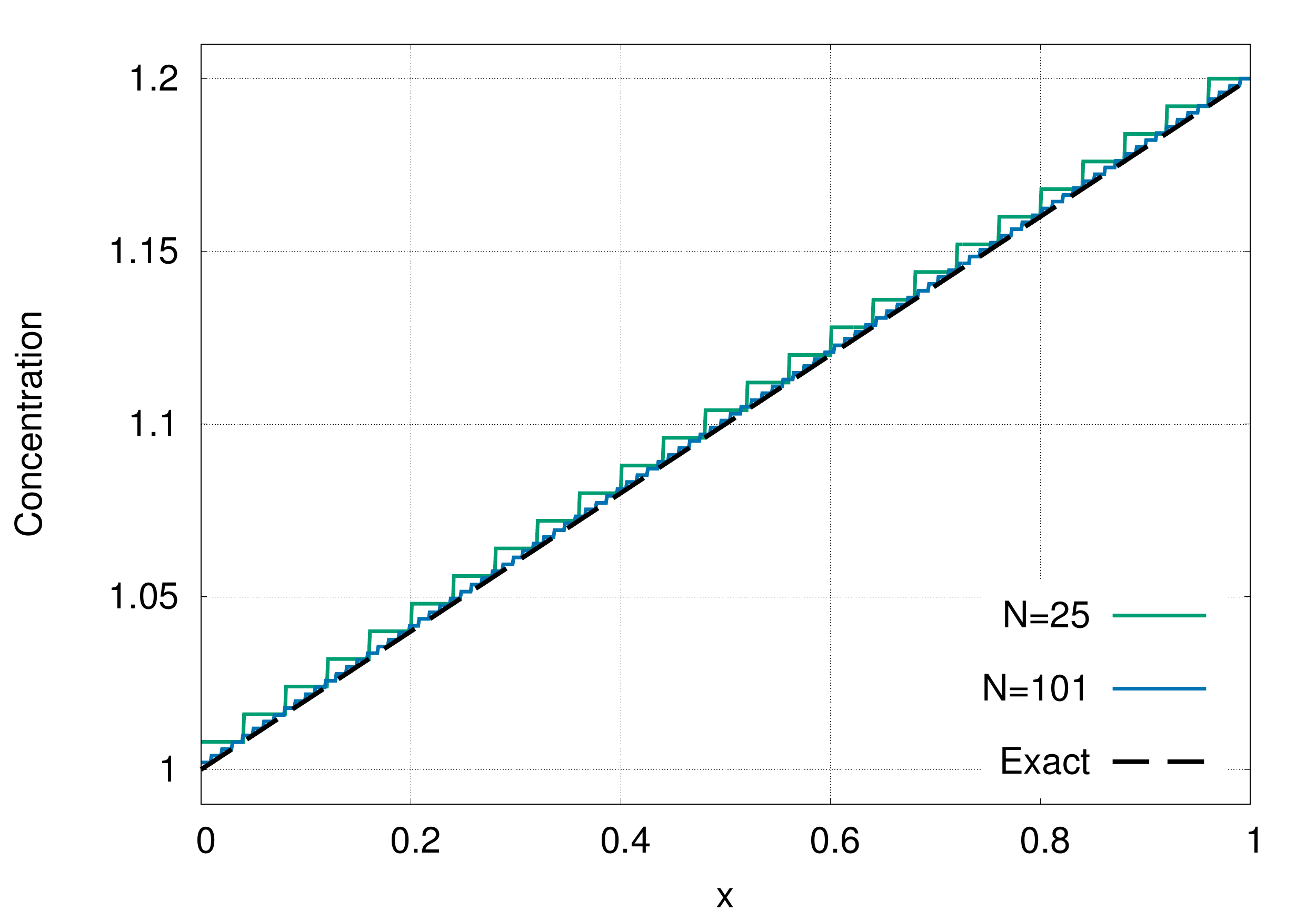}
	}
	\caption{Pure transport problem --- inflow case: Fracture concentration on uniform $N\times N$ meshes compared to exact solution.}
	\label{fig:explicit_velocity_inflow}
\end{figure}

For the inflow case the concentration front from the top and bottom boundary moves with speed 1, so that $\jump{\bfu\cdot\bfn c^*} = 2c$ with $c=0$ for $t<0.5$ and $c=1$ for $t>0.5$. Hence, the exact solution for $t<0.5$ reads
\begin{align}
c_\Gamma(x,t) = 
\begin{cases}
1, & x < \tfrac{u_\Gamma}{w} t, \\
0, & x \ge \tfrac{u_\Gamma}{w} t.
\end{cases}
\end{align}
At steady-state, $c=1$ in the matrix, and the exact steady-state solution is
\begin{align}
c_\Gamma^{\textup{ss}}(x) = 1+\tfrac{2w}{u_\Gamma}x.
\end{align}
The approximation, $c_{h}$, along the fracture is plotted at $t=0.05$ and $t=1.5$ (steady-state) in Fig.~\ref{fig:explicit_velocity_inflow} and compared to the exact solution.

\paragraph{Outflow case}

\begin{figure}[tbp]
	\centering
	\subfloat[$t=0.05$.]{
		\includegraphics[width=0.48\textwidth]{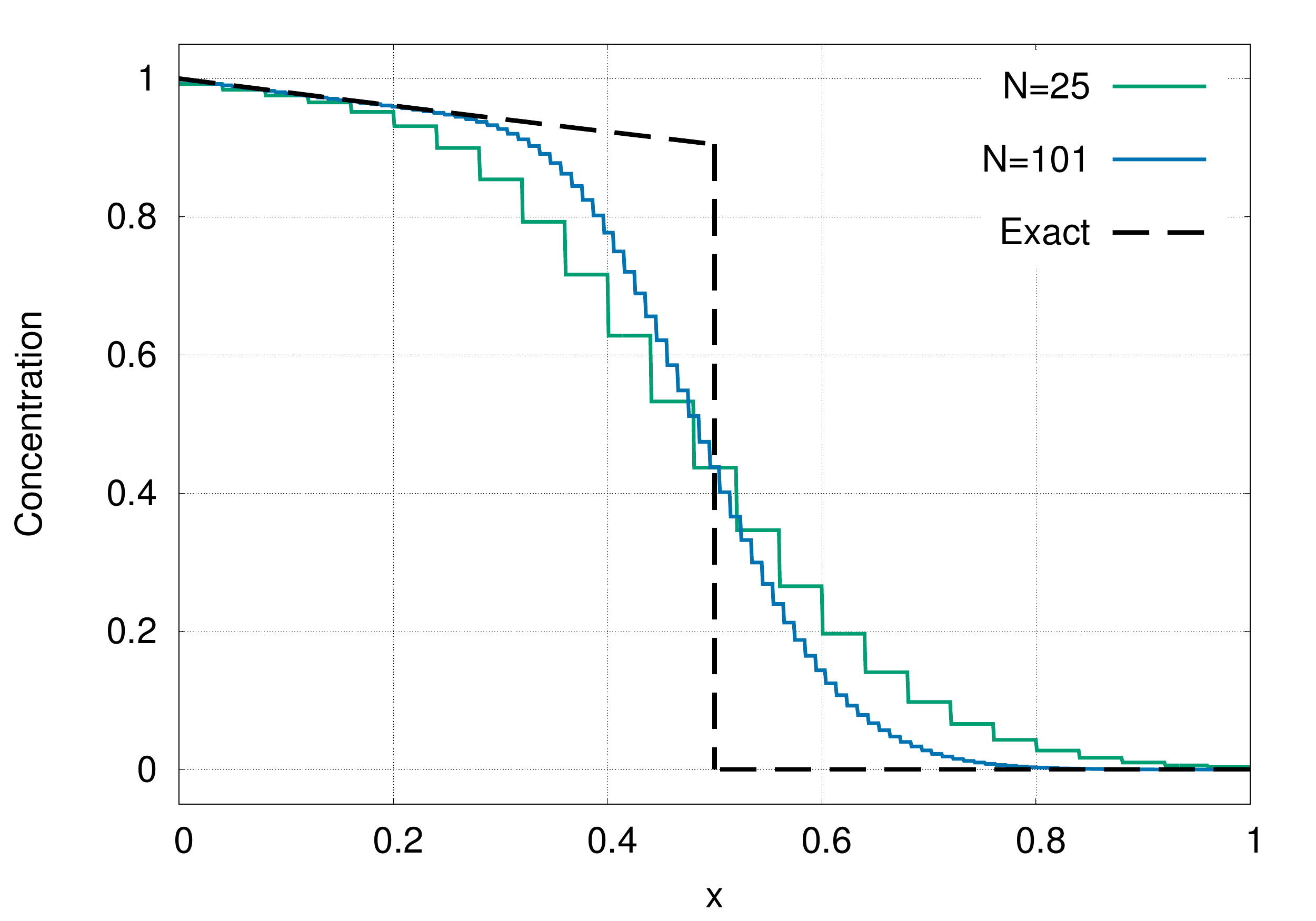}
	}
	\subfloat[Steady-state solution, $t=1.5$.]{
		\includegraphics[width=0.48\textwidth]{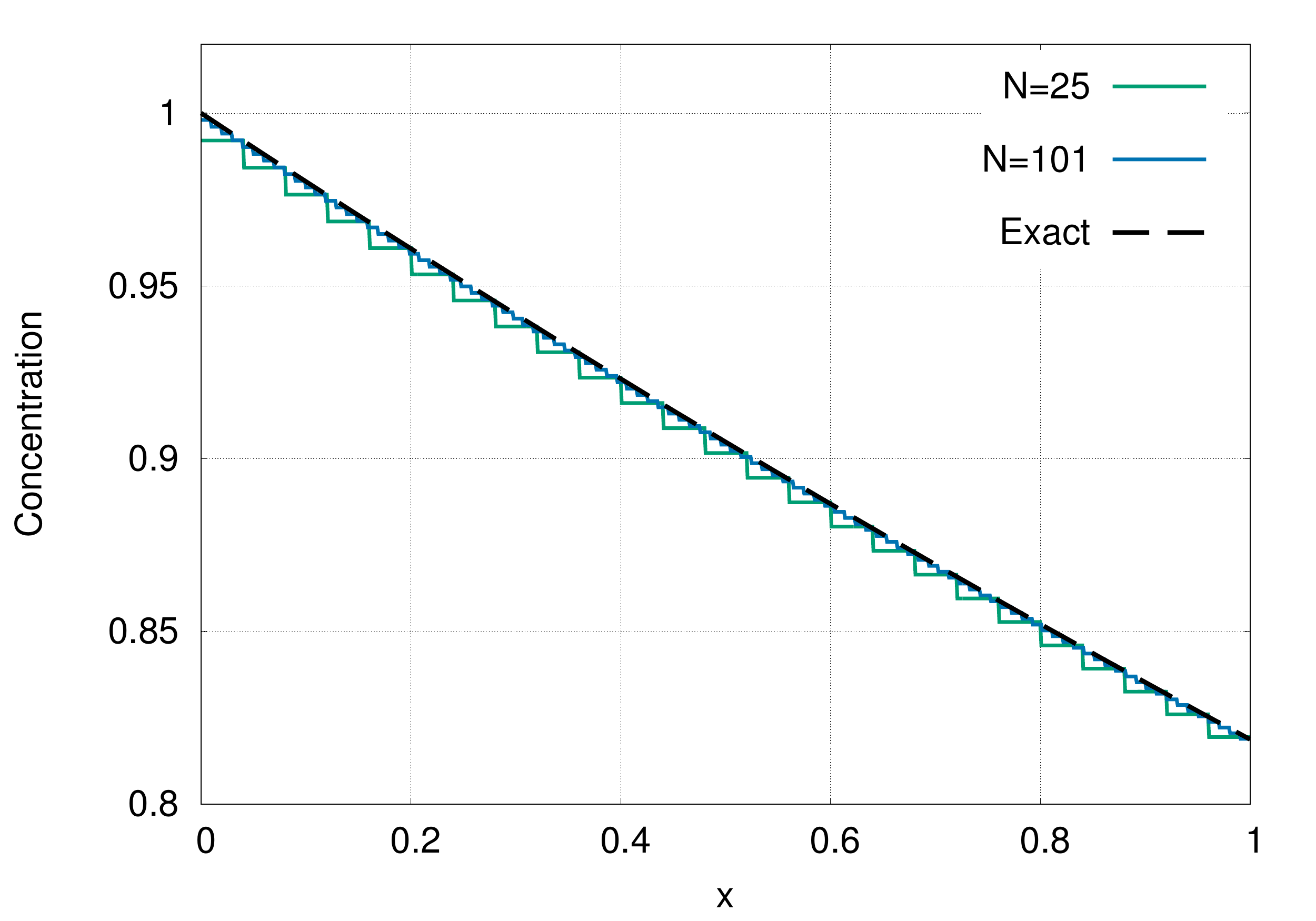}
	}
	\caption{Pure transport problem --- outflow case: Fracture concentration on uniform $N\times N$ meshes compared to exact solution.}
	\label{fig:explicit_velocity_outflow}
\end{figure}

For the outflow case, $\jump{\bfu\cdot\bfn_{\Gamma}c^*} = -2c_\Gamma$, so that the exact solution reads
\begin{align}
c_\Gamma(x,t) =
\begin{cases}
\exp\{-\tfrac{2}{u_\Gamma}x\}, & x<u_\Gamma t, \\
0,        & x\ge u_\Gamma t.
\end{cases}
\end{align}
The approximation, $c_{h}$, along the fracture is plotted at $t=0.05$ and $t=1.5$ (steady-state) in Fig.~\ref{fig:explicit_velocity_outflow} and compared to the exact solution.

We observe that we get the correct steady-state solution in both cases, and that the velocity of the concentration front is correct. Due to the low order method we get significant numerical diffusion as expected, but we see that the front gets sharper as $N$ increases. 

%\clearpage

\subsection{Benchmark cases for the pure pressure problem}
\label{sec:results_pressure}

In this section we consider two of the benchmark cases defined in \cite{flemisch2018bfs} for the pure pressure problem. We only consider the cases where all fractures have higher permeability than the surrounding matrix as our model only applies to such problems. We employ the exact same problem setup as \cite{flemisch2018bfs} and compare our results with those given therein\footnote{All results reported in \cite{flemisch2018bfs} are public available at https://git.iws.unistuttgart.de/benchmarks/fracture-flow.}. 
The six methods considered in this study are listed in Table \ref{tab:benchmark_methods}. We only give a very brief description of these methods here, and refer to \cite{flemisch2018bfs} and references therein for further details. Our method will be denoted EFEM.

\begin{table}[tbp]
	\caption{List of the participating methods in the benchmark paper \cite{flemisch2018bfs}.}
	\label{tab:benchmark_methods}
	\begin{center}
		\begin{tabular}{ll}
			\hline
			\textbf{Method} & \textbf{Description} \\
			\hline
			Box         & Vertex-centered finite-volume method \\
			TPFA        & Control volume finite difference method with two-point flux approximation \\
			MPFA        & Control volume finite difference method with multi-point flux approximation \\
			EDFM        & Embedded discrete fracture-matrix model \\
			Flux-Mortar & Mortar discrete fracture-matrix model \\
			P-XFEM      & Primal extended finite element method \\
			D-XFEM      & Dual extended finite element method \\
			\hline
		\end{tabular}
	\end{center}
\end{table}

\subsubsection{Benchmark 1: Regular fracture network}

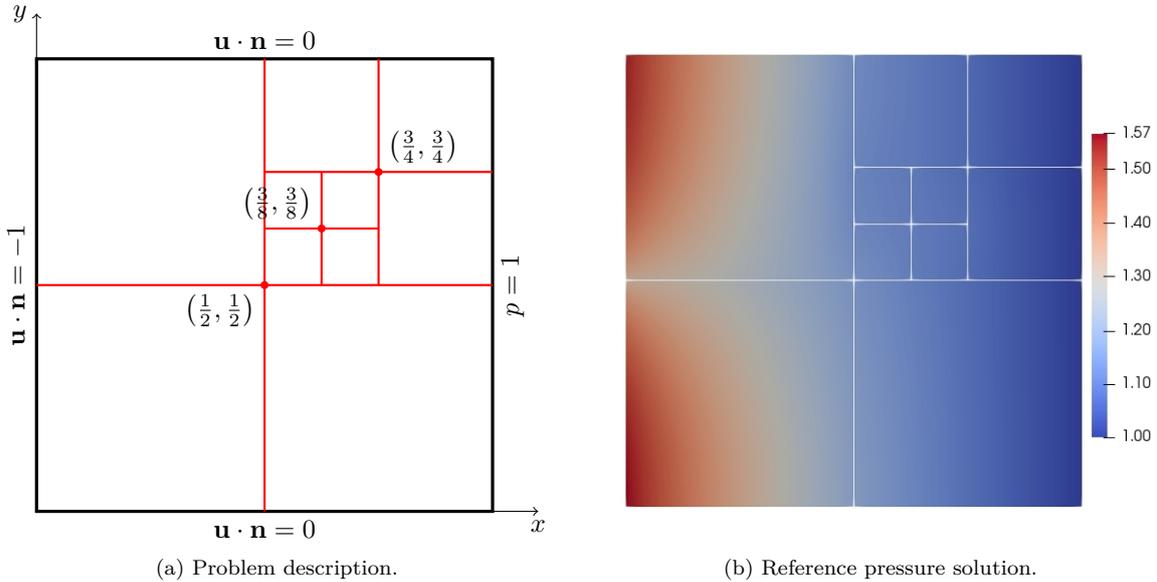
\begin{figure}[tbp]
	\centering
	\subfloat[Problem description.]{
	\begin{tikzpicture}[scale=6]
		\draw[->] (0,0) -- (1.1,0) node[below]{$x$};
		\draw[->] (0,0) -- (0,1.1) node[left]{$y$};
		\draw[very thick] (0,0) -- (1,0) -- (1,1) -- (0,1) -- cycle;
		\draw[thick,red] (0,0.5) -- (1,0.5);
		\draw[thick,red] (0.5,0.75) -- (1,0.75);
		\draw[thick,red] (0.5,0.625) -- (0.75,0.625);
		\draw[thick,red] (0.5,0) -- (0.5,1);
		\draw[thick,red] (0.75,0.5) -- (0.75,1);
		\draw[thick,red] (0.625,0.5) -- (0.625,0.75);
		\draw (0,0.5) node[left,rotate=90,anchor=north,yshift=5mm]{$\bfu\cdot\bfn=-1$};
		\draw (0.5,0) node[below]{$\bfu\cdot\bfn=0$};
		\draw (0.5,1) node[above]{$\bfu\cdot\bfn=0$};
		\draw (1,0.5) node[right,rotate=90,anchor=north,yshift=0mm]{$p=1$};
		\fill[red] (0.625,0.625) circle[radius=0.25pt] node[black,above left]{$\lp\tfrac{3}{8},\tfrac{3}{8}\rp$};
		\fill[red] (0.75,0.75) circle[radius=0.25pt] node[black,above right]{$\lp\tfrac{3}{4},\tfrac{3}{4}\rp$};
		\fill[red] (0.5,0.5) circle[radius=0.25pt] node[black,below left]{$\lp\tfrac{1}{2},\tfrac{1}{2}\rp$};
	\end{tikzpicture}
	\label{fig:benchmark_setup}
	}
	\hspace*{3mm}
	\subfloat[Reference pressure solution.]{
	\raisebox{5mm}{
	\includegraphics[width=0.435\textwidth]{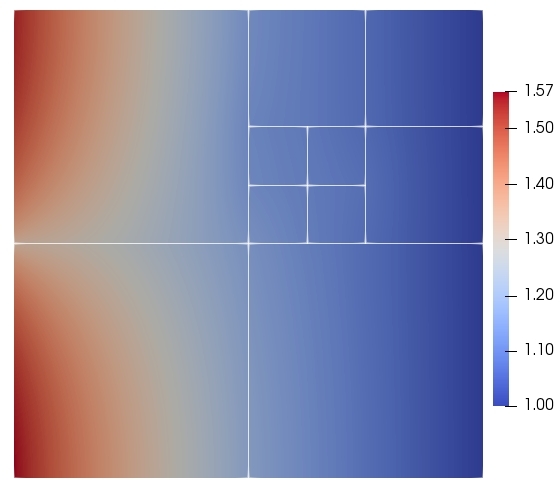}}
	\label{fig:benchmark_ref_pressure}
	}
	\caption{Benchmark 1 (pressure problem).}
\end{figure}

Benchmark 1 is a regular fracture network embedded in the unit square, $\Omega=(0,1)^2$, see Fig.~\ref{fig:benchmark_setup}. The top an bottom boundary faces have homogeneous Neumann conditions (no flow); the left boundary face has a constant inflow flux, $\bfu\cdot\bfn=-1$; and the right boundary face has Dirichlet condition $p_\tD=1$. The rock properties are $\bfK=\bfI$, $\kappa_{\Gamma}=10^4$ and $w=10^{-4}$.

A reference solution is obtained in \cite{flemisch2018bfs} by using a mimetic finite difference (MFD) method on a very fine reference mesh where the fractures are resolved by 10 elements in their normal direction. The reference mesh is coarser away from the fractures and has a total of 1175056 elements. Hence, the fractures are not modeled as a lower dimensional embedding, but as a continuous model with $\bfK=\kappa_\Gamma\bfI$ in the fracture elements. We denote by $K_{\textup{ref}}$ the reference elements and by $p_{\textup{ref}}$ the reference solution. The reference solution is displayed in Fig.~\ref{fig:benchmark_ref_pressure}.

\newcommand{\errm}{err_{\tM}}
\newcommand{\errf}{err_{\tF}}
\newcommand{\pref}{p_{\textup{ref}}}
\newcommand{\Kref}{K_{\textup{ref}}}
To measure the error in the pressure solution, we define two error functions, $\errm$ and $\errf$, measuring the error in the matrix and fractures, respectively,
\begin{align}
\errm^2 &= \frac{1}{\vert\Omega\vert(\Delta \pref)^2}
\sum_{f=\Kref\cap K}\vert f\vert (p_h\vert_{f_m} - \pref\vert_{\Kref})^2
\approx \frac{1}{\vert\Omega\vert(\Delta \pref)^2} \Vert p_h-\pref\Vert^2_{L^2(\Omega)}, \\
\errf^2 &= \frac{1}{\vert\Gamma\vert(\Delta \pref)^2}
\sum_{e=\lp \Kref\cap K \rp\cap \Gamma}\vert e\vert (p_h\vert_{e_m} - \pref\vert_{\Kref})^2
\approx \frac{1}{\vert\Gamma\vert(\Delta \pref)^2} \Vert p_h-\pref\Vert^2_{L^2(\Gamma)},
\end{align}
where $\Delta \pref=\max_{\Omega}\pref - \min_{\Omega}\pref$ and $f_m$ and $e_m$ denotes the midpoints of $f$ and $e$, respectively. Observe that these are $L^2$ errors where the integrals are approximated by the midpoint rule.

%\begin{align}
%\Vert p_{\textup{ref}} \Vert_{L^2(\Kh)} &= 2.131, \\
%\Vert p_{\textup{ref}} \Vert_{L^2(\Gamma)} &= 2.003.
%\end{align}

\begin{figure}[tbp]
	\centering
	\includegraphics[width=0.31\textwidth]{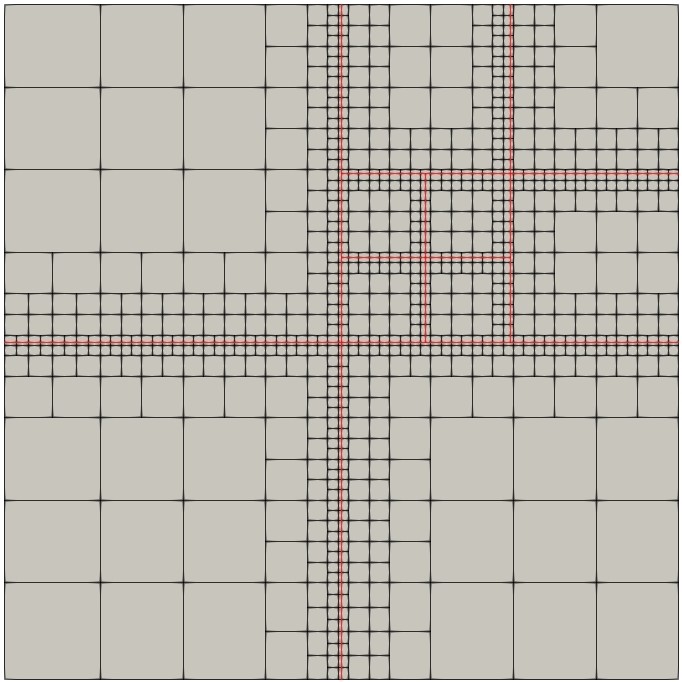}
	\hspace{1mm}
	\includegraphics[width=0.31\textwidth]{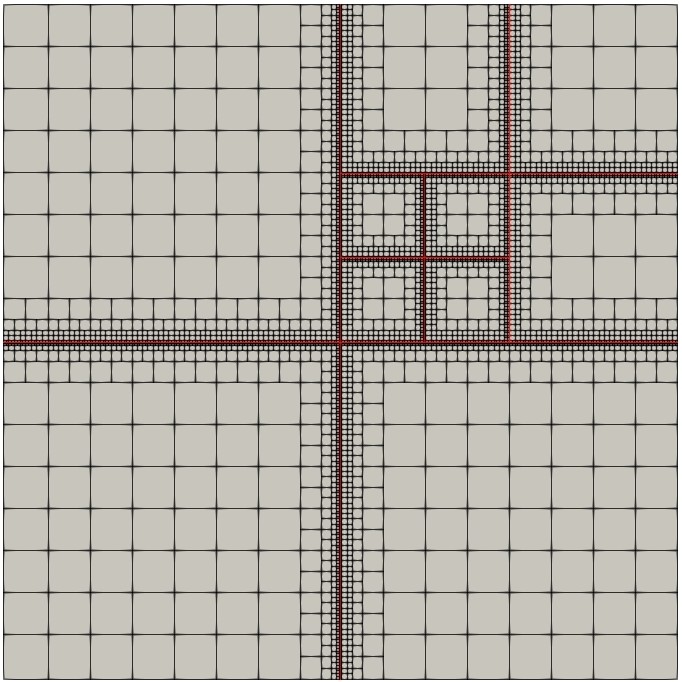}
	\hspace{1mm}
	\includegraphics[width=0.31\textwidth]{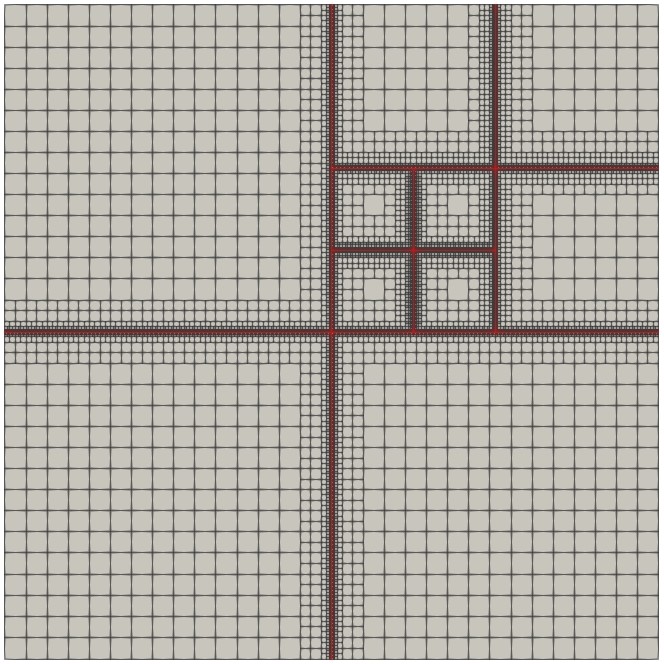}
	\caption{Benchmark 1 (pressure problem). Locally refined meshes with $\ndof$ equal to 932, 4532 and 19579, respectively. Fractures are included as red lines. We avoid fractures along element faces.}
	\label{fig:benchmark_lr_mesh}
\end{figure}

\begin{figure}[tbp]
	\centering
	\includegraphics[width=0.9\textwidth]{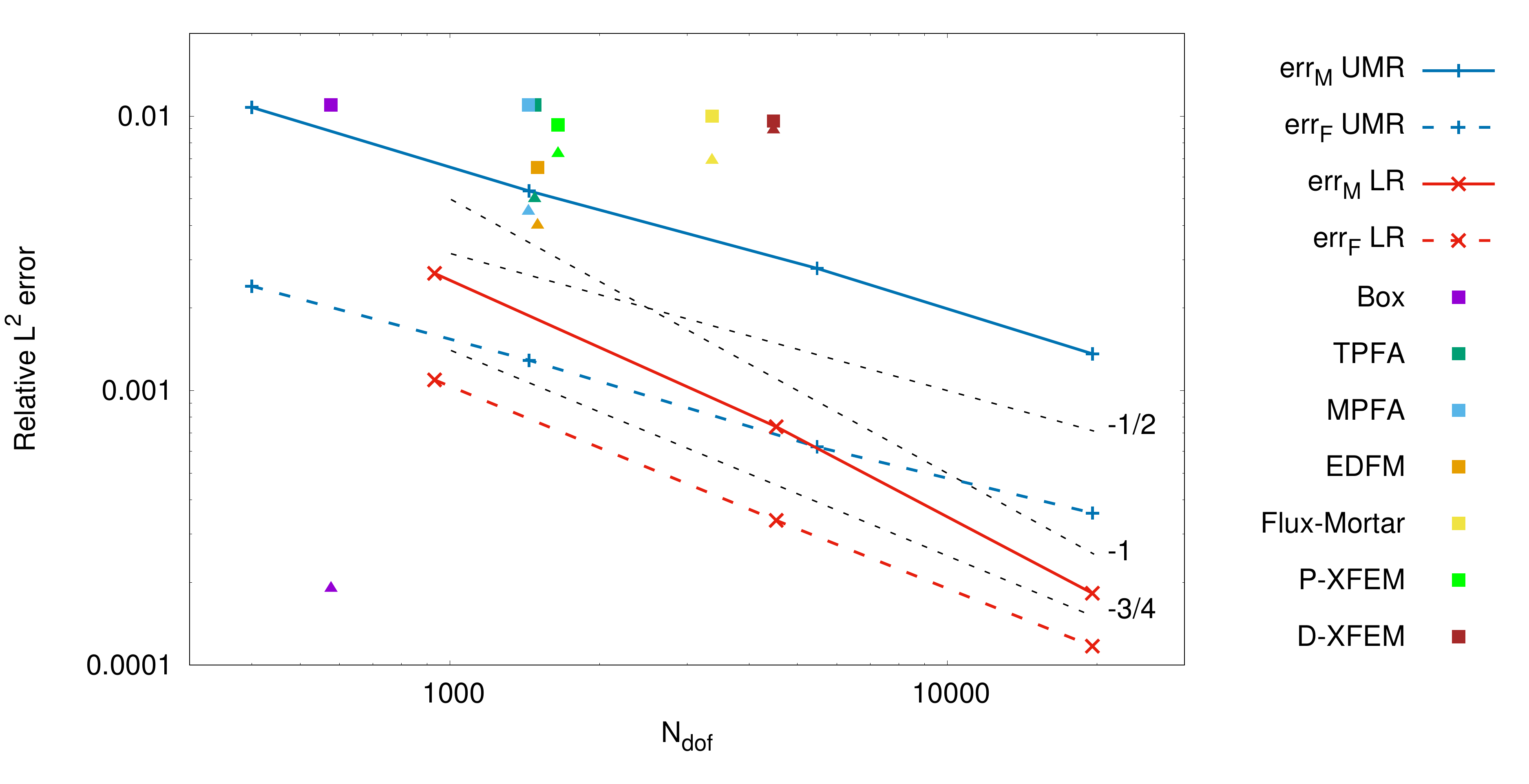}
	\caption{Benchmark 1 (pressure problem). Errors against $\ndof$. Solid lines denote the matrix error $\errm$, and dashed lines denote the fracture error $\errf$ for UMR (blue) and LR (red) meshes. Black dashed lines are straight lines of labeled slope. Filled squares and triangles denote $\errm$ and $\errf$, respectively, for the methods reported in \cite{flemisch2018bfs}.}
	\label{fig:benchmark_conv}
\end{figure}

We solve the pressure problem on both uniform $N\times N$ meshes with $N=\{19,37,73,139\}$, denoted UMR$N$, and three locally refined meshes, denoted LR$i$, for $i=1,2,3$, where the local refinement is based on the a priori estimate \eqref{eq:apriori_l2_pressure} such that $h_{\Gamma}\lesssim h^2$ where $h_{\Gamma}$ is the element size in the vicinity of the fractures. The LR meshes are shown in Fig.~\ref{fig:benchmark_lr_mesh}. The errors ares plotted against $\ndof$ in Fig.~\ref{fig:benchmark_conv}. We see that the error for the uniform meshes has convergence order $\ndof^{-1/2}$ in accordance with the error estimate, while the error is lower and converge faster for the LR meshes.

\begin{figure}[tbp]
	\centering
	\subfloat[Line $y=0.7$.]{
	\includegraphics[width=0.49\textwidth]{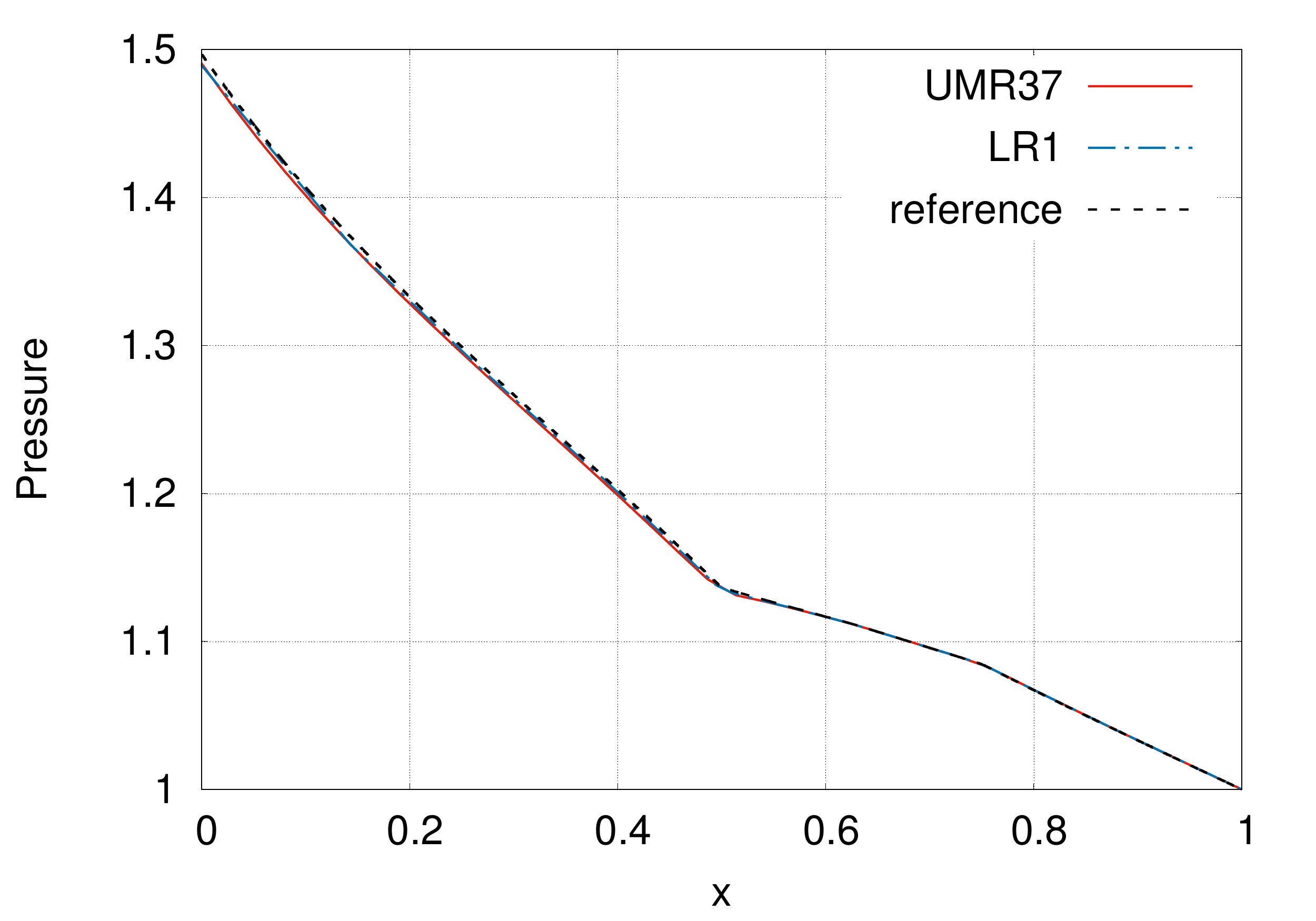}}
	%\hspace{1mm}
	\subfloat[Line $x=0.5$.]{
	\includegraphics[width=0.49\textwidth]{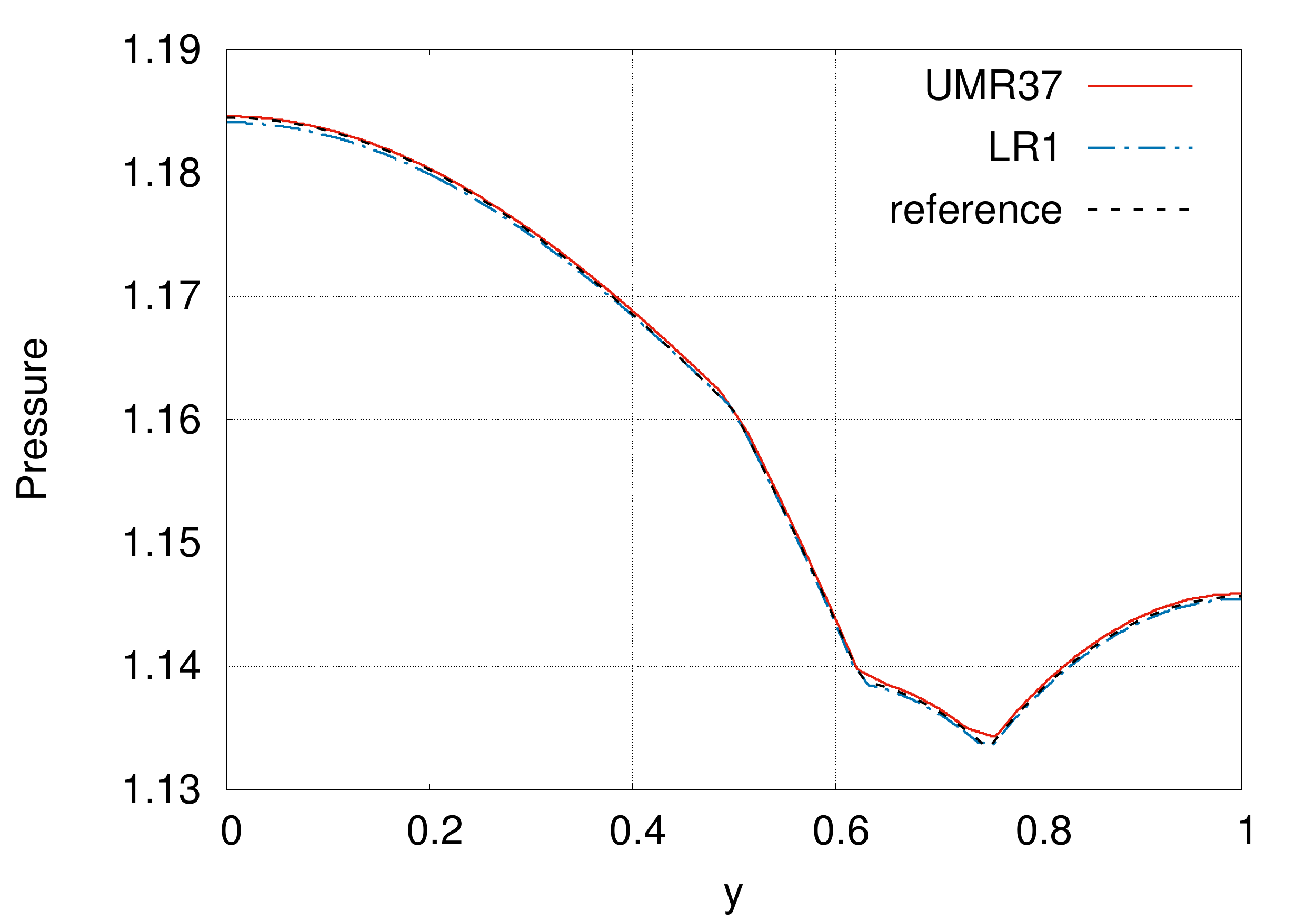}}
	\caption{Benchmark 1 (pressure problem). Pressure solution along two lines.}
	\label{fig:benchmark_pressure_lines}
\end{figure}

In Table~\ref{tab:benchmark_error} we compare our results with the ones reported in \cite{flemisch2018bfs}, and we observe that the results are in good agreement. Furthermore, in Fig.~\ref{fig:benchmark_pressure_lines} the pressure along the lines $y=0.7$ and $x=0.5$ are plotted, similarly to the results reported in \cite{flemisch2018bfs}. We observe a good match with the reference solution, and in particular we see that the LR mesh gives better accuracy close to the fractures.

\begin{remark}
	\label{remark:benchmark_l2norm}
	We would like to point out that the error functions in \cite{flemisch2018bfs} are defined in a similar way, but instead of using $p_h$ directly, the projection of $p_h$ onto piecewise constant functions on the computational mesh is used. In the case where the pressure approximation is piecewise constant this is equivalent to what we do. However, for higher order polynomial approximations, the error functions used in \cite{flemisch2018bfs} would give unfavorable results.
\end{remark}

\begin{table}[tbp]
	\caption{Benchmark 1 (pressure problem). Comparison with the results in \cite{flemisch2018bfs} for error, matrix density (nnz denotes number of non-zero entries in the system matrix) and matrix condition number. \emph{No.\ of elements} for the methods in \cite{flemisch2018bfs} are listed as matrix elements plus fracture elements. In light of Remark~\ref{remark:benchmark_l2norm} we would like to point out that Box and EDFM have continuous pressure approximation, so that the errors associated with these would likely be smaller.}
	\label{tab:benchmark_error}
	\begin{center}
		\begin{tabular}{lrrrrrr}
			\hline
			\textbf{Method} & $\ndof$ & \textbf{No.~of elements} & $\errm$ & $\errf$ & \textbf{nnz/}$\ndof^2$ & $\Vert\cdot\Vert_2\textbf{-cond}$ \\
			\hline
			EFEM UMR37  & 1444 & 1369 & 5.3e-3 & 1.3e-3 & 5.9e-3 & 3.3e4 \\
			EFEM LR1    &  932 &  757 & 2.7e-3 & 1.1e-3 & 7.0e-3 & 1.1e4 \\
			\hline
			Box			&  577 & 1078 + 74  & 1.1e-2 & 1.9e-4 & 1.1e-2 & 2.2e3 \\
			TPFA		& 1481 & 1386 + 95  & 1.1e-2 & 4.4e-3 & 2.7e-3 & 4.8e4 \\
			MPFA	    & 1439 & 1348 + 91  & 1.1e-2 & 4.5e-3 & 8.0e-3 & 5.8e4 \\
			EDFM		& 1501 & 1369 + 132 & 6.5e-3 & 4.0e-3 & 3.3e-3 & 5.6e4 \\
			Flux-Mortar	& 3366 & 1280 + 75  & 1.0e-2 & 6.9e-3 & 1.8e-3 & 2.4e6 \\
			P-XFEM		& 1650 &  961 + 164 & 9.3e-2 & 7.3e-3 & 8.0e-3 & 9.3e9 \\
			D-XFEM		& 4474 & 1250 + 126 & 9.6e-3 & 8.9e-3 & 1.3e-3 & 1.2e6 \\
			\hline
		\end{tabular}
	\end{center}
\end{table}

%\clearpage
\subsubsection{Benchmark 4: A realistic case}

\begin{figure}[tbp]
\centering
%\subfloat[Problem description. Red lines represent fractures, the blue lines represent the lines for which the solution is plotted along in Figs.~\ref{fig:pressure_complex_y500} and \ref{fig:pressure_complex_x625}, and the small black box represents the box plotted in Fig.~\ref{fig:benchmark4_mesh}.]{
\begin{tikzpicture}[scale=1]
\DTLforeach*{coordinates}{\xs=START_X, \xe=END_X, \ys=START_Y, \ye=END_Y}
{\draw[red] (\xs/100,\ys/100) -- (\xe/100,\ye/100);}
\draw[thick] (0,0) -- (7,0) -- (7,6) -- (0,6) -- (0,0); 
\draw (0,0) node[below]{$(0,0)$};
\draw (7,6) node[above]{$(700,600)$};
\draw (0,3) node[rotate=90,anchor=north,yshift=6mm]{$p=1013250$};
\draw (7,3) node[rotate=90,anchor=north,yshift=0mm]{$p=0$};
\draw (3.5,0) node[below]{$\bfu\cdot\bfn=0$};
\draw (3.5,6) node[above]{$\bfu\cdot\bfn=0$};
\draw[blue] (0,5) -- (7,5);
\draw[blue] (6.25,0) -- (6.25,6);
\draw (3.45,3.37) -- (3.75,3.37) -- (3.75,3.55) -- (3.45,3.55) -- cycle;
\end{tikzpicture}
%}
%\subfloat[Pressure solution on $\mathcal{M}_6^6$.]
%{\raisebox{6mm}{\includegraphics[width=0.45\textwidth]{fig/realistic_case/pressure_lr6-6}}
%\label{fig:benchmark4_pressure}}
\caption{Benchmark 4 (pressure problem). Problem description with boundary conditions. The red lines represent fractures, the blue lines represent the lines for which the solution is plotted along in Figs.~\ref{fig:benchmark4_pressure_y500} and \ref{fig:benchmark4_pressure_x625}, and the small black box represents the box plotted in Fig.~\ref{fig:benchmark4_mesh}.}
\label{fig:benchmark4_ill}
\end{figure}
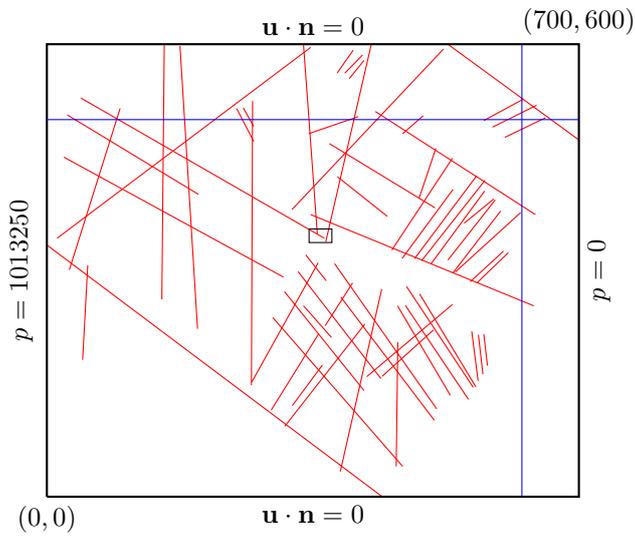

We consider next benchmark 4 in \cite{flemisch2018bfs}. The geometry represents a real set of fractures from an interpreted outcrop in the Sotra island, near Bergen in Norway, see \cite{flemisch2018bfs} for more details. For this case $\Omega=(0,700)\times(0,600)$ and $\Gamma$ is the union of 64 straight fractures in a complex pattern, see Fig.~\ref{fig:benchmark4_ill}. The matrix permeability is set to $\bfK=10^{-14}\bfI~\textup{m}^2$. All fractures have permeability $\kappa_\Gamma=10^{-8}~\textup{m}^2$, and aperture $w=10^{-2}~\textup{m}$. We apply homogeneous Neumann boundary conditions (no flow) on the top and bottom boundary faces, $p_\tD=1013250~\textup{Pa}$ on the left face, and $p_\tD=0~\textup{Pa}$ on the right face.

\begin{figure}[tbp]
	\centering
	\includegraphics[width=0.85\textwidth]{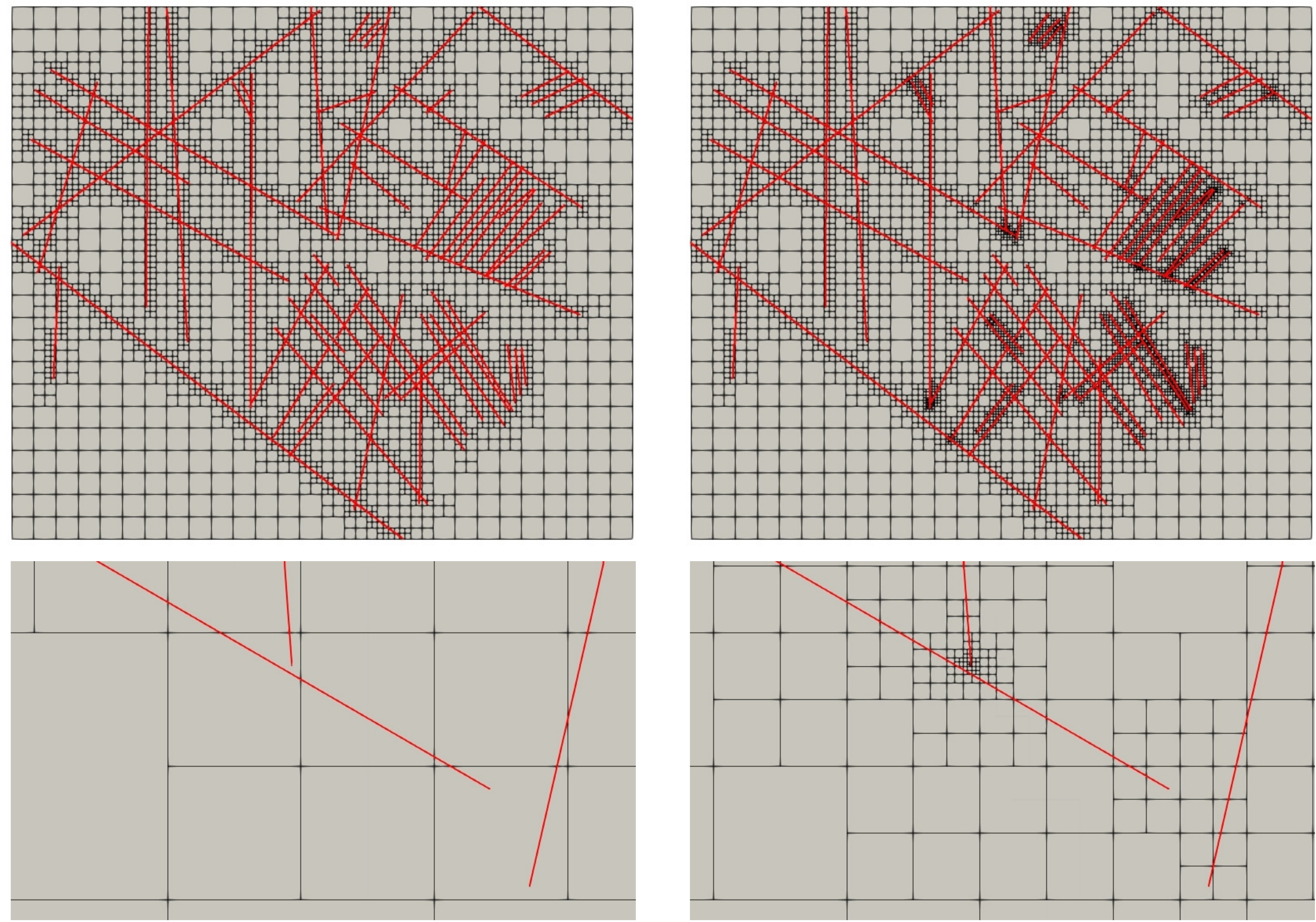}
%	\includegraphics[width=0.4\textwidth]{fig/realistic_case/mesh_2-2}
%	\hspace{3mm}
%	\includegraphics[width=0.4\textwidth]{fig/realistic_case/mesh_2-2-res} \\
%	\vspace{1mm}
%	\hspace{0.05mm}
%	\includegraphics[width=0.395\textwidth]{fig/realistic_case/mesh_2-2_zoom}
%	\hspace{4.1mm}
	%\includegraphics[width=0.395\textwidth]{fig/realistic_case/mesh_2-2-res_zoom}
	\caption{Benchmark 4 (pressure problem). Computational meshes, $\mathcal{M}_2^2$ [left] and $\mathcal{M}_2^{2,\textup{r}}$ [right]. The top row displays the whole domain $\Omega$, while the bottom row displays the mesh on the small rectangle in the middle of Fig.~\ref{fig:benchmark4_ill}. }
	\label{fig:benchmark4_mesh}
\end{figure}

Denote by $\mathcal{M}_0$ the $7\times 6$ base mesh with $h=100$. Then let $\mathcal{M}_i^j$ be the mesh where $\mathcal{M}_0$ is first refined globally $i$ times, and then recursively refined locally around the fractures $j$ times. With this $h_\text{min} = 100\cdot2^{-(i+j)}$ and  $h_\text{max} = 100\cdot2^{-i}$. Some of the fractures are very close without intersecting. It is important to resolve this geometrical aspect so that each vertex patch only see one fracture unless the fractures are connected. We denote by $\mathcal{M}_i^{j,\textup{r}}$ the mesh where $\mathcal{M}_i^j$ is further locally refined to resolve close non-connected fractures. The DoFs and number of elements are reported in Table~\ref{tab:benchmark4_dofs}, along with the corresponding numbers for the methods considered in \cite{flemisch2018bfs}. Fig.~\ref{fig:benchmark4_mesh} displays two of the meshes and illustrates the effect of resolving the geometry.
%Our finest mesh is $\mathcal{M}_6^6$ with 2~421~489 elements, and the corresponding pressure approximation is shown in Fig.~\ref{fig:benchmark4_pressure}.

\begin{figure}[tbp]
	\centering
\subfloat[$\mathcal{M}_2^{2}$.]
{\includegraphics[width=0.4\textwidth]{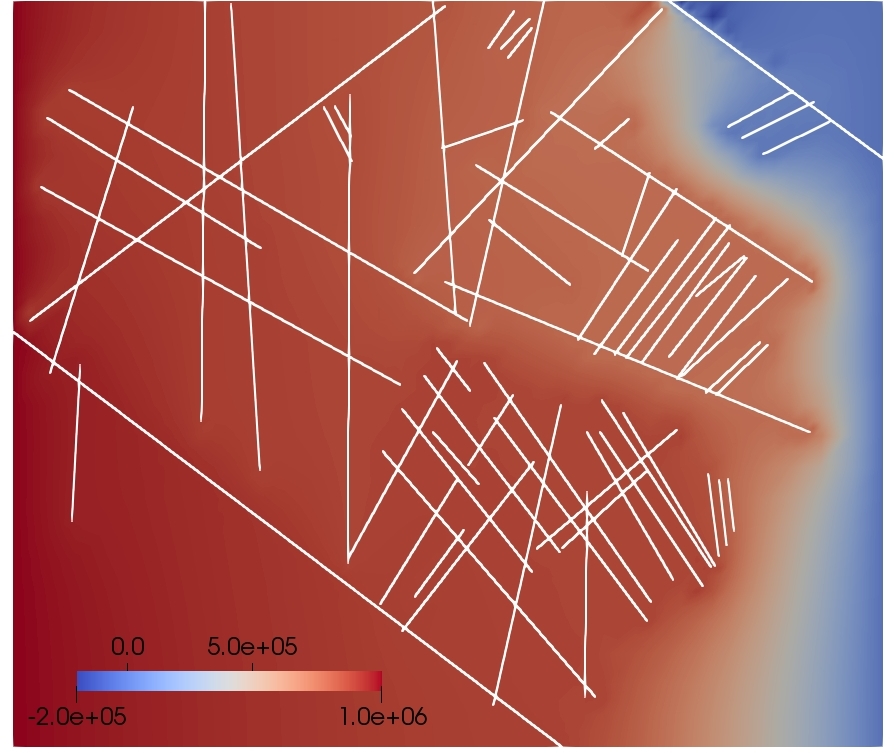}}
\subfloat[$\mathcal{M}_2^{2,\textup{r}}$.]
{\includegraphics[width=0.4\textwidth]{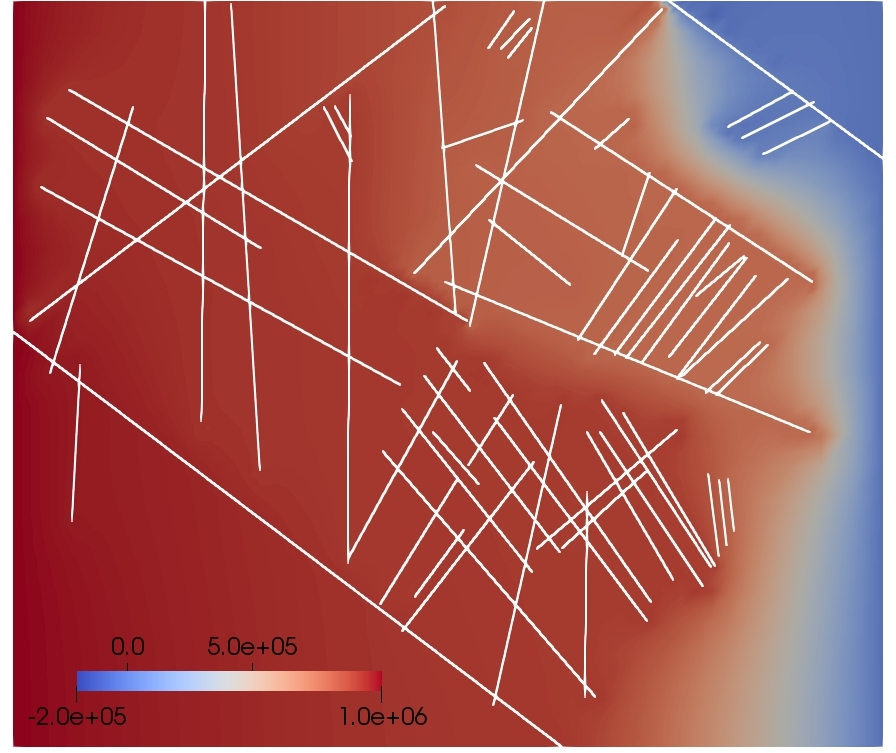}} \\
\subfloat[$\mathcal{M}_2^{4}$.]
{\includegraphics[width=0.4\textwidth]{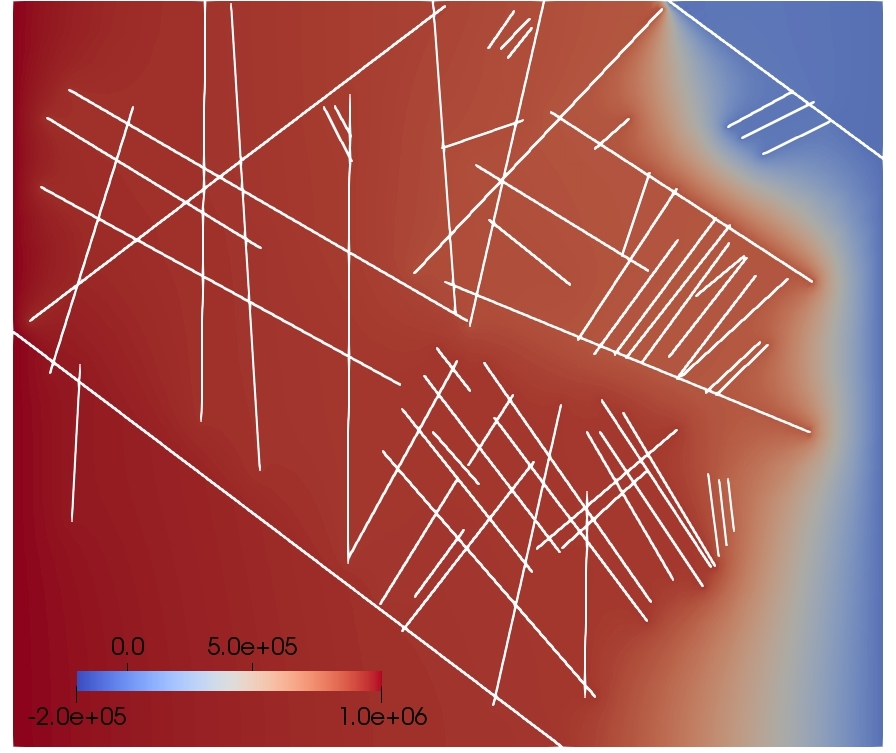}}
\subfloat[$\mathcal{M}_2^{4,\textup{r}}$.]
{\includegraphics[width=0.4\textwidth]{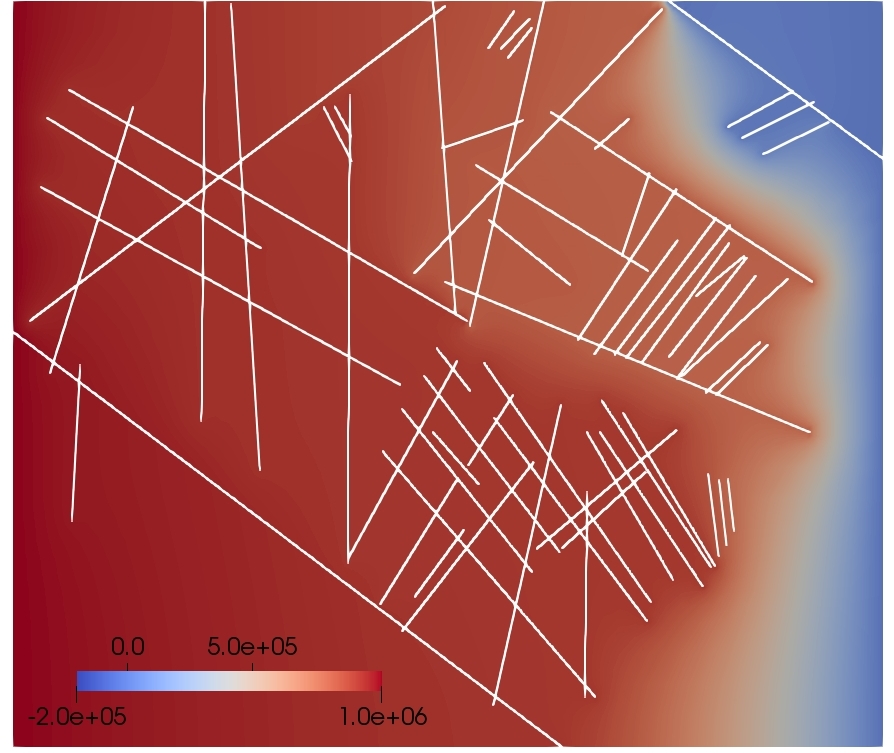}}
\caption{Benchmark 4 (pressure problem). Pressure approximations.}
\label{fig:benchmark4_pressure_sol}
\end{figure}

\begin{figure}[tbp]
	\centering
	\includegraphics[width=0.49\textwidth]{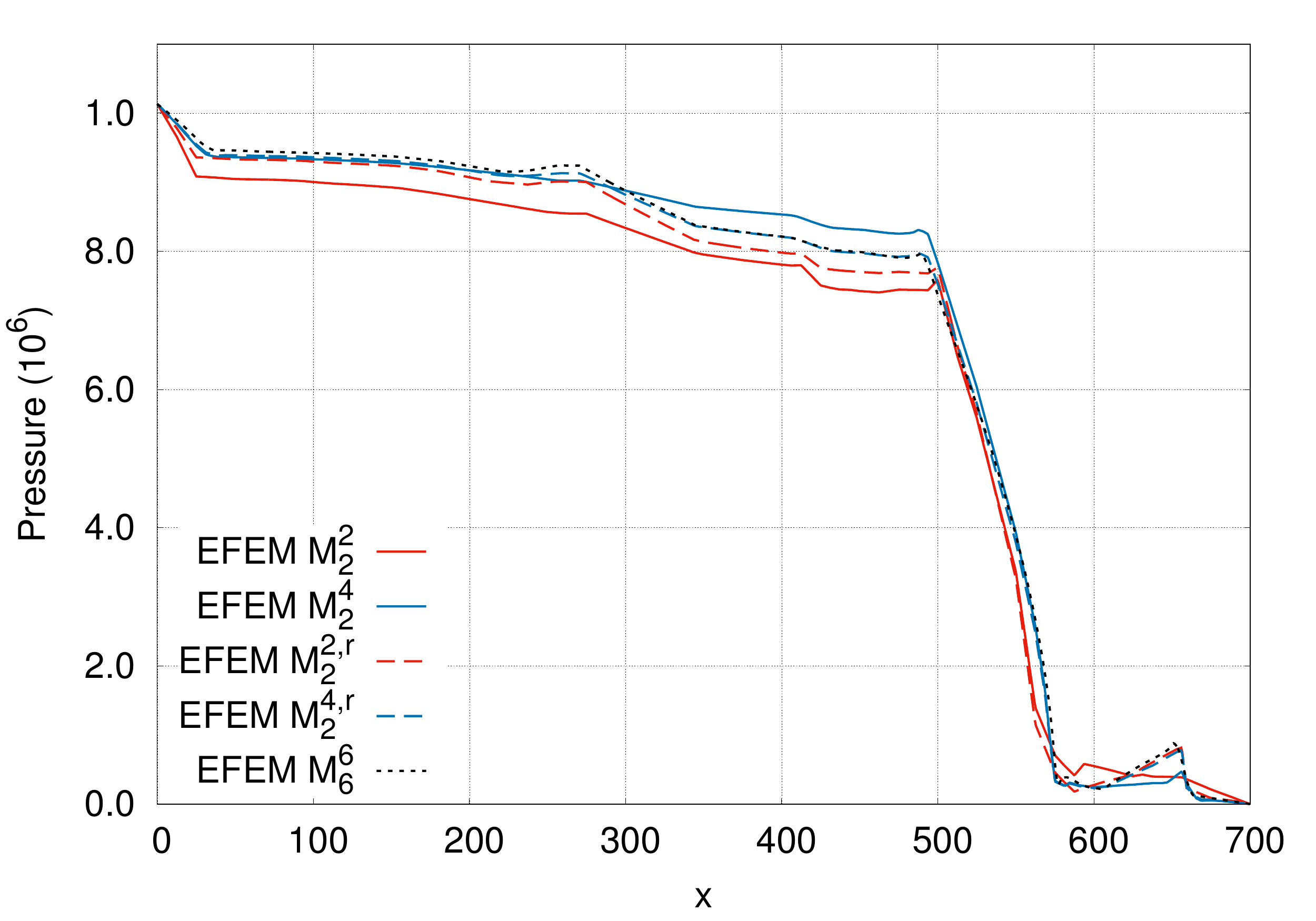}
	\includegraphics[width=0.49\textwidth]{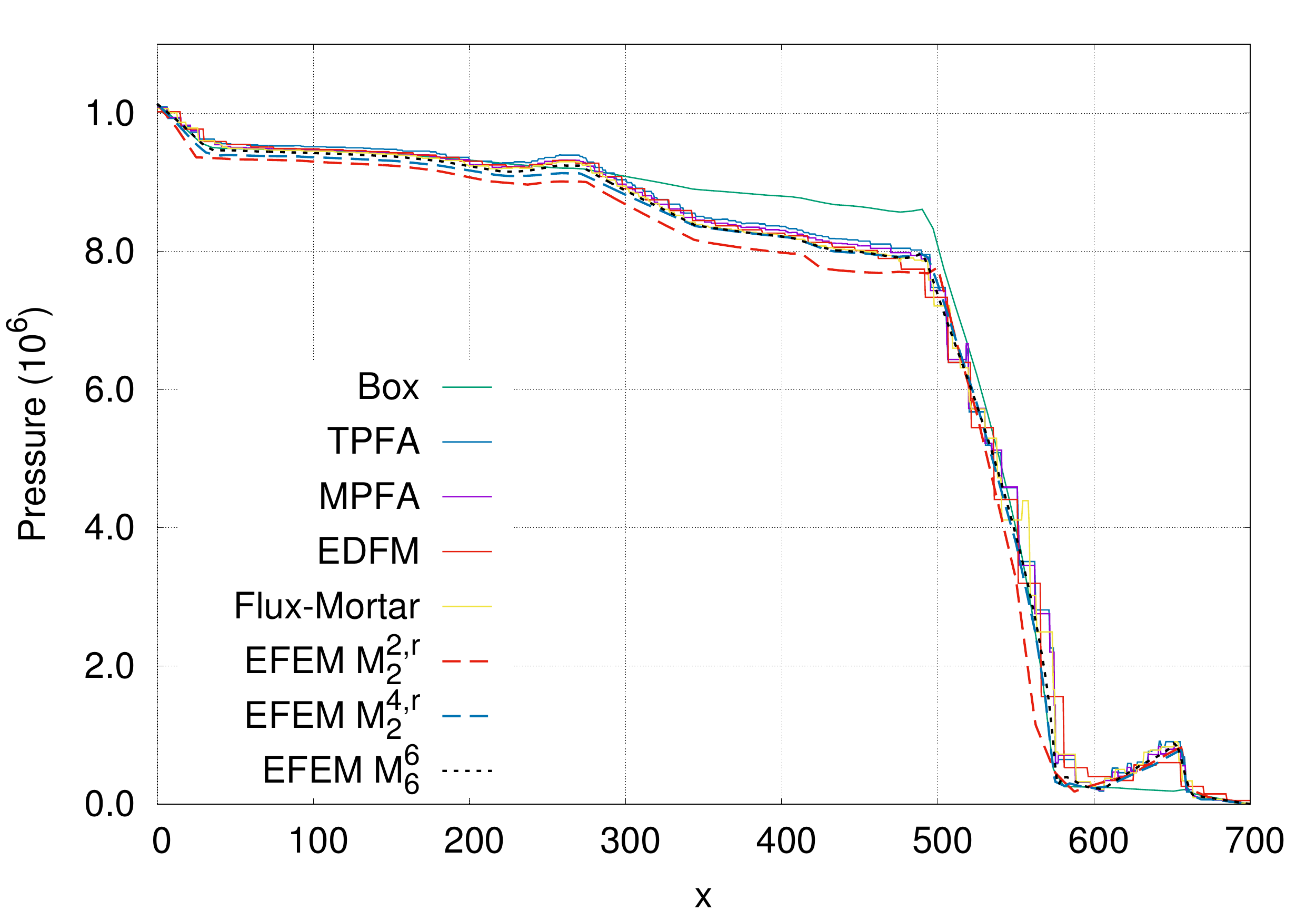}
	\caption{Benchmark 4. Pressure solution along the line $y=500$. The solution on $\mathcal{M}_6^6$ is considered as a reference solution.}
	\label{fig:benchmark4_pressure_y500}
\end{figure}

\begin{figure}[tbp]
	\centering
	\includegraphics[width=0.49\textwidth]{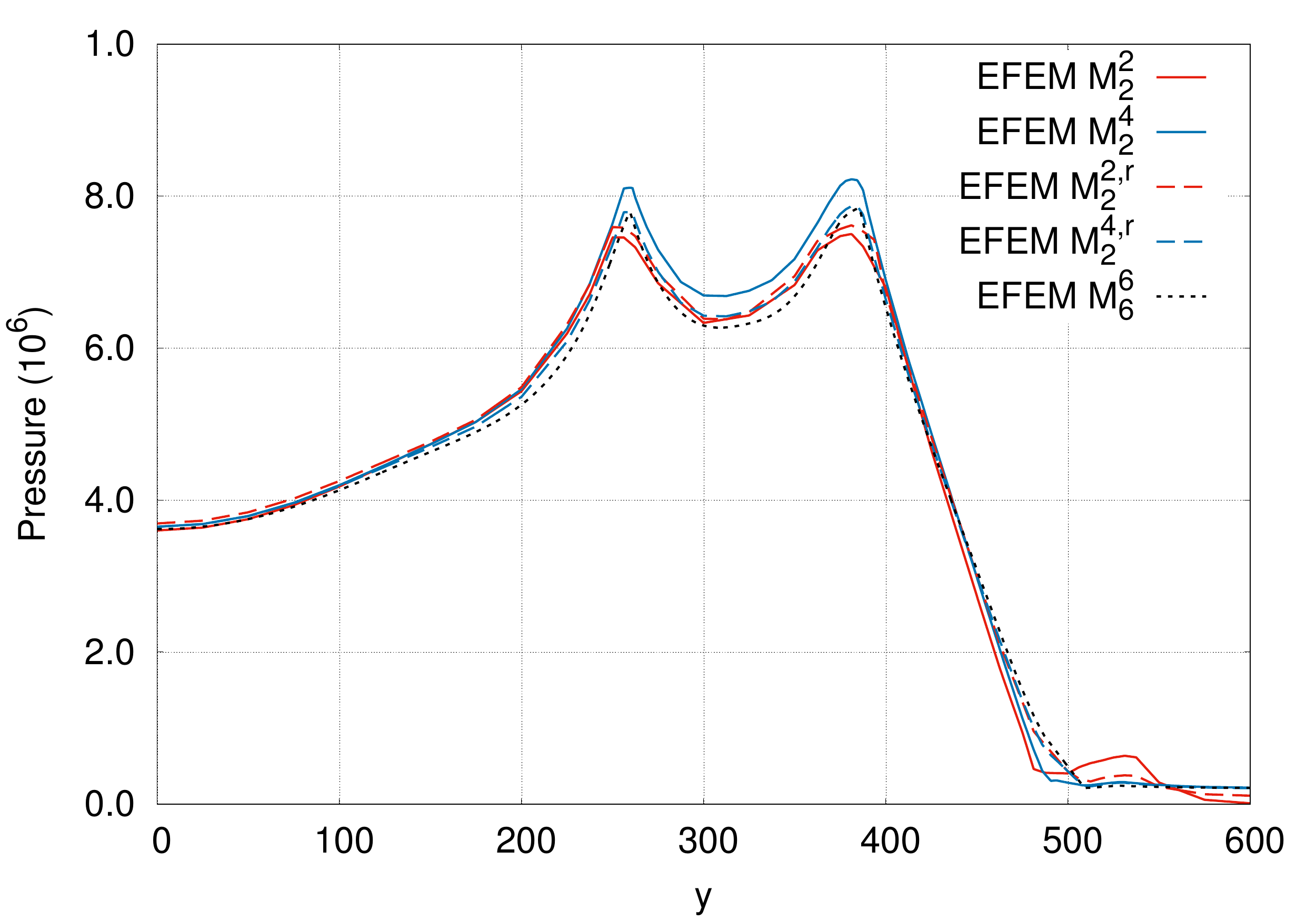}
	\includegraphics[width=0.49\textwidth]{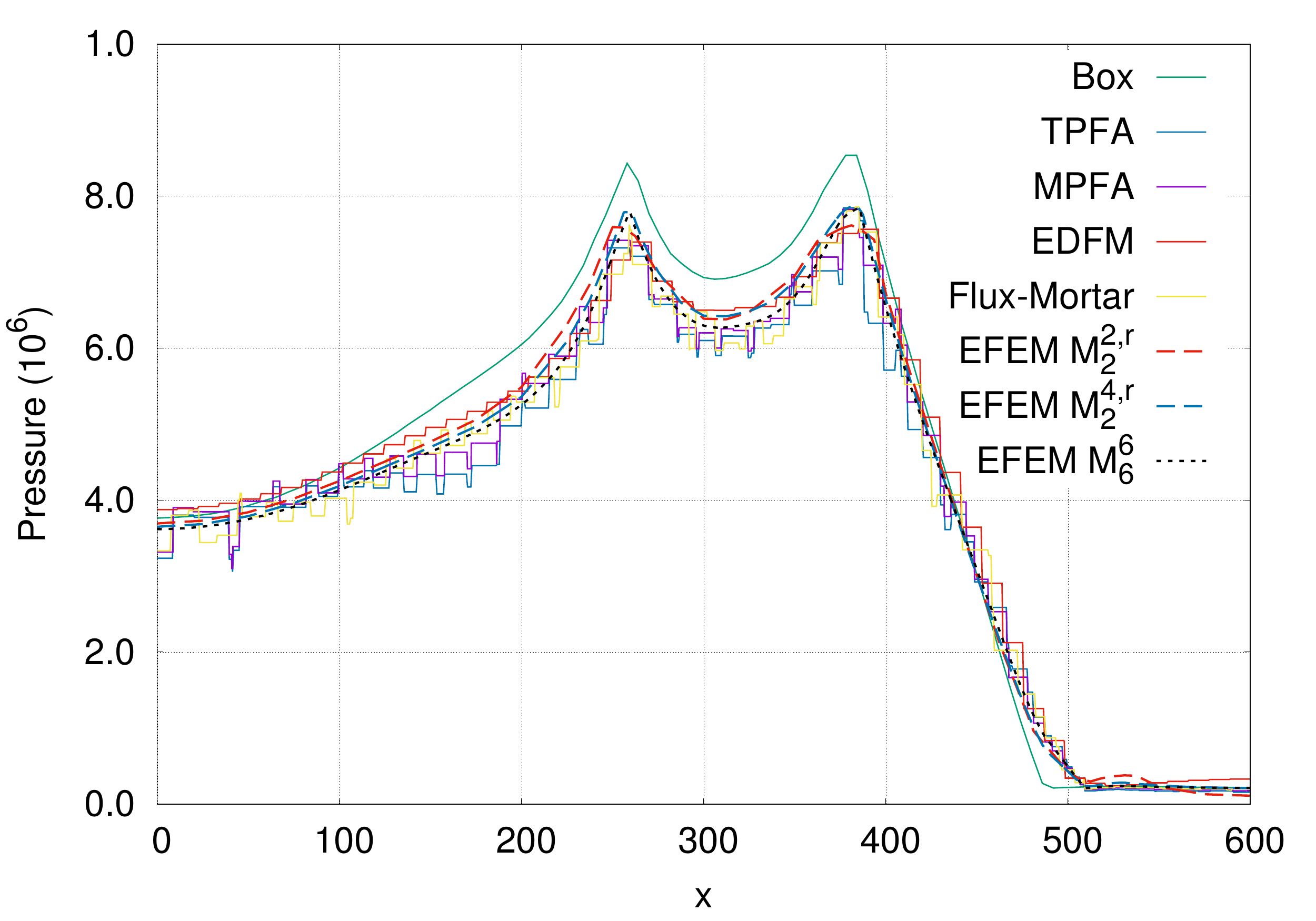}
	\caption{Benchmark 4 (pressure problem). Pressure solution along the line $x=625$. The solution on $\mathcal{M}_6^6$ is considered as a reference solution.}
	\label{fig:benchmark4_pressure_x625}
\end{figure}

Fig.~\ref{fig:benchmark4_pressure_sol} shows the pressure approximation on four different meshes. Observe that we have some unstabilities in the solution in the upper right corner that vanish as we refine. This is due to the fact that one of the fractures cuts the corner and that we have zero Dirichlet condition at the right end of this fracture.
Figs.~\ref{fig:benchmark4_pressure_y500} and \ref{fig:benchmark4_pressure_x625} plots the pressure approximation along the lines $y=500$ and $x=625$. Results on the highly refined mesh $\mathcal{M}_6^6$, where all geometry is sufficiently resolved, are included as a reference solution. We clearly see the effect of resolving the geometry, and we observe similar results to the ones reported in \cite{flemisch2018bfs}.

\begin{filecontents*}{complex_network.csv}
FID,START_X,START_Y,END_X,END_Y
1,269.611206,152.05243,356.9240112,310.14123
2,249.5117187,514.990780001,272.218872,470.97082
3,258.3590698,515.574580001,271.9851684,490.9682
4,270.6622924,524.702640001,269.1347046,147.78143
5,355.8302002,348.479800001,337.5810733205,600
6,366.9730835,338.132990001,426.9185141723,600
7,198.237915,222.724420001,175.1561889,597.603030001
8,151.2785034,261.724610001,154.4623059774,600
9,29.5026855,300.724610001,96.3599853,514.82739
10,386.0808105,33.3621800002,440.585083,275.191830001
11,459.6350708,40.2413900001,461.751709,204.812620001
12,297.180603,237.62103,468.1018066,40.2413900001
13,312.5264892,272.01678,417.3016967,140.7832
14,330.5181884,298.47522,439.5266723,156.6582
15,340.5723877,320.70019,367.5598755,286.304380001
16,492.9725952,312.762820001,576.5811157,419.6546
17,505.6726684,309.05859,576.0520019,405.367190001
18,537.4227905,297.94598,623.3187866,376.68463
19,322.5338745,380.76941,521.8778076,593.552180001
20,344.9320678,481.56122,409.8867798,503.959410001
21,371.8098755,468.12219,510.6787109,383.009210001
22,432.2849731,510.678830001,642.8280029,374.04999
23,527.528634971,600,700,473.015615092
24,0,333.73321,441.2443847,0
25,13.4389038,342.692380001,347.171875,595.791990001
26,22.3981933,450.203790001,311.3347778,291.176630001
27,26.8778076,506.199220001,199.343811,400.92779
28,44.7963867,528.597410001,365.0905151,342.692380001
29,378.5294189,309.095210001,512.918518,116.470640001
30,461.4027099,253.099610001,530.8370971,134.38922
31,347.171875,374.04999,640.5881958,253.099610001
32,490.5203857,268.77844,564.4343872,145.58844
33,47.0361938,181.425410001,53.7556152,306.85541
34,382.4152832,424.151000001,447.8997192,371.76343
35,587.9967651,394.78222,549.1029663,362.635190001
36,589.9812011,393.59161,527.6716919,313.8194
37,597.125,378.90722,533.6248169,295.960200001
38,533.6248169,448.75738,453.8527832,326.91638
39,511.7966919,461.85419,489.5715942,395.17901
40,565.3748779,425.34161,483.6184692,315.40698
41,534.4185791,407.482240001,467.3466186,315.803830001
42,627.2874756,527.3388,574.8999023,498.763610001
43,644.3532104,519.00439,586.4093017,490.03241
44,655.8626098,502.335630001,602.6812133,476.53863
45,415.355896,585.679380001,391.9401855,561.47003
46,417.3402099,578.535580001,397.8933105,554.326230001
47,403.0526733,592.029420001,382.0183105,561.86682
48,495.1278686,505.113580001,468.1403198,481.30121
49,533.6248169,254.84381,420.9121093,159.196590001
50,508.6217041,221.10943,441.152771,159.59363
51,418.5308838,229.04681,312.961914,93.3154300004
52,362.5714111,174.6748,322.883789,120.69983
53,357.8088989,216.3468,295.102478,114.74658
54,402.2589111,283.41882,366.1433105,226.66559
55,337.5681762,253.256220001,374.4776001,211.18744
56,386.7808838,264.765620001,509.8123169,101.25281
57,473.2996826,278.65643,561.0092163,144.909240001
58,471.7122192,253.653200001,554.6593017,129.034240001
59,559.0249023,219.125,567.3593139,153.64044
60,567.7561035,214.759400001,573.7092895,162.37182
61,574.8999023,215.553040001,579.6624145,173.88104
62,557.0404663,285.006410001,600.6968994,325.48761
63,565.3748779,283.022030001,607.0468139,323.503230001
\end{filecontents*}
\DTLloaddb[noheader=false]{coordinates}{complex_network.csv}

\begin{table}[tbp]
	\caption{Benchmark 4 (pressure problem). DoF, number of elements, matrix density (nnz denotes number of non-zero entries in the system matrix) and matrix condition number for EFEM on different meshes compared to the corresponding methods in \cite{flemisch2018bfs}. \emph{No.\ of elements} for the methods in \cite{flemisch2018bfs} are listed as matrix elements plus fracture elements.}
	\label{tab:benchmark4_dofs}
	\begin{center}
		\begin{tabular}{lrrrr}
			\hline
			\textbf{Method} & $\ndof$ & \textbf{No.~of elements} & \textbf{nnz}/$\ndof^2$ & $\Vert\cdot\Vert_2\textbf{-cond}$ \\
			\hline
			EFEM $\mathcal{M}_2^2$     &  5349 &  4623 & 1.3e-3 & 3.5e5 \\
			EFEM $\mathcal{M}_2^{2,\textup{r}}$ &  9185 &  7629 & 7.2e-4 & 9.3e6 \\
			EFEM $\mathcal{M}_2^4$     & 33337 & 27666 & 2.0e-4 & 1.5e7$^*$ \\
			EFEM $\mathcal{M}_2^{4,\textup{r}}$ & 33924 & 28104 & 1.9e-4 & 6.4e7$^*$ \\
			\hline
			Box			&  5563 & 10807 + 1386 & 1.2e-3 & 9.3e5 \\
			TPFA		& 8481  & 7614 + 867   & 4.9e-4 & 5.3e6 \\
			MPFA	    & 8588  & 7614 + 867   & 1.6e-3 & 4.9e6 \\
			EDFM		& 3599  & 2491 + 1108  & 1.4e-3 & 4.7e6 \\
			Flux-Mortar	& 25258 & 8319 + 1317  & 2.0e-4 & 2.2e17 \\
			\hline
			\\
			\multicolumn{5}{l}{\small $^*$ Estimate of the 1-norm condition number based on MATLABs \texttt{condest} command.}
		\end{tabular}
	\end{center}
\end{table}

%\clearpage

\subsection{Coupled problems}
\label{sec:results_coupled}

We now couple the benchmark cases with the transport problem and solve for concentration. The coupled problem was not considered in \cite{flemisch2018bfs}.

\subsubsection{Benchmark 1: Regular fracture network}

For the transport problem we set $c_0=0$, apply a boundary concentration $c_\tB=1$ on the inflow boundary ($x=0$), and let $T=0.5$. We consider the same computational meshes as for the pressure problem. In addition, we construct a reference mesh with 182674 elements where the fractures are fully resolved, i.e., $h=w$ at the fractures and $h\approx 16w$ away from the fractures. A reference solution is then obtained by a standard FV method on a continuous fracture model (no lower dimensional embedding), i.e., Eq.~\eqref{eq:FV-IE_pp_matrix} with $\Kh^\tM=\Kh$ and $\bfK=\kappa_\Gamma\bfI$ in the elements in the fracture domain. Time steps for the reference solution is $\Delta t=\powten{1}{-5}$.

\begin{figure}[tbp]
	\centering
	\subfloat[Fracture 1 ($y=0.5$)]
	{\includegraphics[width=0.33\textwidth]{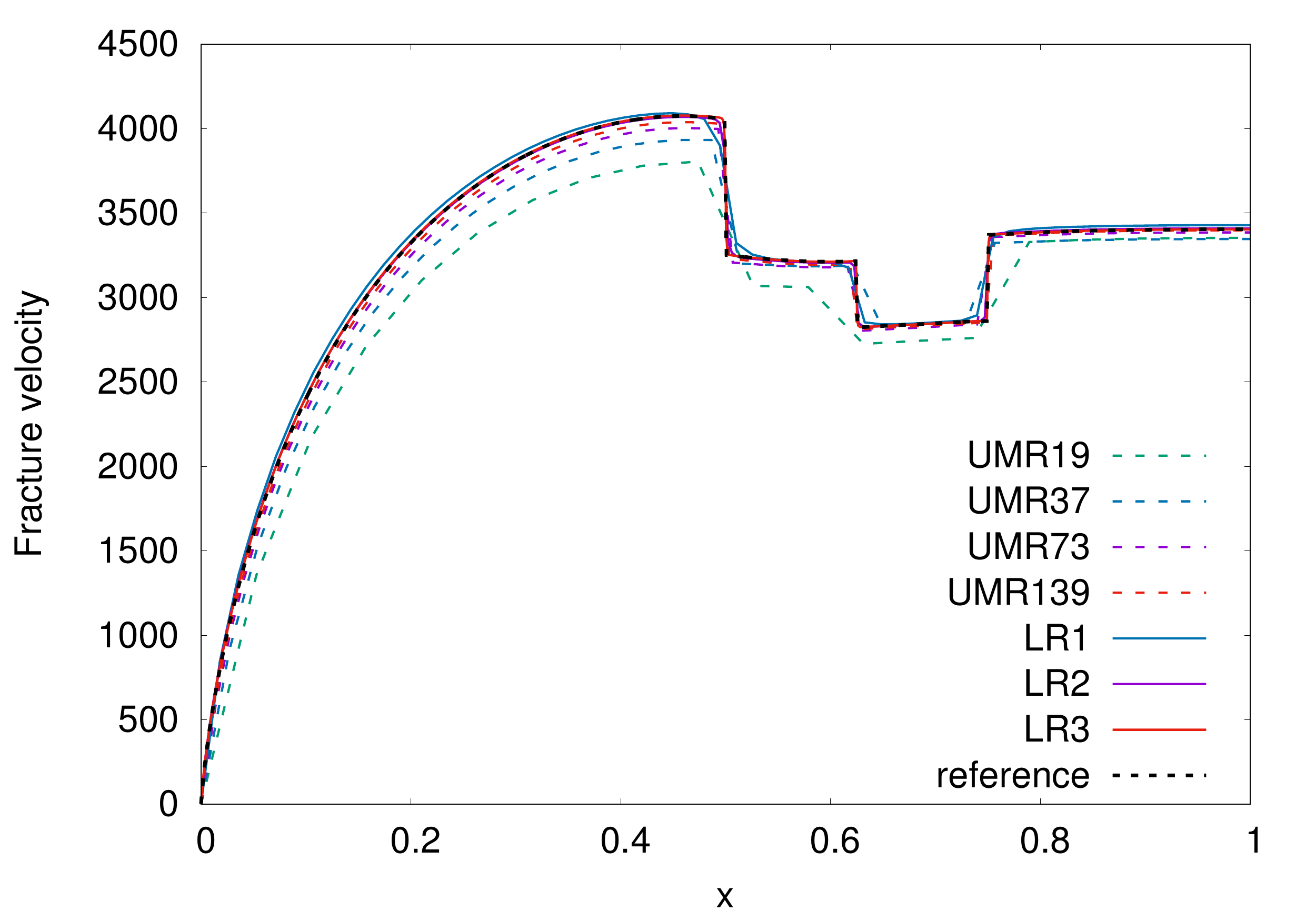}}
	\subfloat[Fracture 3 ($y=0.75$)]
	{\includegraphics[width=0.33\textwidth]{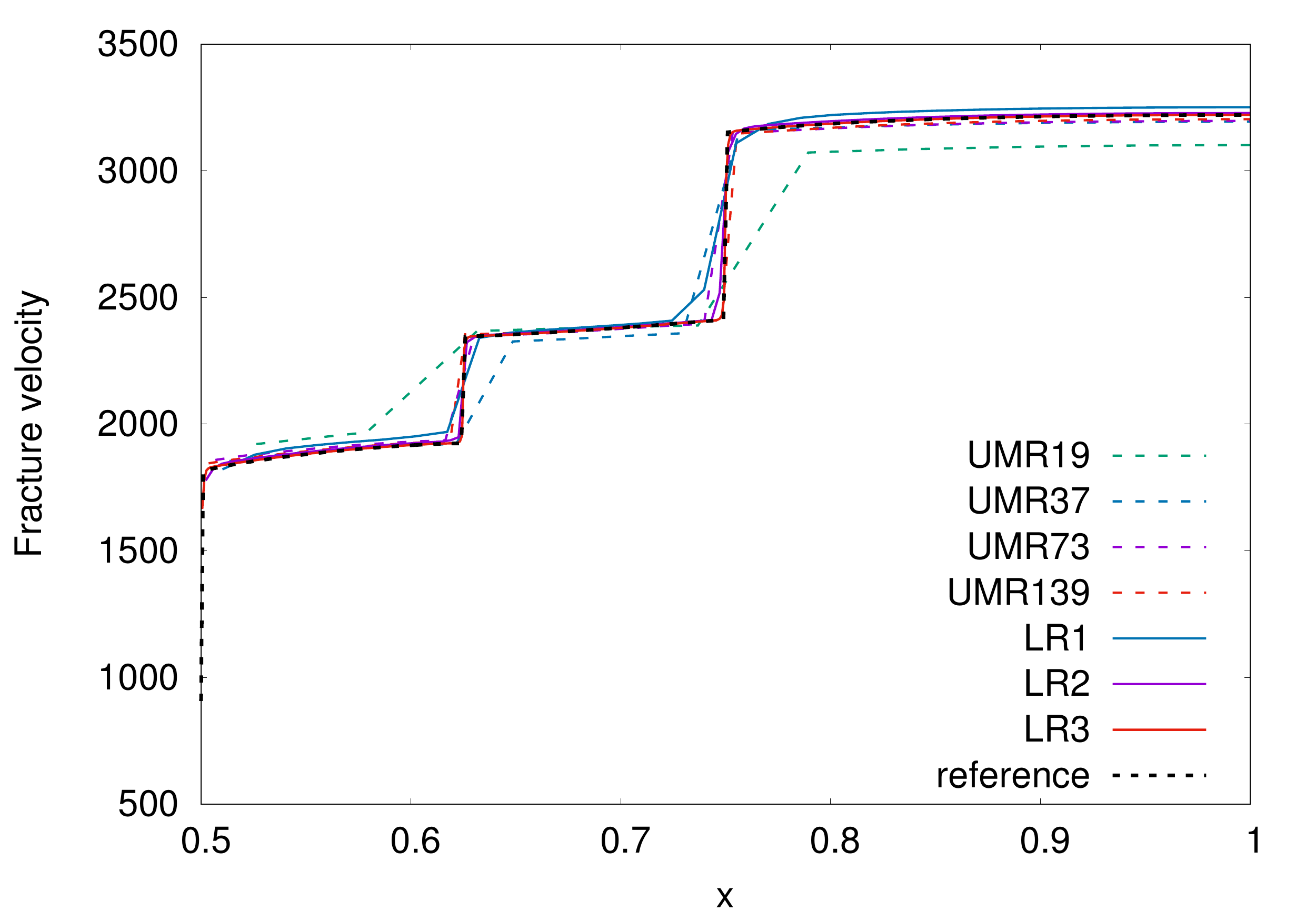}}
	\subfloat[Fracture 5 ($y=0.625$)]
	{\includegraphics[width=0.33\textwidth]{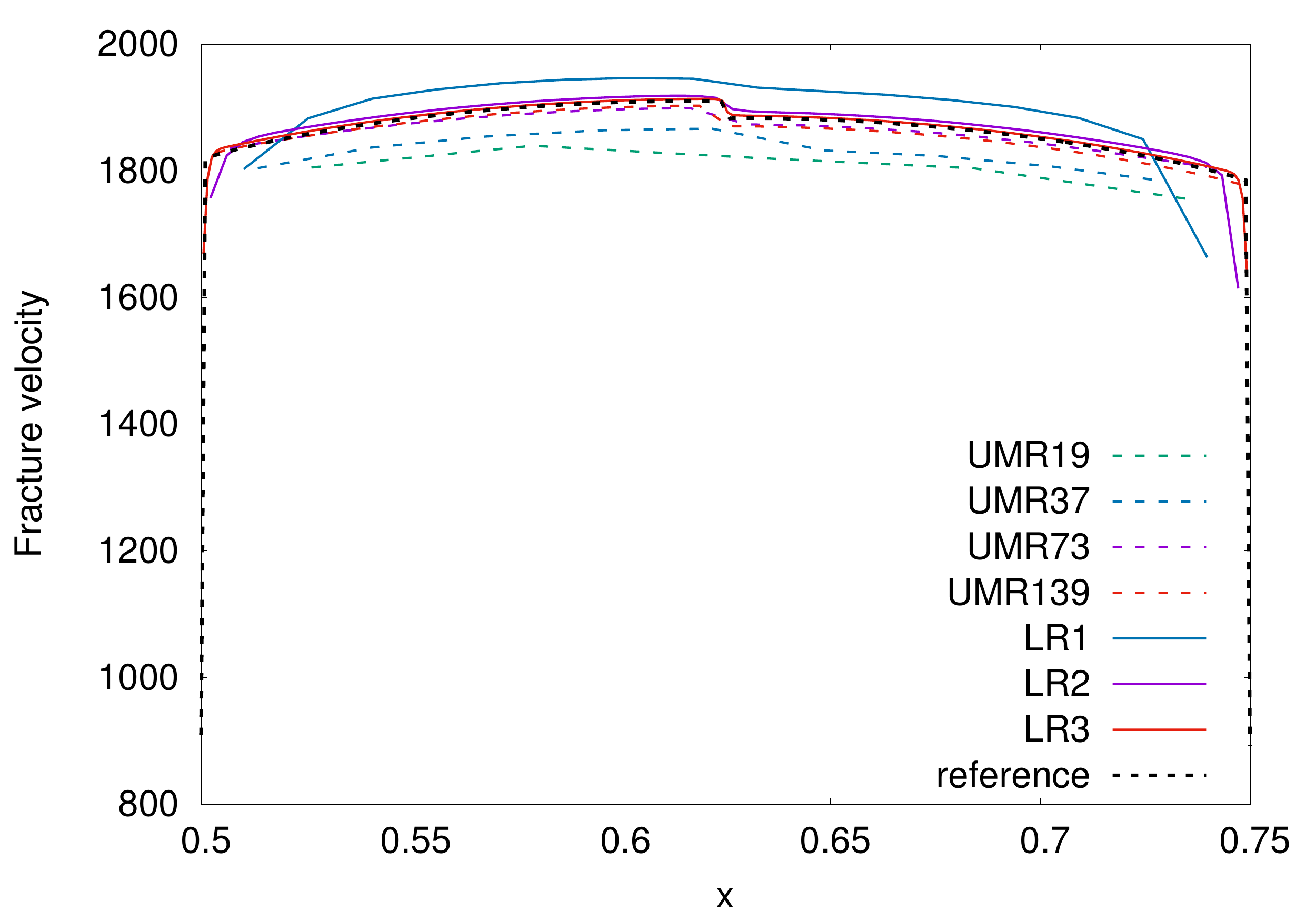}} \\
	\subfloat[Fracture 2 ($x=0.5$)]
	{\includegraphics[width=0.33\textwidth]{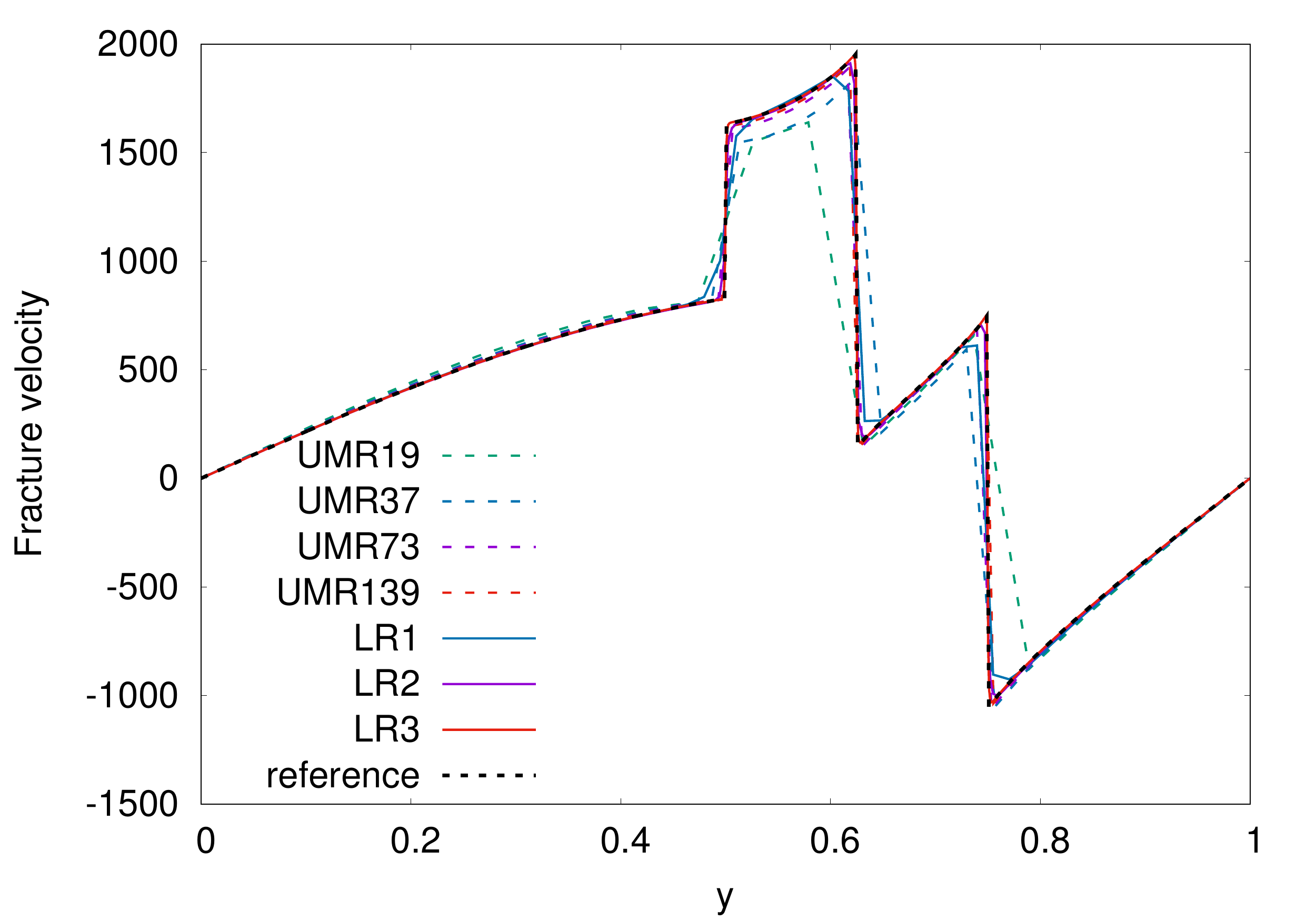}}
	\subfloat[Fracture 4 ($x=0.75$)]
	{\includegraphics[width=0.33\textwidth]{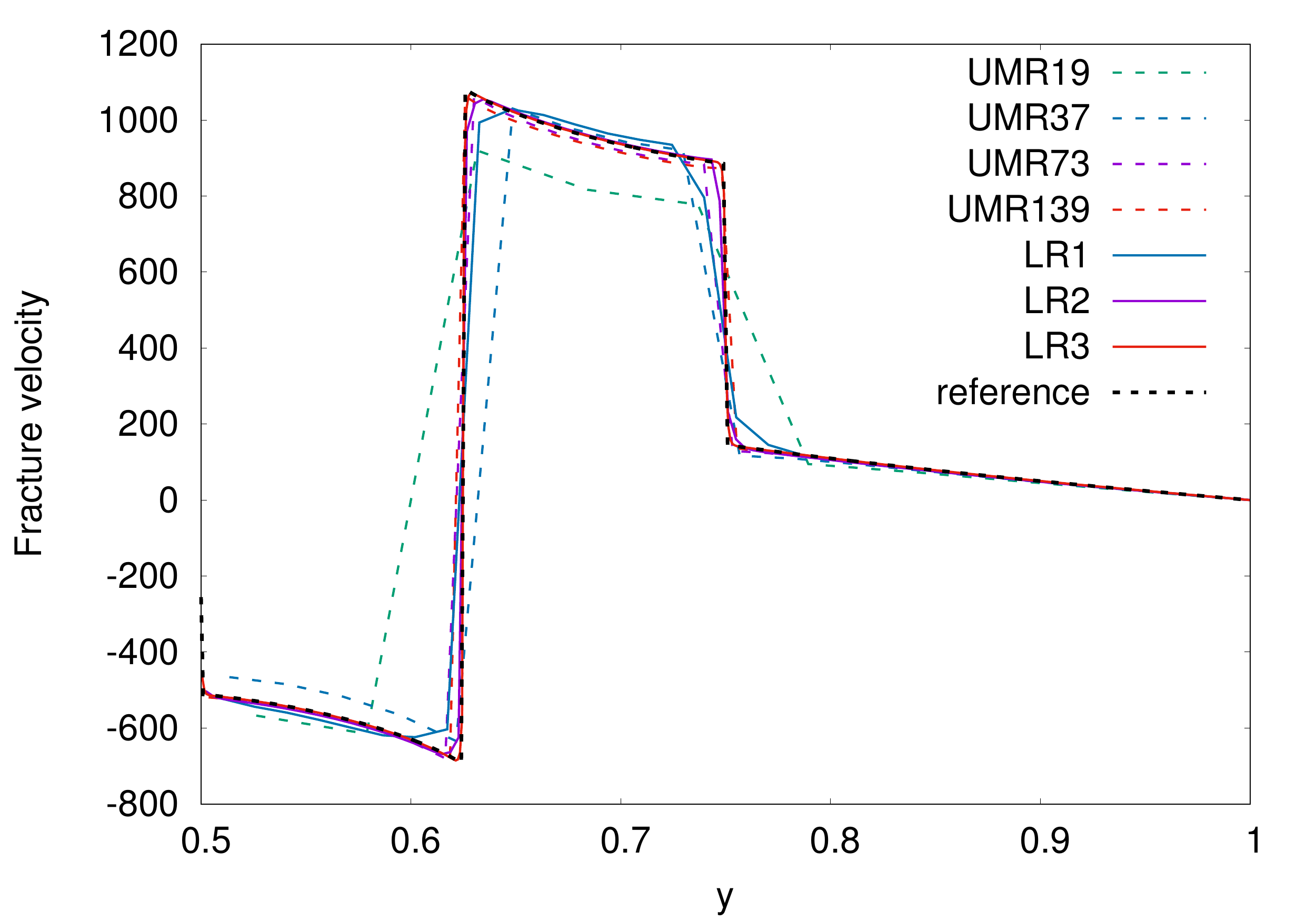}}
	\subfloat[Fracture 6 ($x=0.625$)]
	{\includegraphics[width=0.33\textwidth]{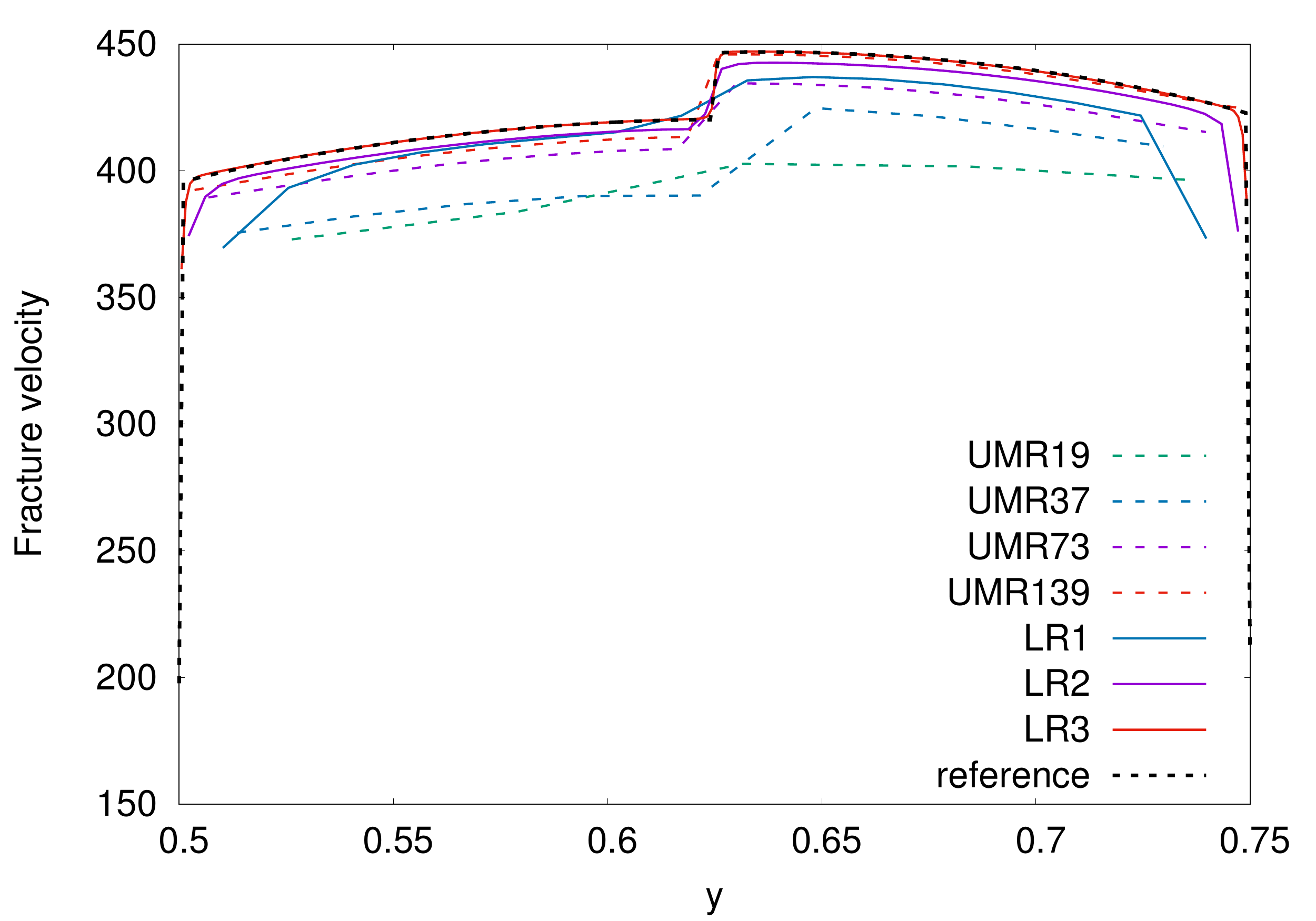}}
	\caption{Benchmark 1 (coupled problem): Fracture velocities for the different meshes.}
	\label{fig:benchmark_velocities}
\end{figure}

Fig.~\ref{fig:benchmark_velocities} shows the velocity approximations in each of the six fractures for the different meshes. We see a very good agreement, in particular for the LR meshes. This is as expected since the velocity is derived from the pressure solution, which was shown earlier to have higher accuracy when refining around the fractures.

For the transport problem, we define a quantity of interest, $\qoi$, as the flux of concentration out of the two fractures on the right boundary face,
\begin{subequations}
\begin{align}
\qoi_1(t) &= (\bfu_\Gamma\cdot\bfn c)\vert_{\bfx=(1.0,0.5)}, \\
\qoi_2(t) &= (\bfu_\Gamma\cdot\bfn c)\vert_{\bfx=(1.0,0.75)}.
\end{align}
\end{subequations}
For the reference solution, these quantities are calculated as
\begin{subequations}
\begin{align}
\qoi^{\textup{ref}}_1(t) &= \int_{0.5-\tfrac{w}{2}}^{0.5+\tfrac{w}{2}} (\bfu\cdot\bfn c)\vert_{x=1} ~\dy, \\
\qoi^{\textup{ref}}_2(t) &= \int_{0.75-\tfrac{w}{2}}^{0.75+\tfrac{w}{2}} (\bfu\cdot\bfn c)\vert_{x=1} ~\dy.
\end{align}
\end{subequations}

\begin{figure}[tbp]
	\centering
	%\subfloat[Reference solution, $t=0.01$.]
	%{\includegraphics[width=0.33\textwidth]{fig/regular_fractures/conc_reference_t001}}
	\subfloat[Reference solution, $t=0.1$.]
	{\includegraphics[width=0.36\textwidth]{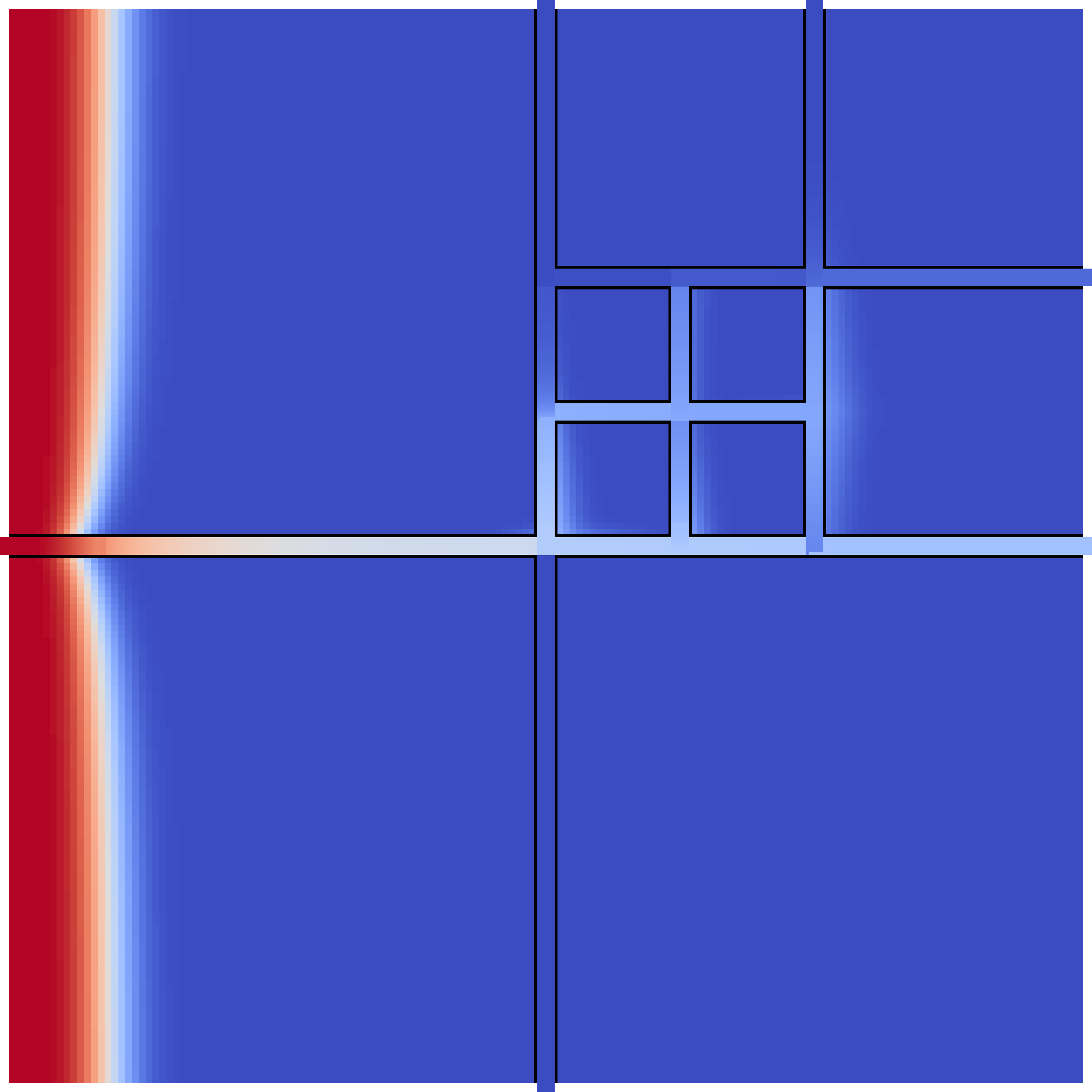}}
	\hspace{3mm}
	\subfloat[Reference solution, $t=0.5$.]
	{\includegraphics[width=0.36\textwidth]{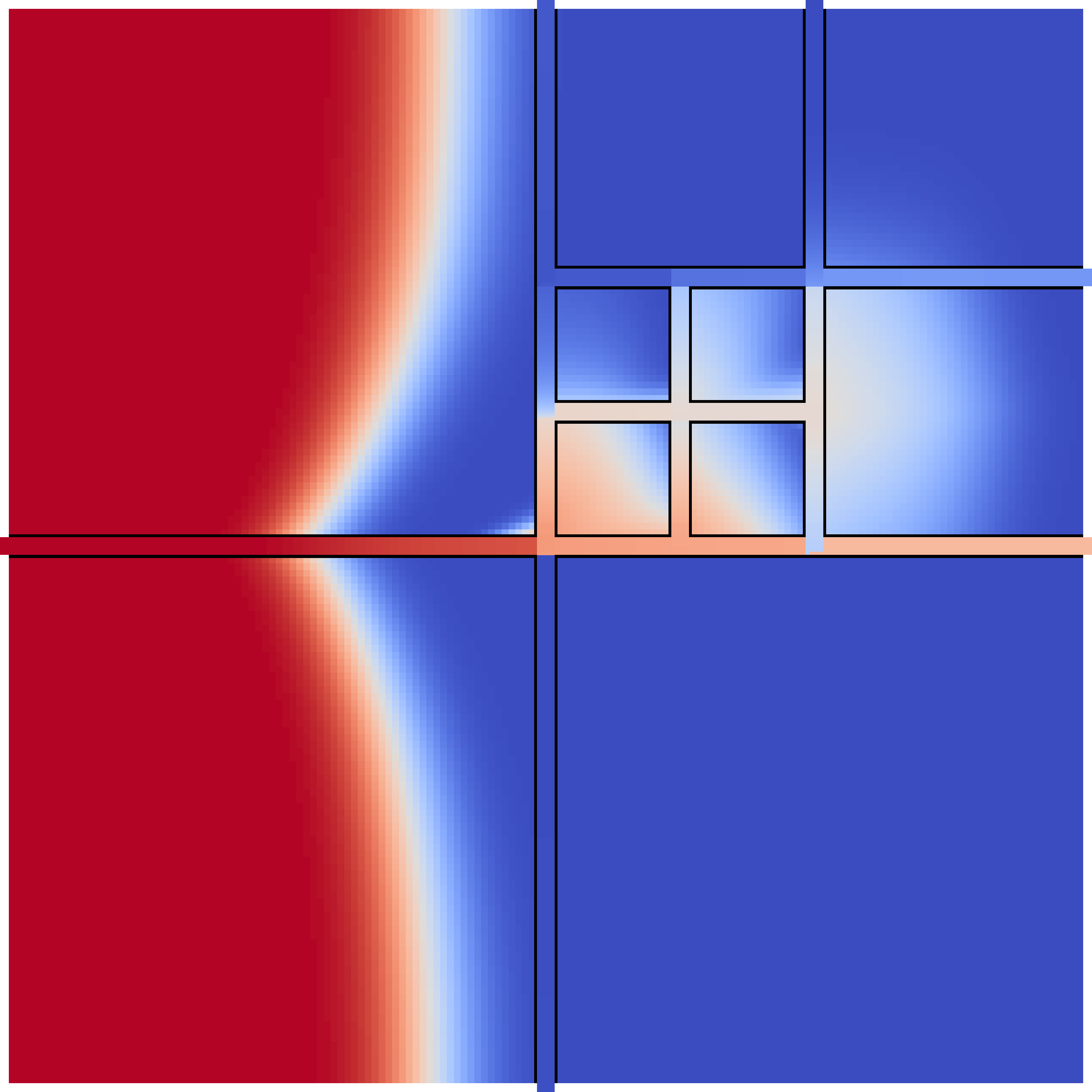}} \\
	%\subfloat[UMR139, $t=0.01$.]
	%{\includegraphics[width=0.33\textwidth]{fig/regular_fractures/conc_umr139_t001}}
	\subfloat[UMR139, $t=0.1$.]
	{\includegraphics[width=0.36\textwidth]{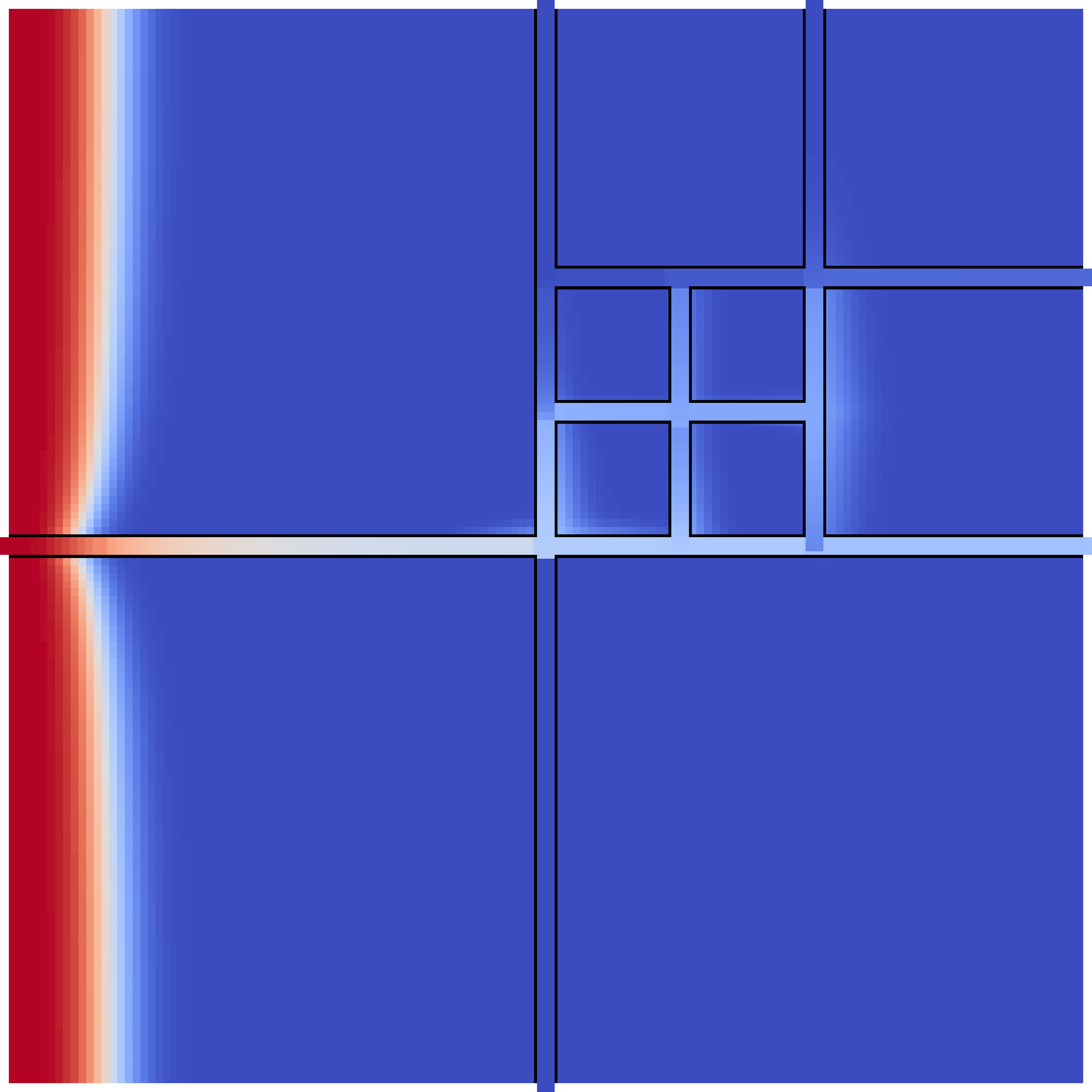}}
	\hspace{3mm}
	\subfloat[UMR139, $t=0.5$.]
	{\includegraphics[width=0.36\textwidth]{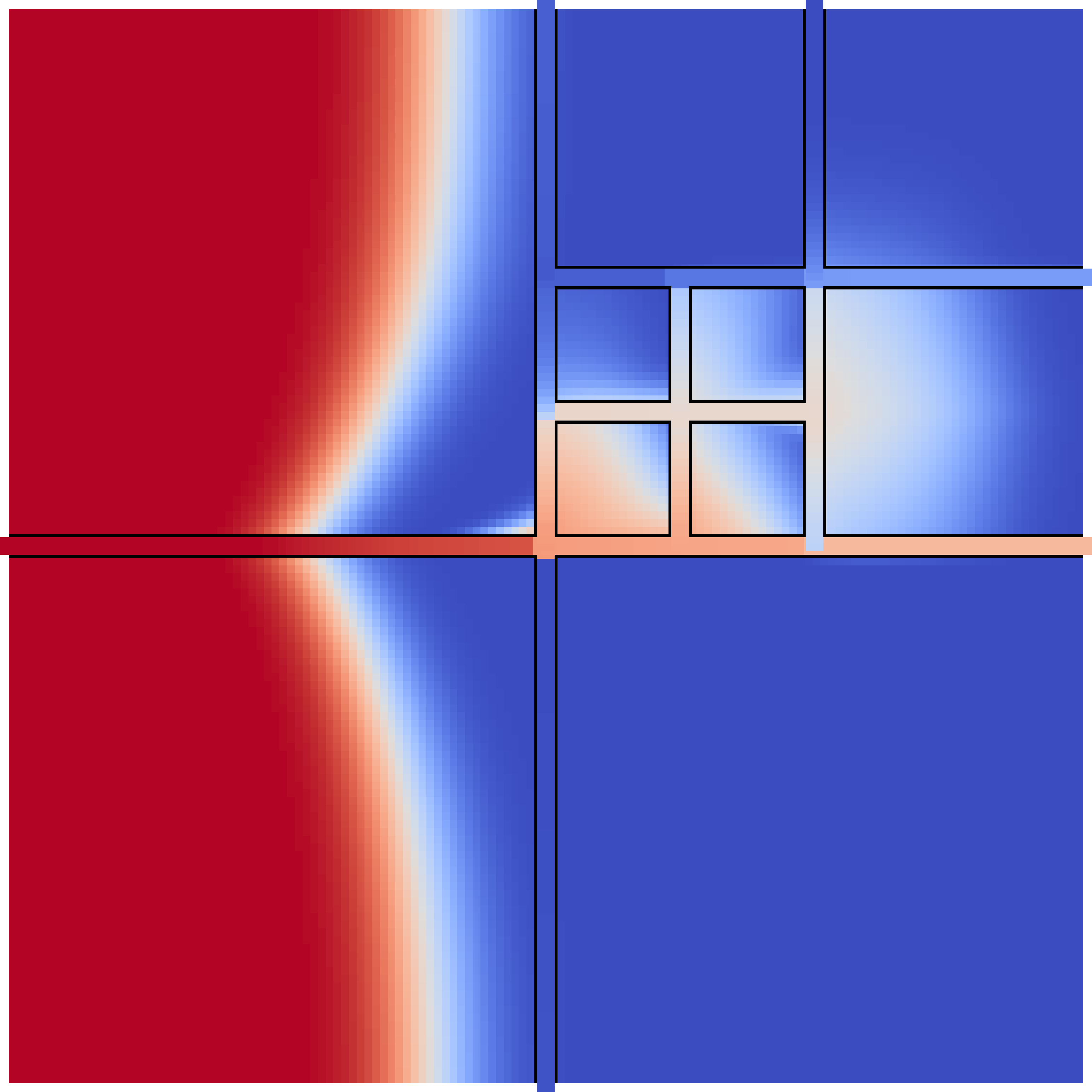}} \\
	%\subfloat[LR3, $t=0.01$.]
	%{\includegraphics[width=0.33\textwidth]{fig/regular_fractures/conc_lr3_t001}}
	\subfloat[LR3, $t=0.1$.]
	{\includegraphics[width=0.36\textwidth]{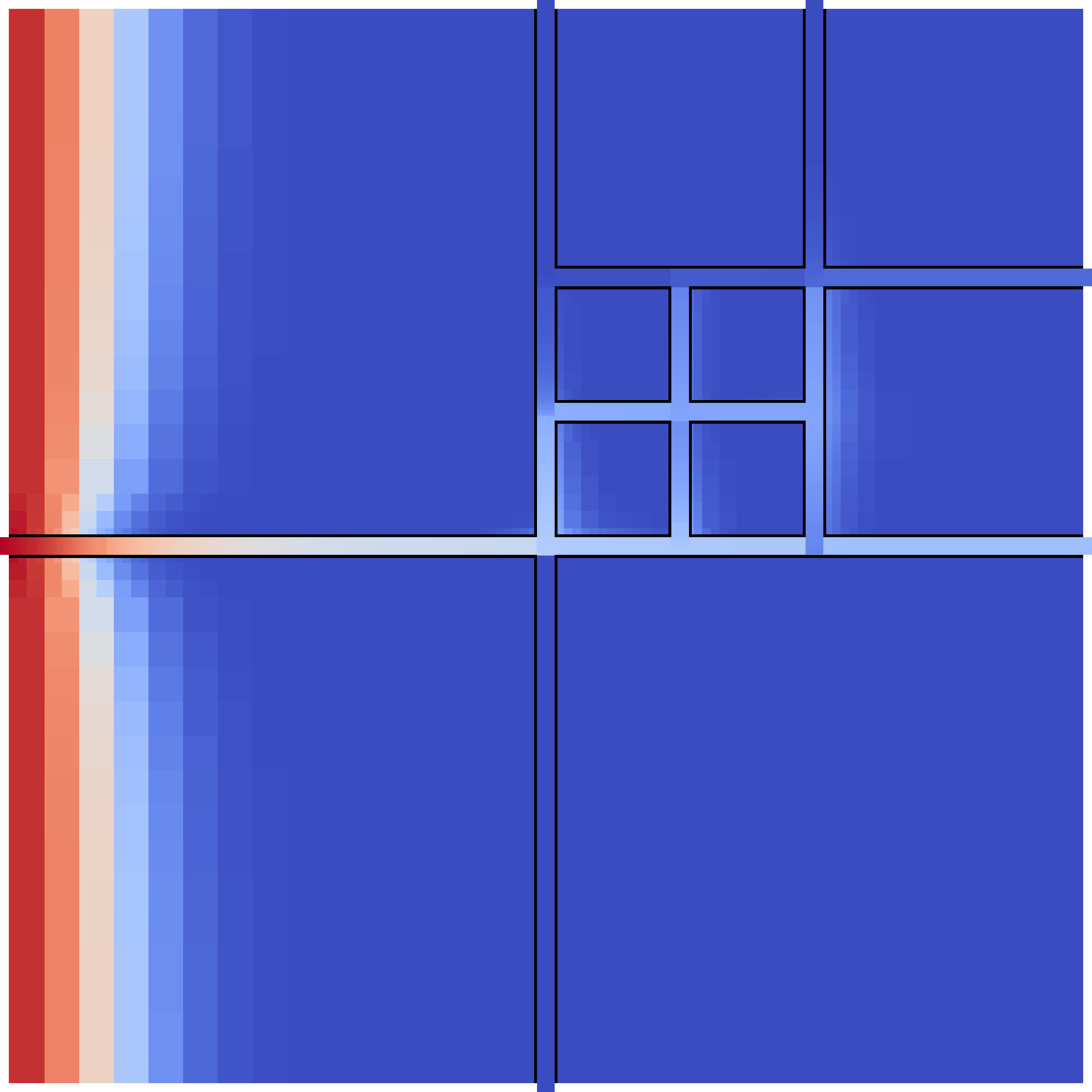}}
	\hspace{3mm}
	\subfloat[LR3, $t=0.5$.]
	{\includegraphics[width=0.36\textwidth]{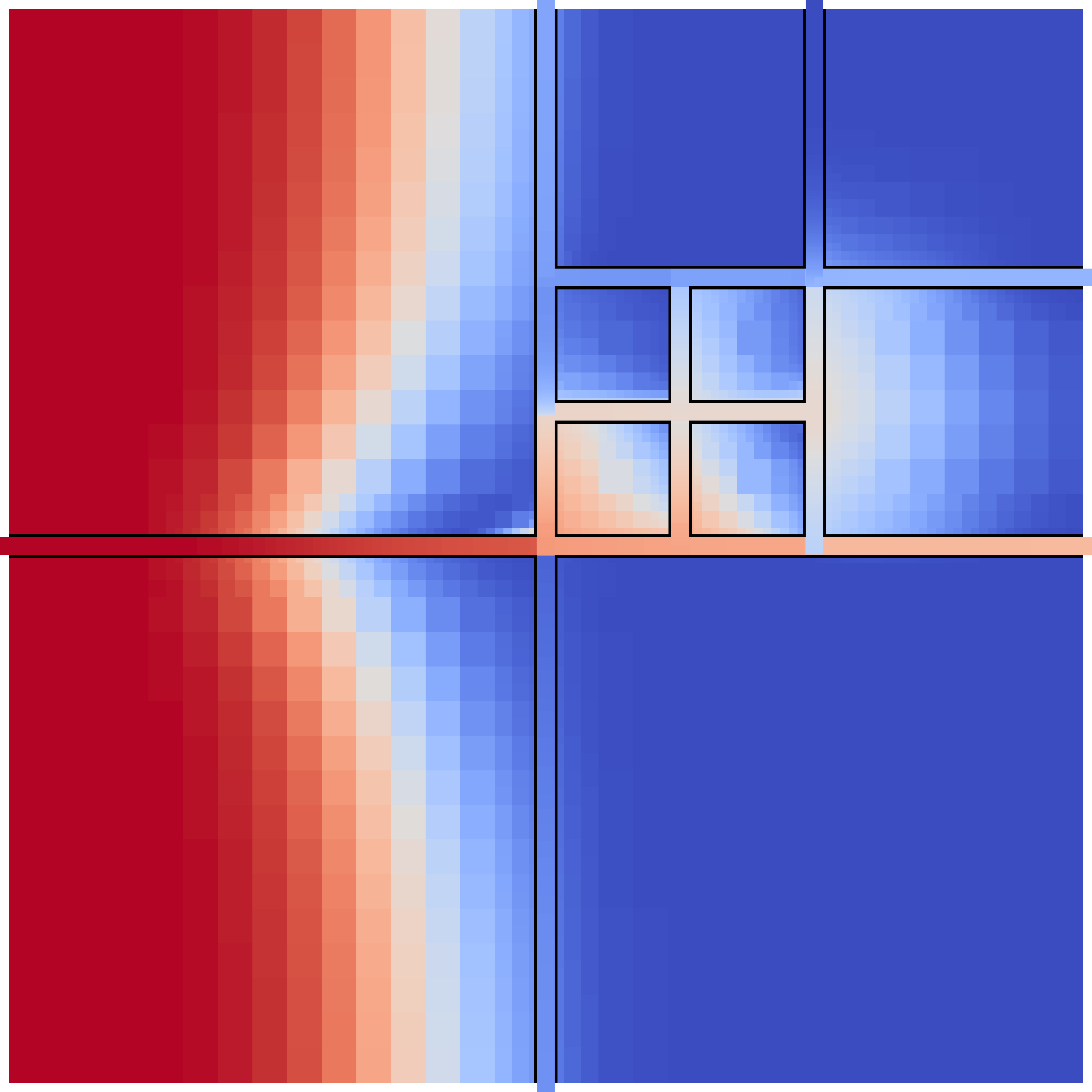}}
	\caption{Benchmark 1 (coupled problem). Concentration solution. The fractures are visualized with a fixed exaggerated width.  Number of degrees of freedom, $\ndof$, is 182674 for the reference solution (continuum model), 19321 for UMR139, and 16252 for LR3.}
	\label{fig:benchmark_conc}
\end{figure}

\begin{figure}[tbp]
	\centering
	\subfloat[Line $y=0.5$.]{
	\includegraphics[width=0.49\textwidth]{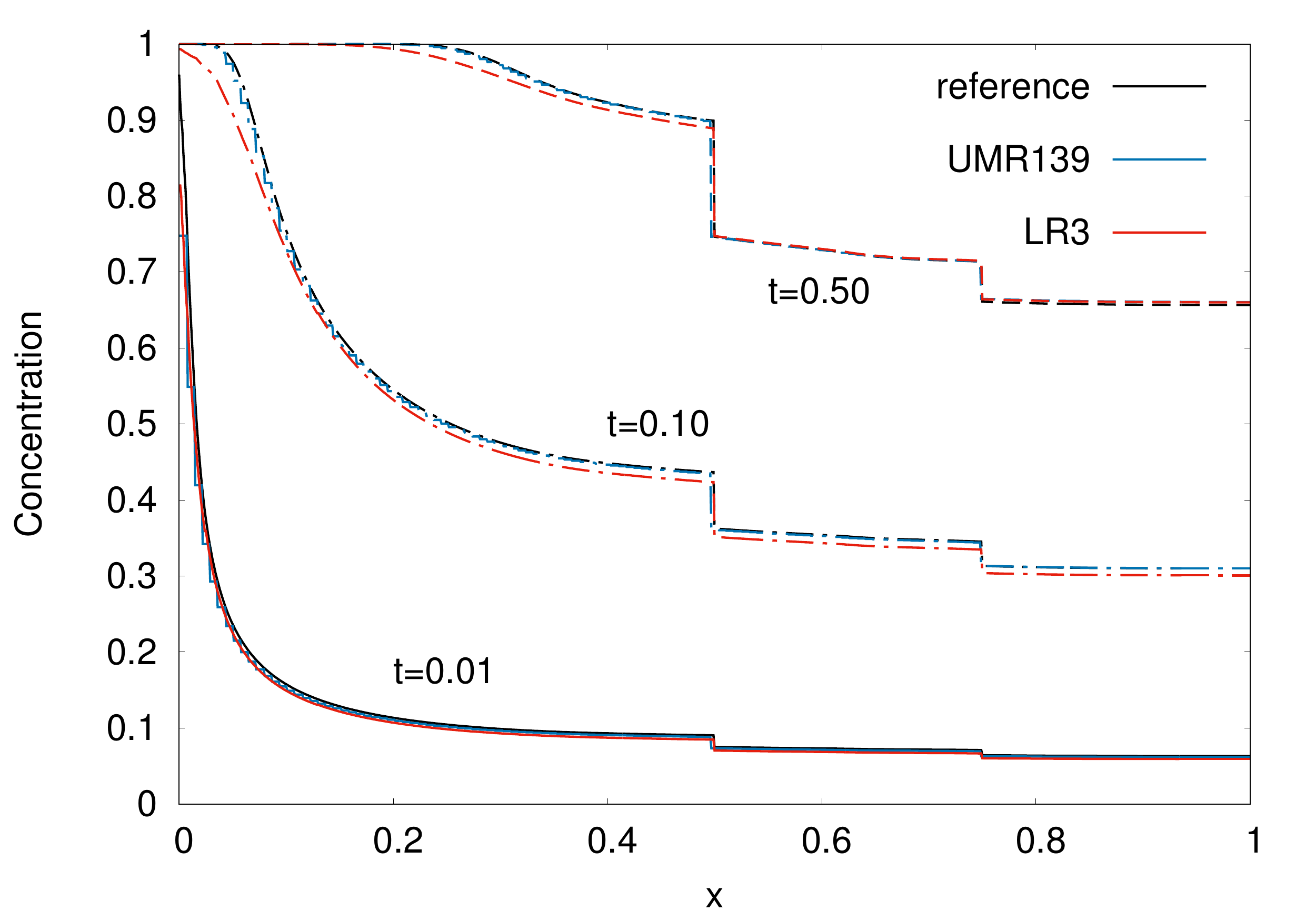}}
	\subfloat[Line $y=0.75$.]{
	\includegraphics[width=0.49\textwidth]{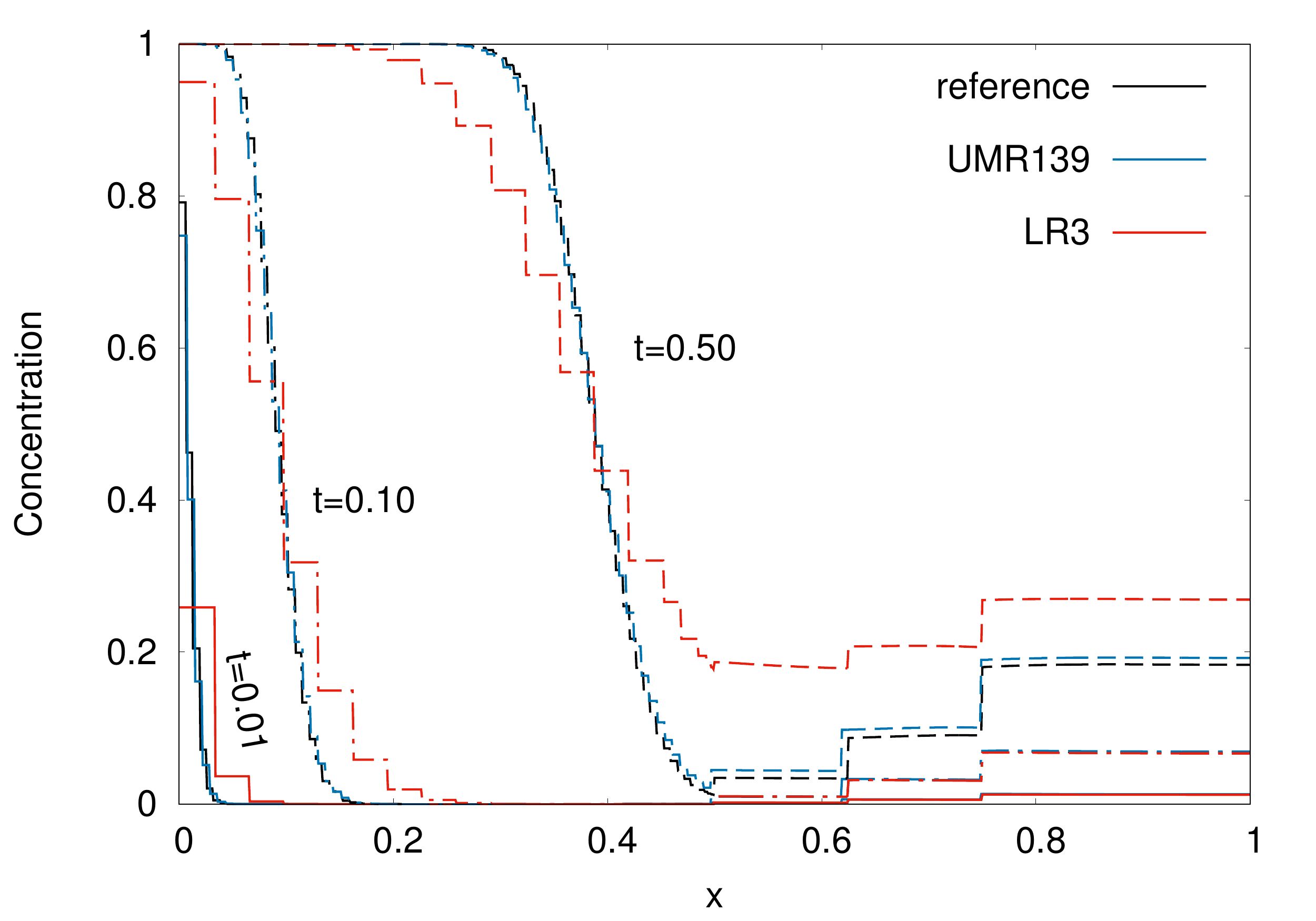}}
	\caption{Benchmark 1 (coupled problem). Concentration solution along two lines at different times.}
	\label{fig:benchmark_conc_fractures}
\end{figure}

\begin{figure}[tbp]
	\centering
	\includegraphics[width=0.49\textwidth]{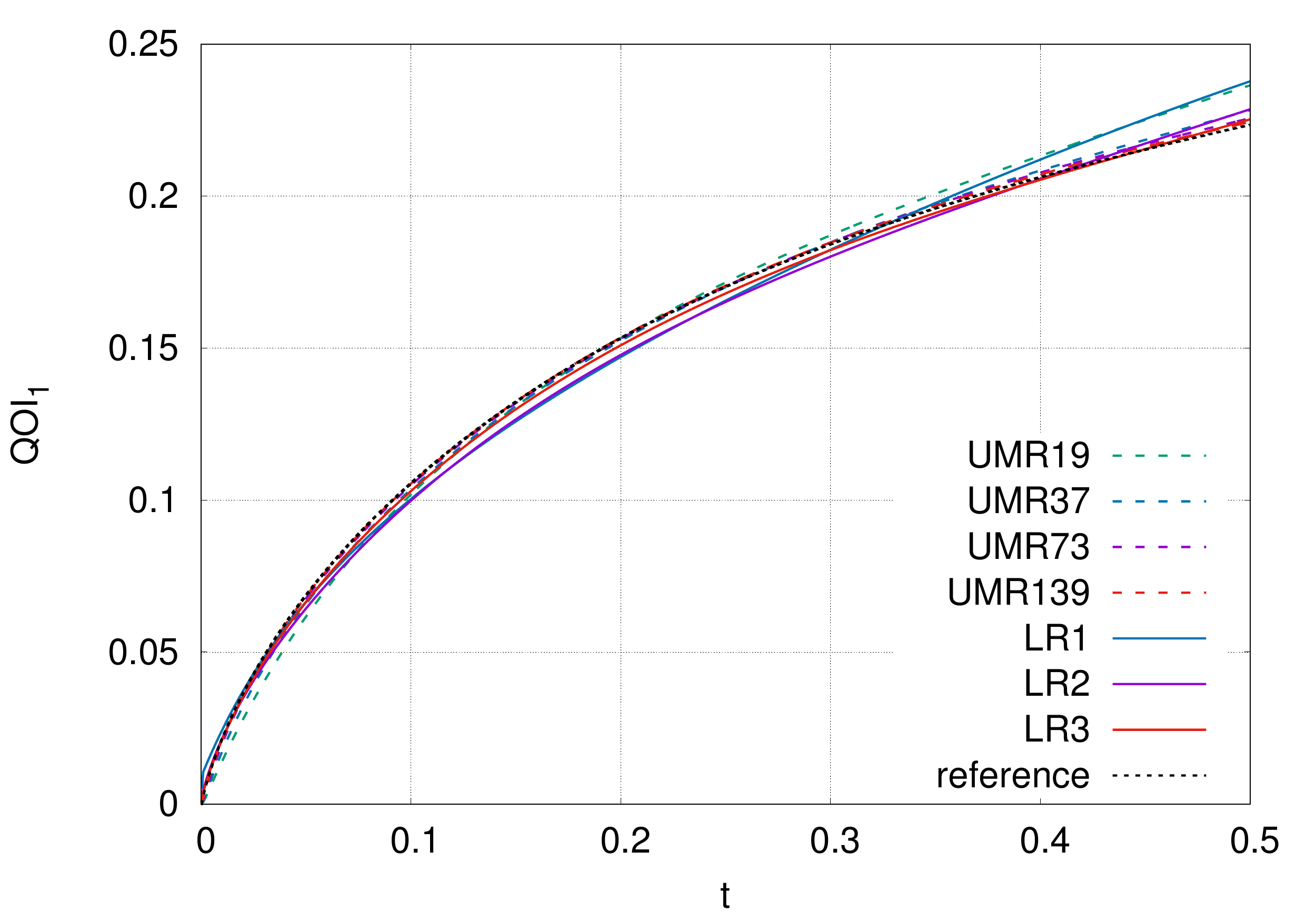}
	\includegraphics[width=0.49\textwidth]{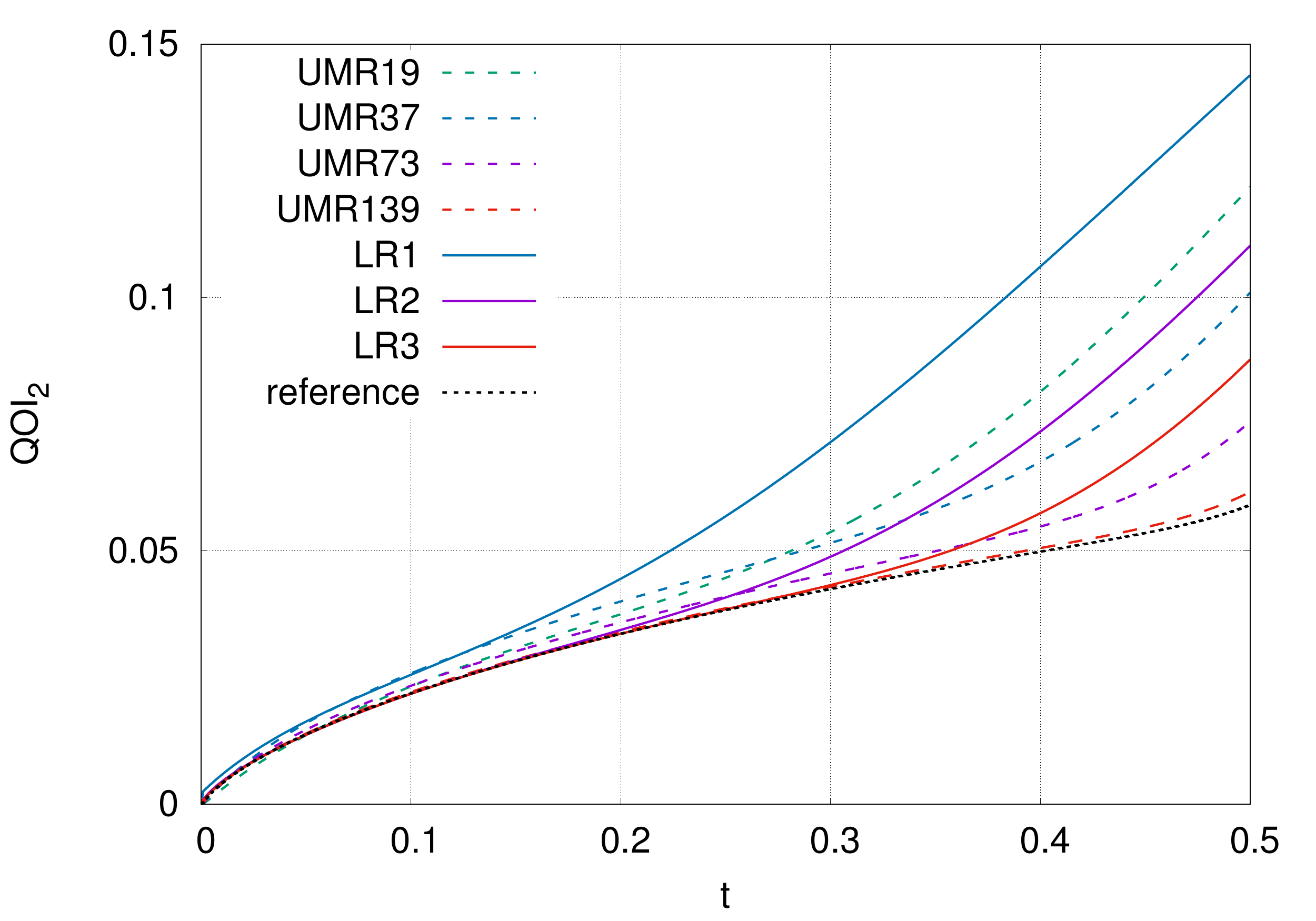}
	\caption{Benchmark 1 (coupled problem). Quantity of interest, $\qoi_i,\,i=1,2$, as functions of time. }
	\label{fig:benchmark_qoi}
\end{figure}

First, we solve the coupled problem on the four uniform meshes with $\Delta t=\powten{1.0}{-3}$, $\powten{5.0}{-4}$, $\powten{2.5}{-4}$, and $\powten{1.25}{-4}$, respectively, and on the three LR meshes with $\Delta t=\powten{5.0}{-4}$, $\powten{2.5}{-4}$, and $\powten{1.25}{-4}$, respectively. The concentration solution on the finest meshes are displayed in Figs.~\ref{fig:benchmark_conc} and \ref{fig:benchmark_conc_fractures}, while $\qoi$ is plotted against time for all meshes in Fig.~\ref{fig:benchmark_qoi}. 

We observe similar solutions for all meshes. At early times and in fracture 1 ($y=0.5$), we have the best results on the LR meshes. However, the LR meshes are relatively coarse in the matrix. This causes large numerical diffusion and with time the concentration front in the matrix reaches the first vertical fracture ($x=0.5$). This explains why the solution on the LR meshes becomes inaccurate at large times. A better meshing for the coupled problem would be to refine both close to the fractures and in the left half of $\Omega$.

\begin{figure}[tbp]
	\centering
	\includegraphics[width=0.49\textwidth]{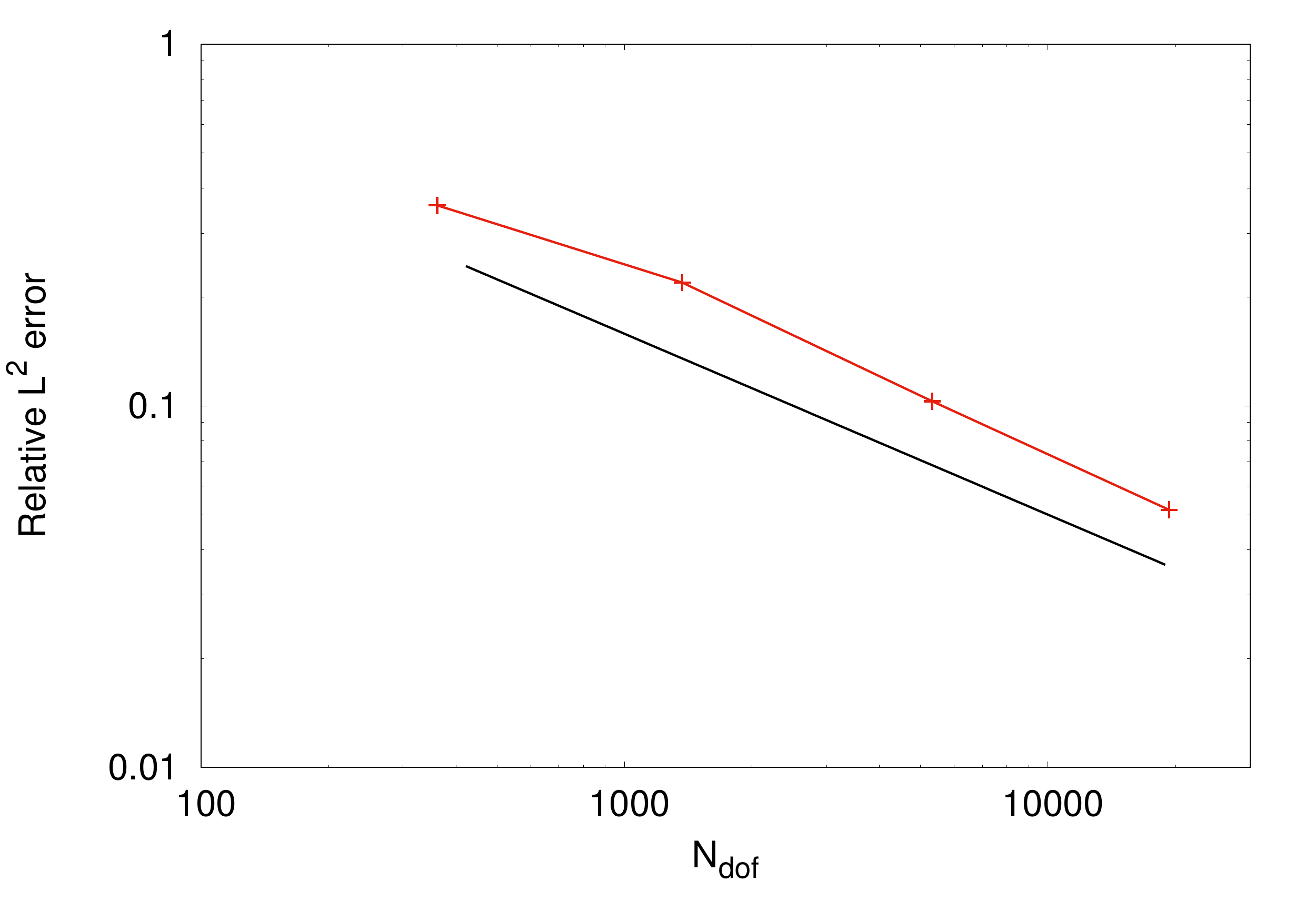}
	\caption{Benchmark 1 (coupled problem). Relative concentration error in fractures, $\frac{\Vert c_{\text{ref}}-c_h\Vert_{L^2(\Gamma)}}{\Vert c_{\text{ref}}\Vert_{L^2(\Gamma)}}$, at $t=T$ against $\ndof$. The reference solution is compared against the solution on the UMR meshes with $\Delta t=\powten{1}{-4}$. The black reference line has slope $-0.5$.}
	\label{fig:benchmark_conc_conv}
\end{figure}

At last, we ran a series of simulations on the UMR meshes with $\Delta t=10^{-4}$ and compared to the reference solution by a $L^2$ norm over the fractures, see Fig.~\ref{fig:benchmark_conc_conv}. We get a higher convergence than what is expected from the lowest order FV method.

%\clearpage
\subsubsection{Benchmark 4: A realistic case}

We now couple benchmark 4 to the transport problem, where $c_0=0$ and an inflow concentration $c_\tB=1$ is set on the left boundary face. Such problem was also considered in \cite{fumagalli2017dve}, but with different permeabilities. We set the end time for the simulations to $T=100~$years\footnote{$1~\textup{year}=365~\textup{days}$.}.

\begin{figure}[tpb]
	\centering
	%\subfloat[$\mathcal{M}_2^{2,\textup{r}}$, $t=5$~years.]
	%{\includegraphics[width=0.32\textwidth]{fig/realistic_case/conc_benchmark4_lr2-2r_t5}}
	%\hspace{1mm}
	\subfloat[$\mathcal{M}_2^{2,\textup{r}}$, $t=15$~years.]
	{\includegraphics[width=0.43\textwidth]{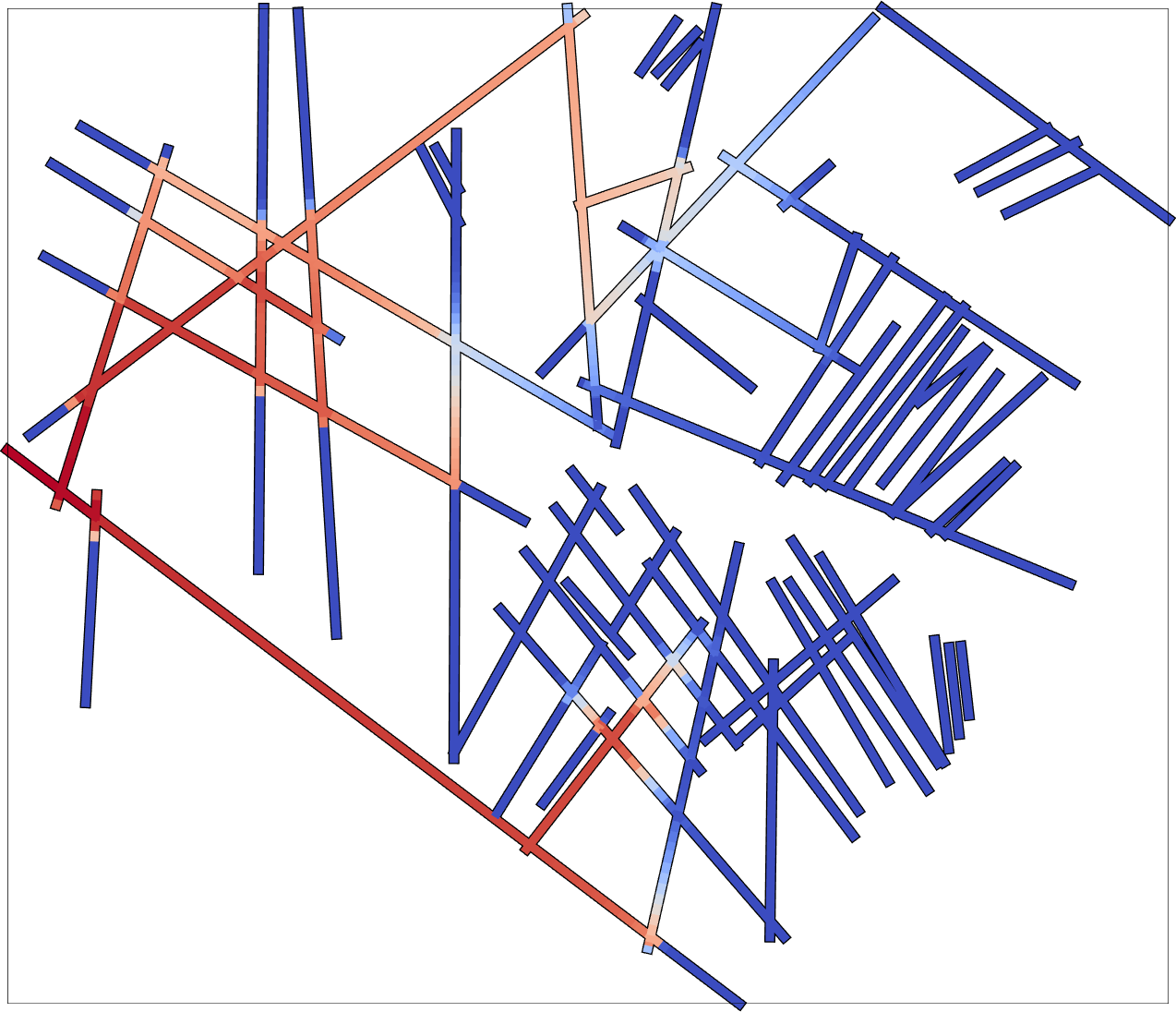}}
	\hspace{3mm}
	\subfloat[$\mathcal{M}_2^{2,\textup{r}}$, $t=100$~years.]
	{\includegraphics[width=0.43\textwidth]{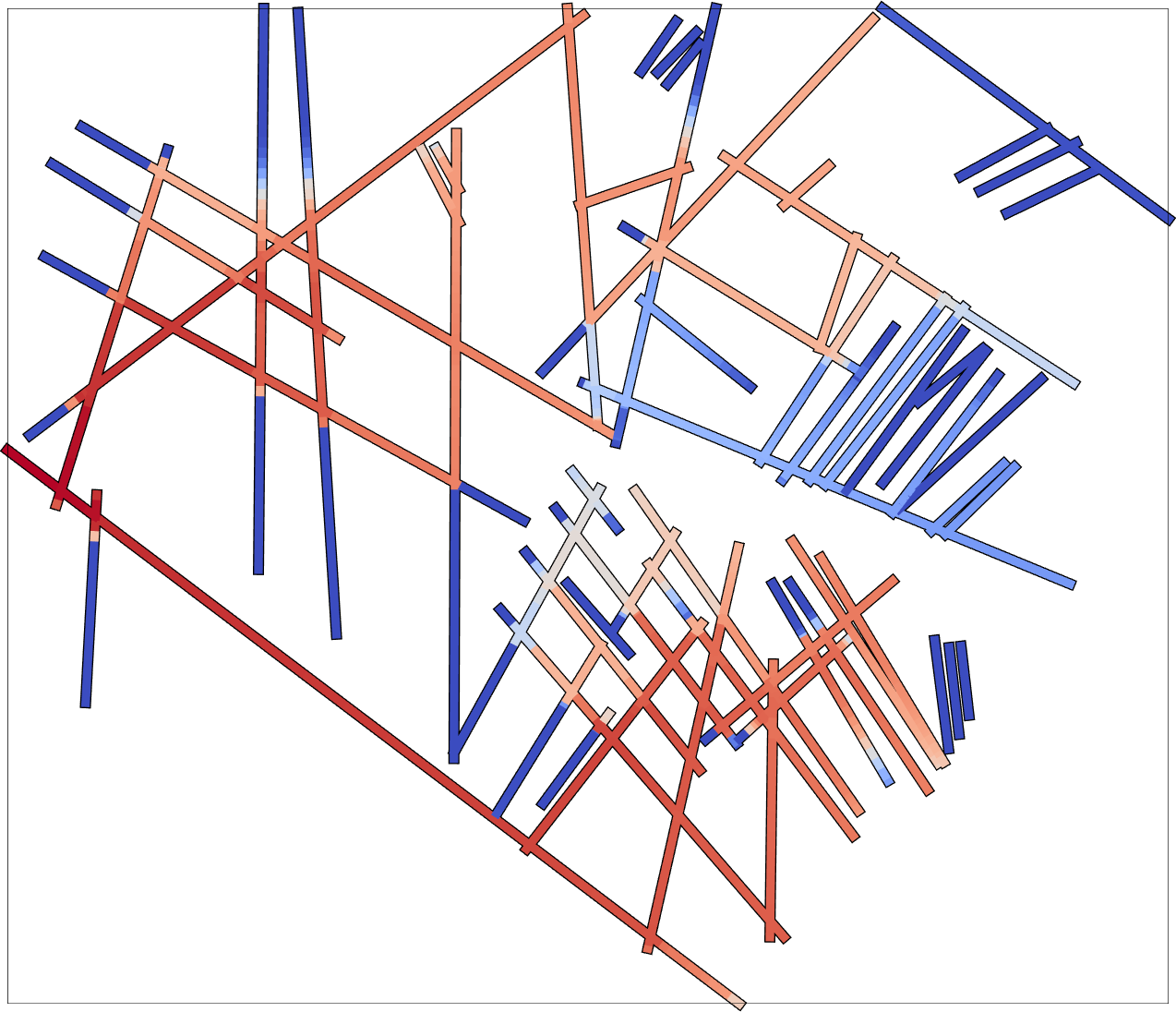}} \\
	%\subfloat[$\mathcal{M}_2^{4,\textup{r}}$, $t=5$~years.]
	%{\includegraphics[width=0.32\textwidth]{fig/realistic_case/conc_benchmark4_lr2-4r_t5}}
	%\hspace{1mm}
	\subfloat[$\mathcal{M}_2^{4,\textup{r}}$, $t=15$~years.]
	{\includegraphics[width=0.43\textwidth]{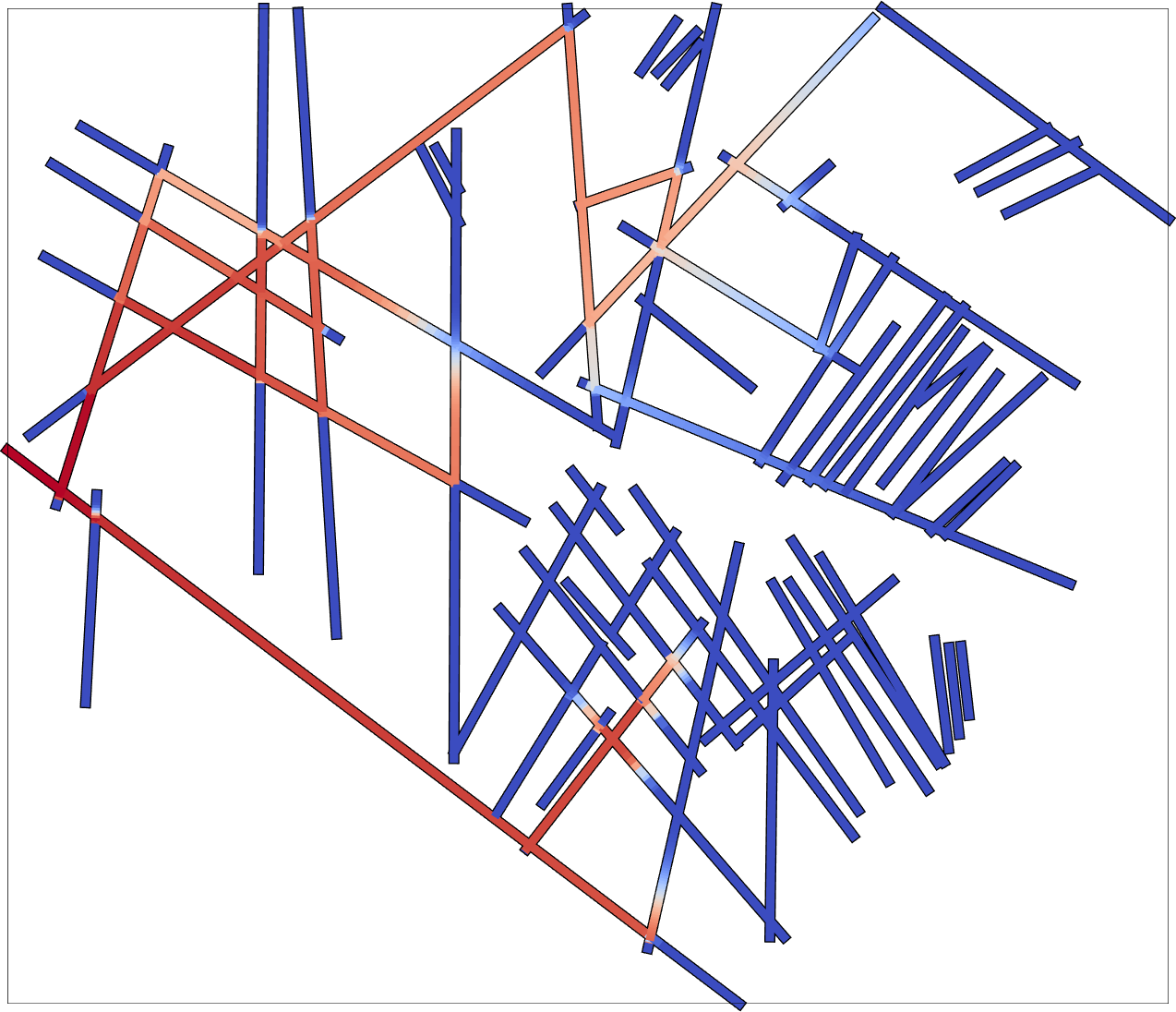}}
	\hspace{3mm}
	\subfloat[$\mathcal{M}_2^{4,\textup{r}}$, $t=100$~years.]
	{\includegraphics[width=0.43\textwidth]{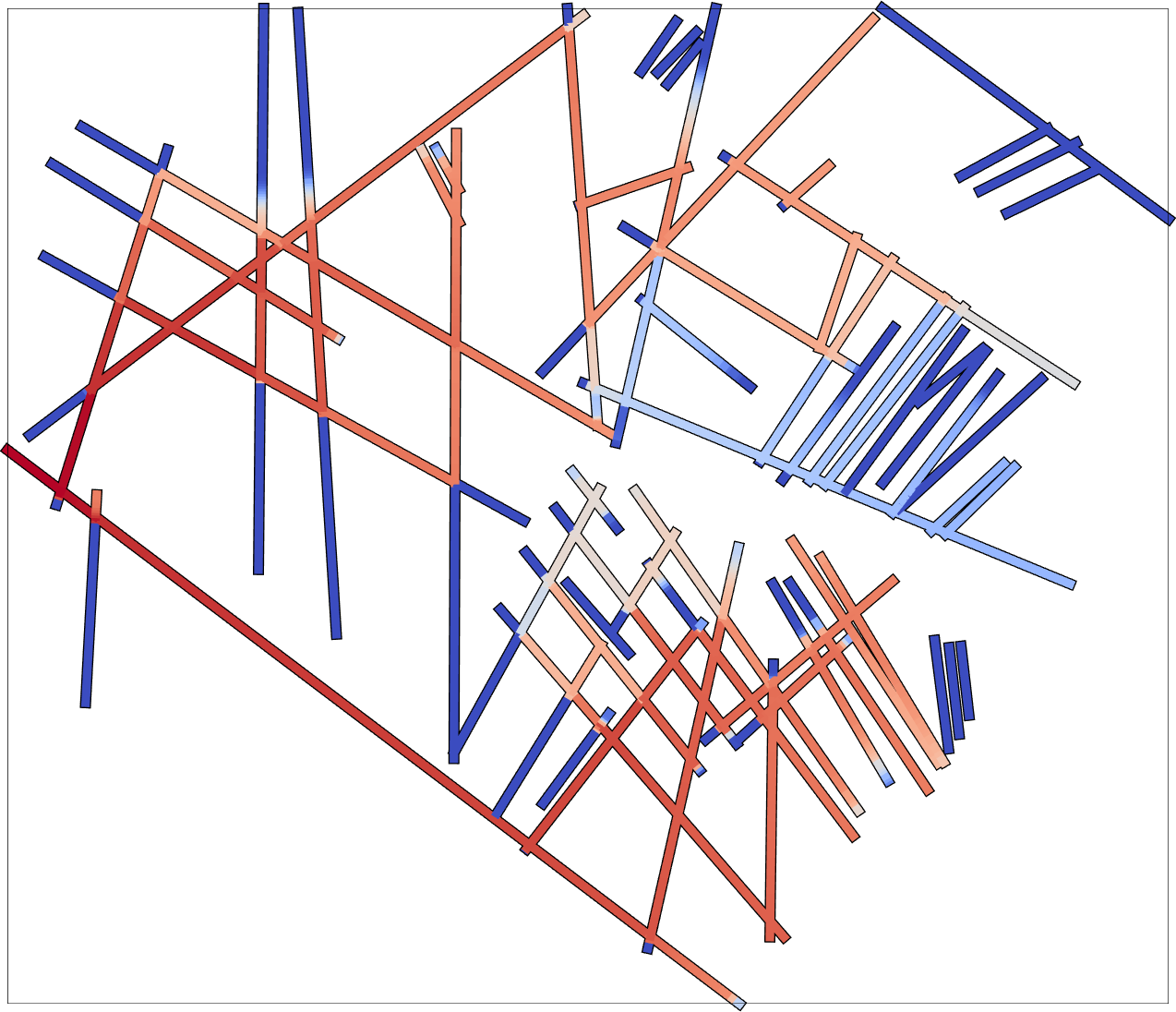}} \\
	%\subfloat[$\mathcal{M}_2^{4}$, $t=5$~years.]
	%{\includegraphics[width=0.32\textwidth]{fig/realistic_case/conc_benchmark4_lr2-4_t5}}
	%\hspace{1mm}
	\subfloat[$\mathcal{M}_2^{4}$, $t=15$~years.]
	{\includegraphics[width=0.43\textwidth]{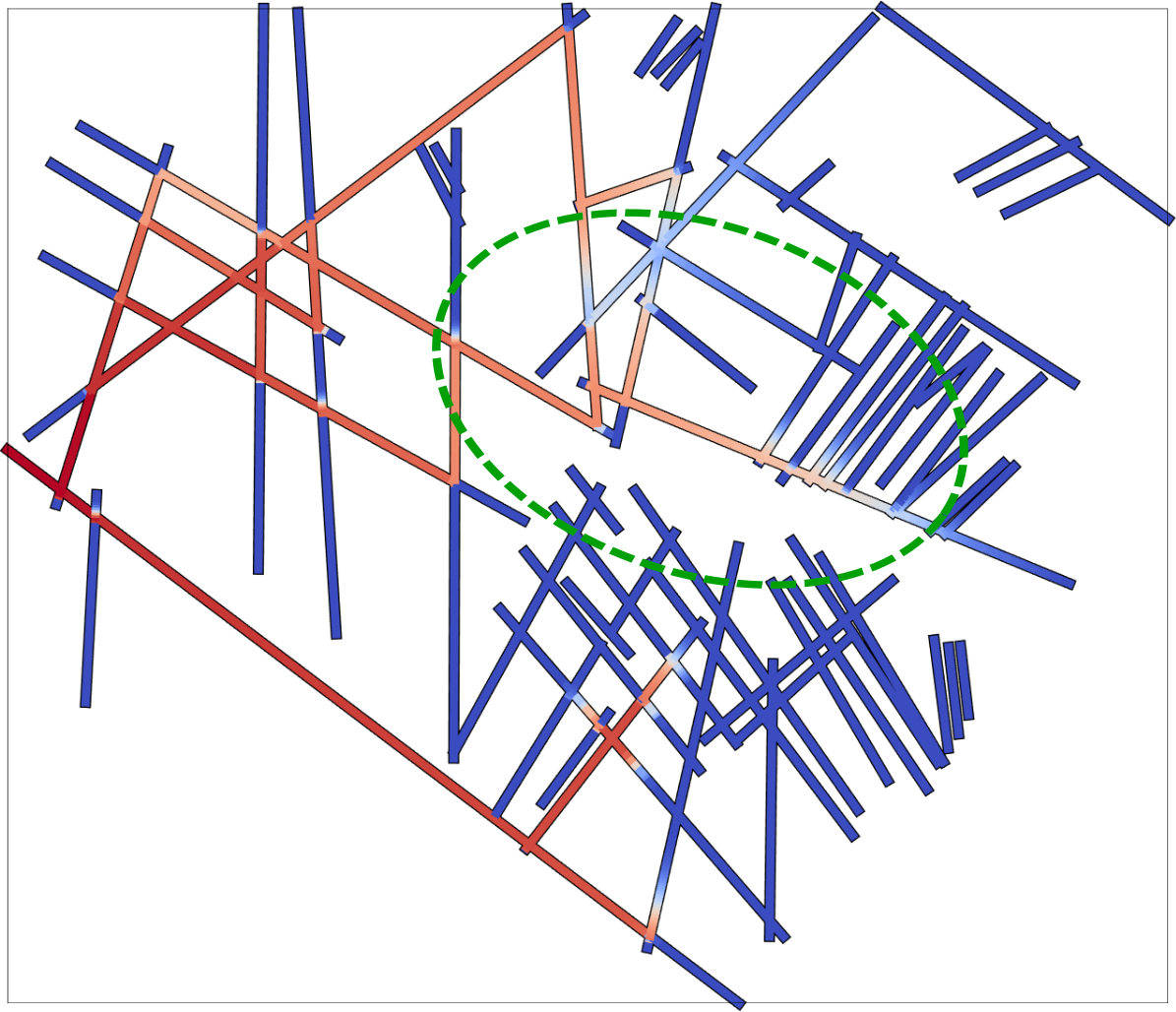}}
	\hspace{3mm}
	\subfloat[$\mathcal{M}_2^{4}$, $t=100$~years.]
	{\includegraphics[width=0.43\textwidth]{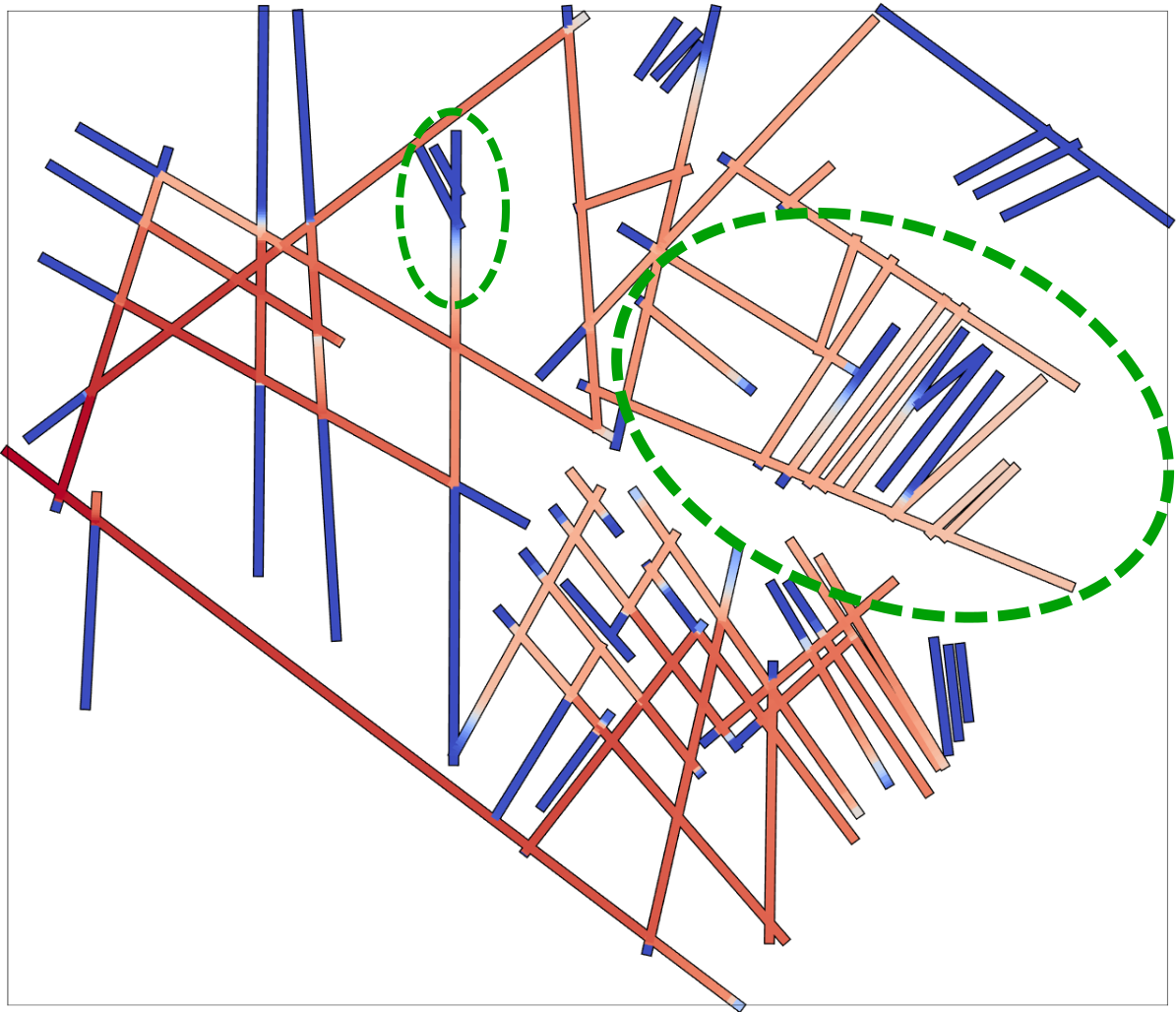}} \\
	\caption{Benchmark 4 (coupled problem). Concentration solutions in fractures. The fractures are displayed with a fixed (exaggerated) width. Regions where the results on the non-resolved mesh $\mathcal{M}_2^{4}$ deviates substantially from the other two meshes are marked with green dashed ellipses.}
	\label{fig:benchmark4_conc}
\end{figure}

\begin{figure}[tbp]
	\centering
	\includegraphics[width=0.49\textwidth]{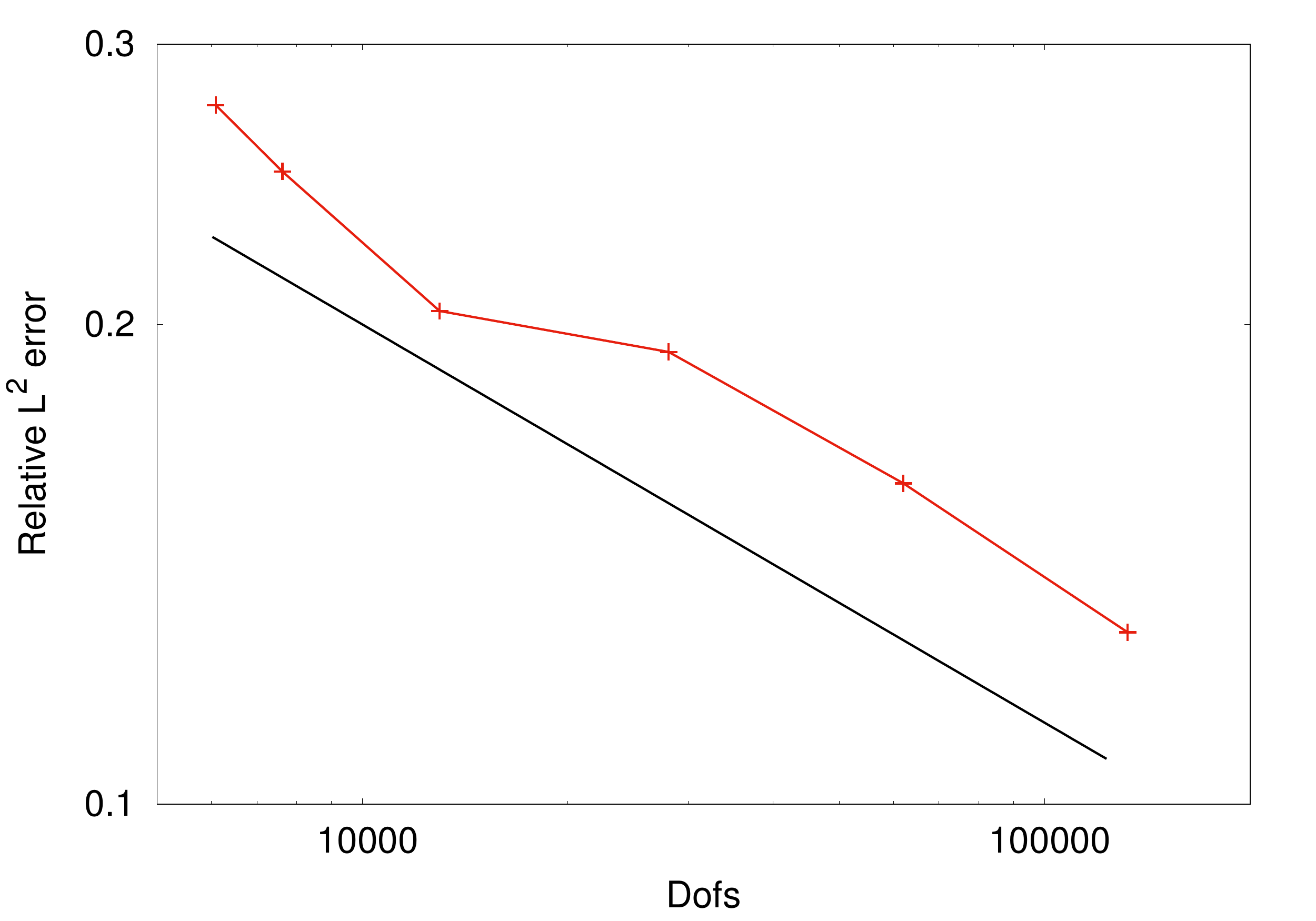}
	\caption{Benchmark 4 (coupled problem). Relative concentration error in fractures, $\frac{\Vert c_{\text{ref}}-c_h\Vert_{L^2(\Gamma)}}{\Vert c_{\text{ref}}\Vert_{L^2(\Gamma)}}$, at $t=T$ against $\ndof$. The solution on $\mathcal{M}_2^{8,\textup{r}}$ is used as reference solution, $c_{\text{ref}}$, and compared against the solutions on $\mathcal{M}_2^{j,\textup{r}}$ for $j=1,2,\ldots 6$. All simulaitons are run with $\Delta t=1~$day. The black reference line has slope $-0.25$.}
	\label{fig:benchmark4_conc_conv_space}
\end{figure}

First, we set $\Delta t=1~$hour (36500 time steps), and consider the meshes $\mathcal{M}_2^{j,\textup{r}}$. Concentration solutions in the fractures are displayed in Fig.~\ref{fig:benchmark4_conc}. Furthermore, convergence of the concentration in the fractures is demonstrated in Fig.~\ref{fig:benchmark4_conc_conv_space}. Due to the high fracture permeability, nearly all transport take place in the fractures, and hence we do not display the matrix solution. We observe that we get reasonable good results even for the coarsest mesh. Furthermore, the results for $\mathcal{M}_2^{4}$ clearly illustrates the importance of resolving the geometry, as the solution is far off in some of the fractures, even compared to $\mathcal{M}_2^{2,\textup{r}}$ which has much less DoFs. Convergence is rather slow due to the low order method.

\begin{figure}[tbp]
\centering
\includegraphics[width=0.49\textwidth]{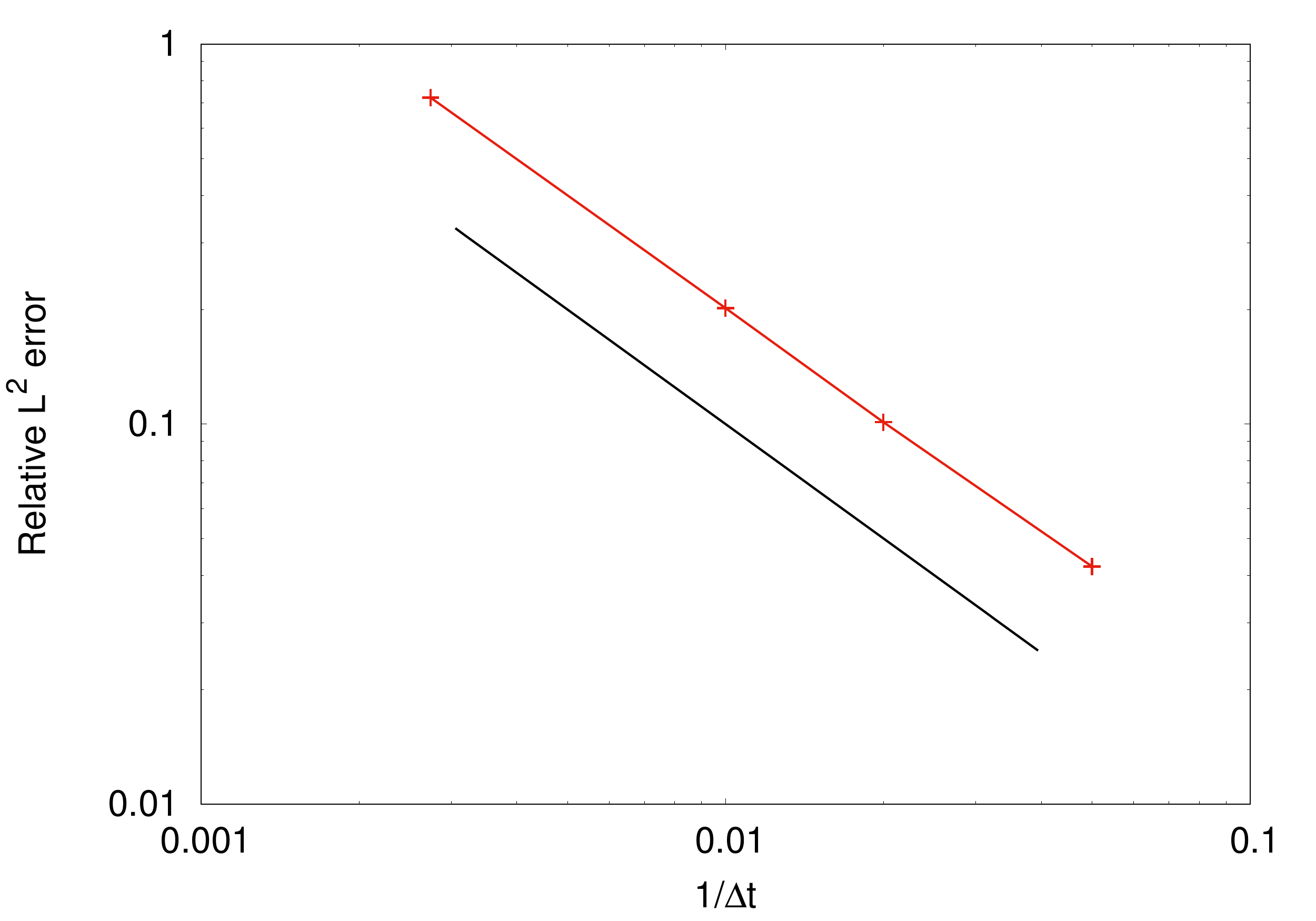}
	\caption{Benchmark 4 (coupled problem). Relative concentration error, $\frac{\Vert c_{\text{ref}}-c_h\Vert_{L^2(\Omega)}}{\Vert c_{\text{ref}}\Vert_{L^2(\Omega)}}$, at $t=T$ against $1/\Delta t$. All simulations are run on $\mathcal{M}_2^{2,\textup{r}}$. The solution with $\Delta t = 1~$hour is used as reference solution, $c_{\text{ref}}$, and compared to the solutions with $\Delta t = \{20,50,100,365\}~$days. The black reference line has slope $-1$.}
	\label{fig:benchmark4_conc_conv}
\end{figure}

Next, we ran a series of simulations on $\mathcal{M}_2^{2,\textup{r}}$ with $\Delta t = \{20,50,100,365\}$~days and compared to the solution with $\Delta t=1~\textup{hour}$ on the same mesh. The $L^2$ error over $\Omega$ at $t=T$ is plotted against $(\Delta t)^{-1}$ in Fig.~\ref{fig:benchmark4_conc_conv}, and we observe linear convergence in time as expected.

%\begin{figure}[tbp]
%	\centering
%	\subfloat[Fracture 12, $t=15~$years.]
%	{\includegraphics[width=0.33\textwidth]{fig/realistic_case/conc_frac12_t15}}
%	\subfloat[Fracture 28, $t=15~$years.]
%	{\includegraphics[width=0.33\textwidth]{fig/realistic_case/conc_frac28_t15}}
%	\subfloat[Fracture 31, $t=15~$years.]
%	{\includegraphics[width=0.33\textwidth]{fig/realistic_case/conc_frac31_t15}} \\
%	\subfloat[Fracture 12, $t=100~$years.]
%	{\includegraphics[width=0.33\textwidth]{fig/realistic_case/conc_frac12_t100}}
%	\subfloat[Fracture 28, $t=100~$years.]
%	{\includegraphics[width=0.33\textwidth]{fig/realistic_case/conc_frac28_t100}}
%	\subfloat[Fracture 31, $t=100~$years.]
%	{\includegraphics[width=0.33\textwidth]{fig/realistic_case/conc_frac31_t100}}
%	\caption{Benchmark 4 (coupled problem). Comparison of concentration solution in selected fractures at different times. Start and end points of displayed fractures: Fracture 12 $(297.2,237.6)\rightarrow(468.1,40.2)$; Fracture 28 $(44.8,528.6)\rightarrow(365.1,342.7)$; Fracture 31 $(347.2,374.0)\rightarrow(640.6,253.1)$.}
%	\label{fig:benchmark4_conc_frac}
%\end{figure}

%\begin{figure}[tbp]
%	\centering
%	\subfloat[Fracture 12, $t=15~$years.]
%	{\includegraphics[width=0.49\textwidth]{fig/realistic_case/image}}
%	\subfloat[Fracture 28, $t=15~$years.]
%	{\includegraphics[width=0.49\textwidth]{fig/realistic_case/image}}
%	\subfloat[Fracture 31, $t=15~$years.]
%	{\includegraphics[width=0.33\textwidth]{fig/realistic_case/image}}
%	\caption{Benchmark 4 (coupled problem).}
%\end{figure}

%\clearpage
\section{Concluding remarks}
\label{sec:conclusion}
This article addresses the \emph{numerical solution of a coupled flow and transport problem} in fractured porous media, where the \emph{fractures are modeled as lower-dimensional interfaces} embedded in the surrounding matrix. The proposed solution strategy includes three main steps: (1) Solving the flow problem with an embedded finite element method (EFEM) \cite{burman2017asf}; (2) Locally conservative flux approximation; and (3) Solving the transport problem with a non-standard lowest order finite volume (FV) method.

\emph{The main contribution of this work is that we couple EFEM with a numerical model for the transport problem}. EFEM allows for complex fracture geometry, where the fractures can cut the elements arbitrarily, and compared to other embedded discrete fracture-matrix models, there are no lower-dimensional elements along the fractures. The numerical model for the transport problem presented in this work, aims to be as flexible as EFEM with respect to meshing. This is resolved by a lowest order upwind FV method where the fracture solution is represented by elements cut by the fracture. The main novelty in our method is how we approximate the coupling term between the matrix and the fractures. 

Furthermore, this work includes the following contributions:
\begin{itemize}
	\itemsep0pt
	\item We apply EFEM to realistic benchmark problems presented in \cite{flemisch2018bfs}. Our results are in good agreement with the results in \cite{flemisch2018bfs} and in some cases EFEM is most accurate. Applying a priori local refinement based on an estimate in \cite{burman2017asf} gives especially good results.
	\item To ensure locally conservative fluxes, we adapt and apply the postprocessing method presented in \cite{odsaeter2016pon}.
	\item We demonstrate the abilities of our coupled solution strategy by numerical examples on a realistic case with a complex fracture network.
\end{itemize}

\emph{A main advantage of our method compared to other discrete fracture-matrix models, is the inherent simplicity}, both in terms of formulation, implementation and meshing. Despite the simplicity, the presented results are promising.

One direction for further work is to consider adaptive routines for the transport solver, perhaps looking at the space-time approach allowing for easy handling of local time refinement. Rigorous convergence proofs and error estimates for the transport solver are also desired. Moreover, one could pursue including more physics into the model, e.g., multi-phase flow. 
To reduce the numerical diffusion associated with the implicit Euler method, one could consider higher order time integration, e.g., \cite{matthai2010sos}.

\section*{Acknowledgements}
LHO is funded by VISTA (Grant No.\ 6355) --- a basic research program funded by Statoil,
conducted in close collaboration with The Norwegian Academy of Science and Letters. MGL was supported by The 
Swedish Foundation for Strategic Research Grant No.\ AM13-0029, and the Swedish Research Council Grants Nos.\ 2013-4708, 2017-03911. The authors appreciate the open access to the results given in the benchmark paper \cite{flemisch2018bfs}. Furthermore, the authors thank Karl Larsson, Department of Mathematics and Mathematical Statistics, Umeå University, for his help with visualization of the results in Figs.~\ref{fig:benchmark_conc} and \ref{fig:benchmark4_conc}.

\renewcommand\thesubsection{A.\arabic{subsection}}
\renewcommand{\theequation}{A.\arabic{equation}}

\section*{Appendix}

\subsection{Discontinuous Galerkin formulation}
\label{app_DG}
We formulate here a lowest order Discontinuous Galerkin method (DG(0)) with upwinding for the transport problem, Eq.~\eqref{eq:transport_problem}. This formulation is equivalent to the finite volume method derived in Section \ref{sec:methods_transport}. We include it here for convenience of the reader more familiar with DG methods. The DG formulation is also more suitable for deriving error estimates.

We start by multiplying Eq.~\eqref{eq:transport_matrix} by a test function $v$, integrate over an element $K\in\Kh^\tM$ and apply Greens formula to obtain
\begin{align}
\int_K\phi\dd{c}{t}v - \int_K c\bfu\cdot\nabla v + \int_{\pK}\bfu\cdot\bfn_K cv = \int_K f(c)v, \quad K\in\Kh^\tM.
\end{align}
Summing the boundary term over all $K\in\Kh^\tM$ we get
\begin{align}
\sum_{K\in\Kh^\tM} \int_{\pK}\bfu\cdot\bfn_K cv =
&- \sum_{F\in\Fh{\tI}^\tM} \int_F \bfu\cdot\bfn_Fc\jump{v} \nonumber \\
&+ \sum_{F\in\Fh{\tI}^{\tF\tM}} \int_F \bfu\cdot\bfn_{\tM\tF} cv \nonumber \\
&+ \sum_{F\in\Fh{\tout}^\tM}\int_F  \bfu\cdot\bfn_F cv 
+ \sum_{F\in\Fh{\tin}^\tM} \int_F \bfu\cdot\bfn_F c_\tB v.
\end{align}
We have denoted by $\bfn_{\tM\tF}$ the unit normal on $F\in\Fh{\tI}^{\tF\tM}$ pointing from the matrix domain towards the fracture.

Similarly, we multiply Eq.~\eqref{eq:transport_fracture} by a test function $v_\Gamma$, integrate over $K\cap\Gamma$, for $K\in\Kh^\tF$, and apply Greens formula to obtain
\begin{align}
\int_{K\cap\Gamma} w\phi_\Gamma \dd{c_\Gamma}{t} v_\Gamma
- \int_{K\cap\Gamma}c_\Gamma\bfu_\Gamma\cdot\nabla_\Gamma v_\Gamma
+ \int_{\partial (K\cap\Gamma)}\bfu_\Gamma\cdot\bfn_{K\cap\Gamma} c_\Gamma v_\Gamma
- \int_{K\cap\Gamma} \jump{\bfu\cdot\bfn c^*}v_\Gamma
= \int_{K\cap\Gamma} f_\Gamma(c_\Gamma) v_\Gamma.
\end{align}
Again, summing the boundary term over all $K\in\Kh^\tF$ we get
\begin{align}
\sum_{K\in\Kh^\tF}\int_{\partial (K\cap\Gamma)}\bfu_\Gamma\cdot\bfn_{K\cap\Gamma} c_\Gamma v_\Gamma
= &-\sum_{F\in\Fh{\tI}^\tF} \lp\bfu_\Gamma\cdot\bft_{\Gamma,F} c_\Gamma\jump{v_\Gamma}\rp\vert_{F\cap\Gamma} \nonumber \\
&+ \sum_{F\in\Fh{\tin}^\tF} \lp\bfu_\Gamma \cdot \bft_{\Gamma,F} c_{\Gamma,\tB}v_\Gamma\rp\vert_{F\cap\Gamma} \nonumber \\
&+ \sum_{F\in\Fh{\tout}^\tF} \lp\bfu_\Gamma \cdot \bft_{\Gamma,F} c_\Gamma v_\Gamma\rp\vert_{F\cap\Gamma},
\end{align}
where $\bft_{\Gamma,F}$ is the unit tangent to $\Gamma$ oriented in the same direction as $\bfn_F$.

To formulate the DG method, we replace $c$ and $c_\Gamma$ by $c_h\in Q_0(\Kh)$ and $\varphi$ and $\varphi_\Gamma$ by $\varphi_h\in Q_0(\Kh)$. Furthermore, we approximate the coupling term as follows,
\begin{align}
\int_{K\cap\Gamma} \jump{\bfu\cdot\bfn c_h^*}v_h \approx \sum_{\tilde{K}\in\Kh^\tM} \int_{\partial\tilde{K}\cap\pK} (\bfn_F\cdot\bfn_K)\bfu\cdot\bfn_F c_h v_h.
\end{align}
This is an equivalent approximation as Eq.~\eqref{eq:coupling_approx}.
Summing this term over all $K\in\Kh^\tF$ we get
\begin{align}
\sum_{K\in\Kh^\tF} \int_{K\cap\Gamma} \jump{\bfu\cdot\bfn c_h^*}v_h \approx
\sum_{F\in\Fh{\tI}^{\tF\tM}} \int_F \bfu\cdot\bfn_{\tM\tF} c_hv_h.
\end{align}

Adding everything up, we end up the following formulation. Find $c_h\in Q_0(\Kh)$ such that
\begin{align}
b(c_h,v_h) = k(v_h), \qquad \forall v_h\in Q_0(\Kh),
\end{align}
where
\begin{align}
b(c,v) =& 
\sum_{K\in\Kh^\tM} \int_K \lp 
\phi\dd{c}{t} cv 
+ \check{q}cv
\rp
+ \sum_{K\in\Kh^\tF} \int_{K\cap\Gamma} \lp
w\phi_\Gamma \dd{c}{t} v
+ \check{q}_\Gamma cv
\rp \nonumber\\
&- \sum_{F\in\lp\Fh{\tI}\setminus\Fh{\tI}^\tF\rp} \int_F \bfu\cdot\bfn_Fc\jump{v}
- \sum_{F\in\Fh{\tI}^\tF} \lp\bfu_\Gamma\cdot\bft_{\Gamma,F}c\jump{v}\rp\vert_{F\cap\Gamma} \nonumber \\
&+ \sum_{F\in\Fh{\tout}^\tM} \int_F \bfu\cdot\bfn_F cv 
+ \sum_{F\in\Fh{\tout}^\tF} \lp\bfu_\Gamma \cdot \bft_{\Gamma,F} c v\rp\vert_{F\cap\Gamma}
\end{align}
and
\begin{align}
k(v) =& 
\sum_{K\in\Kh^\tM} \int_K \hat{q}c_\textup{w} v
+ \sum_{K\in\Kh^\tF} \int_{K\cap\Gamma} \hat{q}_\Gamma c_\textup{w} v \nonumber \\
&+ \sum_{F\in\Fh{\tin}^\tM} \int_F \bfu\cdot\bfn_F c_\tB v
+ \sum_{F\in\Fh{\tin}^\tF} \lp\bfu_\Gamma \cdot \bft_{\Gamma,F} c_{\Gamma,\tB} v\rp\vert_{F\cap\Gamma}.
\end{align}
We have used the definition of $f(c)$ and $f_\Gamma(c_\Gamma)$, see Eq.~\eqref{eq:rhs_def}.
For the concentration on faces $F\in\Fh{}$, we use the upwind approximation
\begin{align}
\bfu\cdot\bfn_Fc_h =
\begin{cases}
\bfu\cdot\bfn_F (c_h)_-, & \text{if } \bfu\cdot\bfn_F \ge 0, \\
\bfu\cdot\bfn_F (c_h)_+, & \text{if } \bfu\cdot\bfn_F < 0.
\end{cases}
\end{align}
Applying the velocity model described by Eqs.~\eqref{eq:velocity}--\eqref{eq:velocity_pp}, we may simplify the forms $b(\cdot,\cdot)$ and $k(\cdot)$ for the coupled problem by a similar approach as described by Eq.~\eqref{eq:FV_simplification}, so that
\begin{align}
b(c,v) =& 
\sum_{K\in\Kh^\tM} \int_K \lp 
\phi\dd{c}{t} cv 
+ \check{q}cv
\rp
+ \sum_{K\in\Kh^\tF} \int_{K\cap\Gamma} \lp
w\phi_\Gamma \dd{c}{t} v
+ \check{q}_\Gamma cv
\rp \nonumber \\
&- \sum_{F\in\Fh{\tI}} \int_F V_h c\jump{v} 
+ \sum_{F\in\Fh{\tout}} \int_F V_h cv
\end{align}
and
\begin{align}
k(v) =& 
\sum_{K\in\Kh^\tM} \int_K \hat{q}c_\textup{w} v
+ \sum_{K\in\Kh^\tF} \int_{K\cap\Gamma} \hat{q}_\Gamma c_\textup{w} v 
+ \sum_{F\in\Fh{\tin}^\tM} \int_F V_h c_\tB v
+ \sum_{F\in\Fh{\tin}^\tF} \int_F V_h c_{\Gamma,\tB} v.
\end{align}
If we apply implicit Euler as time integrator, we end up with a scheme that is equivalent to the FV-IE scheme in Eq.~\eqref{eq:FV-IE_pp}.

\subsection{Interpretation of concentration solution}
\label{app_interpretation}
The numerical method given by Eq.~\eqref{eq:FV-IE_pp} defines a constant solution on $K\in\Kh^\tF$. However, $K$ originally contains both a fractured domain represented as a lower-dimensional interface, $K\cap\Gamma$, and a matrix domain $K\setminus\Gamma$. Let $\{K_j\}_{j\in\{1,\ldots,N_K\}}$ be a partition of $K$ into $N_K$ subelements defined by $\Gamma$. For a single fracture cutting $K$, we have $n_K=2$, but for intersecting and bifurcating fractures in $K$, $n_K$ can be larger. For a terminating fracture, we have $n_K=1$.
We interpret the solution on $K$ by assigning a constant value to each subelement $K_j$ and one value on the fracture intersection $K\cap\Gamma$. Let $\mathring{c}_h$ denote the interpreted solution. Then we define
\begin{align}
\mathring{c}_h\vert_{K} &= c_h\vert_K, \quad \forall K\in\Kh^\tM, \\
\mathring{c}_h\vert_{K\cap\Gamma} &= c_h\vert_K, \quad \forall K\in\Kh^\tF.
\end{align}

It is not as easy to define the interpreted solution on the fracture subelements, but the following algorithm can be used. Denote by $\Kh^{\tM,i}$ the (non-empty) subsets of $\Kh^\tM$ with elements contained in $\Omega_i$, i.e.,
\begin{align}
\Kh^{\tM,i} = \{ K\in\Kh^\tM : K \subset \Omega_i \}.
\end{align} 
Furthermore, let $\Kh^{\tF,i}$ be the set of all fracture subelements contained in $\Omega_i$, i.e.,
\begin{align}
\Kh^{\tF,i} = \{ K_j\in K : K_j\subset \Omega_i, \, K\in\Kh^\tF \}.
\end{align}
At last, denote by $\mathcal{N}_K(F)$ the neighbor of $K$ that shares face $F$.

Algorithm \ref{alg:interpret_conc} recursively assigns values to each element $\tilde{K}\in\Kh^{\tF,i}$. Once a value is assigned to $\tilde{K}$, we move $\tilde{K}$ from $\Kh^{\tF,i}$ to $\Kh^{\tM,i}$. In this way we mark $\tilde{K}$ as assigned and also allow for its value to be further assigned to another neighbor in the next cycle.
The algorithm works under the mild assumption that there is at least one matrix element for each subdomain, i.e., for all $\Omega_i,\,i=1,2,\ldots,n_d$, there is an $K\in\Kh^{\tM}$ such that $K\subset\Omega_i$. Moreover, we remark that the algorithm is sensitive to the ordering of elements. 
We refer to Fig.~\ref{fig:fracture_conc_interpretation} for an illustrative example of how this algorithm works. Finally, we emphasize that this algorithm is purely for the interpretation of the results, and not part of the numerical method.

%\begin{algorithm}
%\caption{Assigning values to fracture subelements}
%\label{alg:interpret_conc}
%\begin{algorithmic}
%	\FOR{$i=1,\ldots,n_d$} 
%		\WHILE{$\Kh^{\tM,i}\neq\emptyset$} 
%			\FOR{$K\in\Kh^{\tM,i}$}
%				\FOR{$F\in\pK$}
%					\STATE{$\tilde{K} = \mathcal{N}_K(F)$}
%					\IF{$\tilde{K}\in \Kh^{\tF,i}$}
%						\STATE{$\mathring{c}_h\vert_{\tilde{K}} = c_h\vert_K$}
%						\STATE{$\Kh^{\tF,i} = \Kh^{\tF,i} \setminus \{\tilde{K}\}$}
%						\STATE{$\Kh^{\tM,i} = \Kh^{\tM,i} \cup \{\tilde{K}\}$}
%					\ENDIF
%				\ENDFOR
%			\ENDFOR
%		\ENDWHILE
%	\ENDFOR
%\end{algorithmic}
%\end{algorithm}

\begin{algorithm}
	\SetAlgoNoLine
	\caption{Assigning values to fracture subelements}
	\label{alg:interpret_conc}
	\For{$i=1,\ldots,n_d$}{ 
		\While{$\Kh^{\tM,i}\neq\emptyset$}{
			\For{$K\in\Kh^{\tM,i}$}{
				\For{$F\in\pK$}{
					$\tilde{K} = \mathcal{N}_K(F)$\\
					\If{$\tilde{K}\in \Kh^{\tF,i}$}{
						$\mathring{c}_h\vert_{\tilde{K}} = c_h\vert_K$\\
						$\Kh^{\tF,i} = \Kh^{\tF,i} \setminus \{\tilde{K}\}$\\
						$\Kh^{\tM,i} = \Kh^{\tM,i} \cup \{\tilde{K}\}$\\
					}
				}
			}
		}
	}
\end{algorithm}

\bibliographystyle{acm}
%\section*{References} % This is not needed on arXiv
\bibliography{references}

\end{document}